\documentclass{tac}
\usepackage[usenames,dvipsnames,svgnames,table]{xcolor}
\usepackage{natbib} 
\usepackage{amsmath}
\usepackage{amssymb}
\usepackage[utf8]{inputenc}
\usepackage[T1]{fontenc}
\usepackage{ifthen}
\usepackage{hyperref}
\usepackage[rotate,all,color]{xy}
\usepackage[figuresright]{rotating}
\usepackage{titleps}
\usepackage{floatpag}
\usepackage{graphicx}

\copyrightyear{2015}
\author{Björn Gohla}
\title{Skew Closed Structure of $\Gray$-Categories}
\thanks{This work was supported by the Portuguese Science Foundation
  (FCT) under the post-doc grant SFRH/BPD/99060/2013. The problem was
  problem was posed to the author by John Bourke.}

\newcommand{\from}{\colon}
\renewcommand{\to}{\longrightarrow}
\newcommand{\To}{\Longrightarrow}
\newcommand{\Tto}{\Rrightarrow}

\newtheorem{thm}{Theorem}    

\newtheorem{lem}[thm]{Lemma}
\newtheorem{clm}[thm]{Claim}
\newtheorem{rem}[thm]{Remark}
\newtheorem{defn}[thm]{Definition}
\newenvironment{prf}{\proof}{}
\newcommand{\defterm}[2][]{{\bf
    #2}\ifthenelse{\equal{#1}{}}{\index{#2}}{\index{#1}}}

\newcommand{\catfont}[1]{\mathsf{#1}}
\def\mathcatdef#1{\expandafter\def\csname#1\endcsname{\catfont{#1}}}
\mathcatdef{Set}
\mathcatdef{Cat}
\mathcatdef{Mod}
\mathcatdef{Span}
\mathcatdef{Map}
\mathcatdef{Psd}
\mathcatdef{Lax}
\mathcatdef{Hom}
\mathcatdef{Trihom}
\mathcatdef{El}
\mathcatdef{GC}
\mathcatdef{Gray}

\newcommand{\gcatfont}[1]{\mathcal{#1}}
\def\mathgcatdef#1{\expandafter\def\csname#1\endcsname{\gcatfont{#1}}}
\mathgcatdef{G}
\mathgcatdef{H}
\mathgcatdef{K}
\mathgcatdef{P}

\def\1{{\pmb{1}}}
\def\fQ{\mathrm{Q}}
\def\pathspc{\overrightarrow}

\def\dblbarspc#1{\overline{\overline{#1}}}
\def\trplbarspc#1{\overline{\overline{\overline{#1}}}}
\def\pair#1{\langle{#1}\rangle}
\newcommand{\lhc}{\triangleleft}
\newcommand{\rhc}{\triangleright}

\renewcommand{\epsilon}{\varon}

\renewcommand{\phi}{\varphi}
\renewcommand{\theta}{\vartheta}
\newcommand{\id}{{\mathrm {id}}}
\newcommand{\op}{{\mathrm {op}}}
\newcommand{\ten}{\mathord{\otimes}}

\newcommand{\BOX}{\hbox {$\sqcap$ \kern -1em $\sqcup$}}
\def\qed {{
    \parfillskip=0pt        
    \widowpenalty=10000     
    \displaywidowpenalty=10000  
    \finalhyphendemerits=0  
    \leavevmode             
    \unskip                 
    \nobreak                
    \hfil                   
    \penalty50              
    \hskip.2em              
    \null                   
    \hfill                  
    \BOX
    \par                    
  }
}

\makeatletter
\def\xyboxmatrix{%
  \def\xyboxmatrixname@{}%
  \def\xyboxmatrixsetup@{}%
  \xyFN@\xyboxmatrix@%
}
\def\xyboxmatrix@{%
  \ifx"\next\DNii@{\xyboxmatrix@prefix}%
  \else\addAT@\ifx\next\DNii@{\xyboxmatrix@setup}%
  \else\DNii@{\xyboxmatrix@ii}%
  \fi\fi\nextii@%
}
\def\xyboxmatrix@prefix"#1"{%
  \def\xyboxmatrixname@{#1}%
  \xyboxmatrix@setup%
}
\def\xyboxmatrix@setup#1#{%
  \def\xyboxmatrixsetup@{#1}%
  \xyboxmatrix@ii%
}
\def\xyboxmatrix@ii#1{%
  \def\xyboxmatrix@contents{#1}%
  \edef\xyboxmatrix@iii{\noexpand\POS*+!C\noexpand\xybox{%
      \noexpand\POS\noexpand\xymatrix\xyboxmatrixsetup@{%
        \noexpand\xyboxmatrix@contents}}="\xyboxmatrixname@"}%
  \xyboxmatrix@iii
}
\makeatother

\let\oldlabelbox=\labelbox
\def\newlabelbox#1{%
  \oldlabelbox{\vcenter{\normalbaselines%
      \let\\=\cr\ialign{\hspace{-.7ex}$\labelstyle##\hfil$\hspace{.7ex}\crcr#1\crcr}}
  }  
}
\let\labelbox\newlabelbox

\def\tria{\NoRules\save\afterPOS{\POS **{} ?<(.1) *=0{}="zero" ; ?>(.9) : %
    "zero"+/v(0,1).7ex/ ::%
    \PATH ~={**\dir{-}} '(0,1) '(1,0) '(0,-1) '(0,1) %
    \restore};p;}

\def\revise#1{\textcolor{blue}{#1}}

\newenvironment{saves}{\save}{\restore}

\def\ttar{%
  \afterPOS{%
    \save%
    \POS
    @(%
    **{}?/-1ex/@+c;?/1ex/@+c%
    \ar@3s1;s0%
    @)%
    \restore%
  }
}

\def\tar{%
  \afterPOS{%
    \save%
    \POS
    @(%
    **{}?/-1ex/@+c;?/1ex/@+c%
    \ar@2s1;s0%
    @)%
    \restore%
  }
}

\def\areq{%
  \afterPOS{%
    \save%
    \POS
    @(%
    **{}?/-1ex/@+c;?/1ex/@+c%
    \ar@{=}s1;s0%
    @)%
    \restore%
  }
}

\def\tarb#1{%
  \afterPOS{%
    \POS @(%
    **{}?/-1ex/@+c;?/1ex/@+c%
    \ar@2s1;s0_{#1}="here"%
    \POS @)%
    ,"here"
  }
}
\def\tara#1{%
  \afterPOS{%
    \POS @(%
    **{}?/-1ex/@+c;?/1ex/@+c%
    \ar@2s1;s0^{#1}="here"%
    \POS @)%
    ,"here"
  }
}

\def\ttarb#1{%
  \afterPOS{%
    \POS @(%
    **{}?/-1ex/@+c;?/1ex/@+c%
    \ar@3s1;s0_{#1}="here"%
    \POS @)%
    ,"here"
  }
}

\def\ttara#1{%
  \afterPOS{%
    \POS @(%
    **{}?/-1ex/@+c;?/1ex/@+c%
    \ar@3s1;s0^{#1}="here"%
    \POS @)%
    ,"here"
  }
}

\def\ttar{%
  \afterPOS{%
    \POS @(%
    **{}?/-1ex/@+c;?/1ex/@+c%
    \ar@3s1;s0|{}="here"%
    \POS @)%
    ,"here"
  }
}

\def\pla#1{%
  \afterPOS{%
    \POS @(%
    **{}?/-1ex/@+c;?/1ex/@+c%
    \ar@{}s1;s0|{#1}="here"%
    \POS @)%
    ,"here"%
  }
}

\def\offlabel#1{%
  \save%
  \afterPOS{%
    \POS%
    *!C+\labelbox{#1}%
    **\dir{.}%
    \restore%
  }%
}

\def\tbdef#1#2#3{%
  \def\t(##1,##2,##3){\POS{0;\tp(##1,##2,##3)}}%
  \def\tp(##1,##2,##3){\POS{p+#1**{}?(##1);p+#2**{}?(##2);p+#3**{}?(##3)}}%
}

\begin{document}
\maketitle{}

\begin{abstract}
  We define a skew-closed ste \citep{street13} for
  $\Gray$-categories extening the mapping space construction of
  \citet{gohla2014}. 
\end{abstract}

\section{Definitions}
\label{sec:definitions}

Let $\Gray\Cat$ be the category $\Gray$-categories and strict
$\Gray$-enriched functors. In \citep{gohla2014} we showed that there is
for $\Gray$-categories $\G$, $\H$ a $\Gray$-category $[\G,\H]$
consisting of pseudo-functors, pseudo-transformations,
pseudo-modifications and perturbations. This construction can be
restricted to a $\Gray$-category $[\G,\H]^p$ of strict $\Gray$-functors,
pseudo-transformations, pseudo-modifications and perturbations. 

We denote the composition operartions in $[\G,\H]^p$ by $*_{i}$ where
$i$ is the dimension of the intervening cell, i.e.\ we have
$\alpha'*_0\alpha\from{}H\To\H''$ for pseudo-transformations
$\alpha\from{}H\To{}H'$ and $\alpha'\from{}H'\To{}H''$ where
$H,H',H''$ are $\Gray$-functors $\G\to\H$; $\alpha$ and $\alpha'$ are
1-cells in $[\G,\H]^p$ that coincide on $H'$. We will denote by
$*_{-1}$ the composition of functors, pseudo-transformations,
pseudo-modifications and perturbations incident on a $\Gray$-category,
i.e.\ from the point of view of the mapping space $[\G,\H]^p$ the
$\Gray$-categories $\G$ and $\H$ have dimension $-1$.

\begin{clm}
  \label{lem:pbrackfunc}
  $[\G,\H]^p$ is functorial $\G$ and $\H$, meaning
  $[\_,\_]^p\from\Gray\Cat^\op\times\Gray\Cat\to\Gray\Cat$ is a functor.
\end{clm}
\begin{prf}
  That $[\_,\H]^p$ is a functor each $\H$ is seen easily, since it
  is defined in terms of pre-composition.

  Let $G\from\H\to\H'$ be a $\Gray$-functor, we define a mapping
  $[\G,G]^p\from[\G,\H]^p\to[\G,\H']$ as follows:
  \begin{itemize}
    \item $(F\from\G\to\H)\mapsto(GF\from\G\to\H')$
    \item $(\alpha\from\fQ^1\G\to\pathspc{\H})\mapsto(\pathspc{G}\alpha\from\fQ^1\G\to\pathspc{\H'})$
    \item $(A\from\fQ^1\G\to\dblbarspc{\H})\mapsto(\dblbarspc{G}A\from\fQ^1\G\to\dblbarspc{\H'})$
    \item $(\sigma\from\fQ^1\G\to\trplbarspc{\H})\mapsto(\trplbarspc{G}\sigma\from\fQ^1\G\to\trplbarspc{\H'})$
  \end{itemize}

  We check that $[\G,\_]^p$ preserves unit and composition, i.e.\ that
  $[\G,\id_{\H}]^p=\id_{[\G,\H]^p}$ and
  $[\G,G']^p[\G,G]^p=[\G,G'G]^p$, by using the fact that
  $\pathspc{(\_)}, \dblbarspc{(\_)}$ and $\trplbarspc{(\_)}$ are
  endo-functors on $\Gray\Cat$ and composition of $\Gray$-functors is
  associative.
  
  Finally, we need to verify that $[\G,G]^p$ itself is a
  $\Gray$-functor, meaning it needs to preserve the composition
  operations $*_i$ of $[\G,\H]^p$ and the tensor.

  Let $\alpha\from{}F\To{}F'$ and $\alpha'\from{}F'\To{}F''$ be
  pseudo-transformations for $\Gray$-functors
  $F,F',F''\from{}\G\to\H$. $\alpha$ and $\alpha'$ are given by
  $\Gray$-functors $\alpha, \alpha'\from\fQ^1\G\to\pathspc{H}$. Their
  composite along $F'$ is defined as
  $\alpha'*_0\alpha=m\fQ^1\pair{\alpha',\alpha}\delta$, hence
  $[\G,G](\alpha'*_0\alpha)=\pathspc{G}m\fQ^1\pair{\alpha',\alpha}\delta$. On
  the other hand
  $[\G,G](\alpha')*_0[\G,G](\alpha)=m\fQ^1\pair{\pathspc{G}\alpha',\pathspc{G}\alpha}\delta=m\fQ^1(\pathspc{G}\times_{d_0,d_1}\pathspc{G})\fQ^1\pair{\alpha',\alpha}\delta$. Thus,
  to prove
  $[\G,G](\alpha'*_0\alpha)=[\G,G](\alpha')*_0[\G,G](\alpha)$, it
  remains we to be shown that
  \begin{equation}\label{eq:mgcomm}
    \begin{xy}
      \xymatrix{
        \fQ^1(\pathspc{\H}\times_{d_0,d_1}\pathspc{\H})\ar[r]^-{m}\ar[d]_-{\fQ^1(\pathspc{G}\times\pathspc{G})}&\pathspc{\H}\ar[d]^-{\pathspc{G}}\\
        \fQ^1(\pathspc{\H'}\times_{d_0,d_1}\pathspc{\H'})\ar[r]_-{m}&\pathspc{\H'}
      }
    \end{xy}
  \end{equation} commutes in $\Gray\Cat$.
  This can be seen when comparing the two paths as composites of
  $\fQ^1$ graph maps \citep{gohla2014}. Their agreement as globular
  maps is obvious. 

  We proceed to check the equality of the 2-cocycle data: According to
  \citep[][section 5]{gohla2014} the 2-cocyle on a 2-by-2 arrangement of
  squares is given by 
  \begin{align*}
    &(\hat{f''}\#_0g'_2\#_0g_0)\\\#_1&(\hat{g'_2}\ten{}g_2)\\\#_1&(\hat{g'_1}\#_0\hat{g_2}\#_0f)
  \end{align*}
  and because applying $\pathspc{G}$ of a strict $\Gray$-functor to
  cells of $\pathspc{\H}$ is given by applying $G$ elementwise it is
  obvious that \eqref{eq:mgcomm} commutes.

  For the remaining operations one can argue similarly.  
  \qed
\end{prf}

First, we need for every triple of $\Gray$-categories
$\G,\H,\K$ a $\Gray$-functor
\begin{equation}
  \label{eq:1}
  L_{\H,\K}^{\G}\from[\H,\K]^p\to[[\G,\H]^p,[\G,\K]^p]^p
\end{equation}
natural in $\H$ and $\K$ and extranatural in $\G$. For brevity we
shall just write $L$. Let $H, \beta, B, \Delta$ be respectively a
0-,1-,2- and 3-cell from $[\H,\K]^p$, then we define $L$ as giving the
action by postcomposition:
\begin{align}
  \begin{split}
    \label{eq:Ldef}
    L(H)&=H*_{-1}(\_)\\ 
    L(\beta)&=\beta*_{-1}(\_)\\ 
    L(B)&=B*_{-1}(\_)\\ 
    L(\Delta)&=\Delta*_{-1}(\_)\,.
  \end{split}
\end{align} 
The composite $*_{-1}$ is defined in section
\ref{sec:horiz-comp-grayc}. Furthermore, let $L(\beta)^2$ be the
family of perturbations indexed by composable pairs $\alpha,\alpha'$
of pseudo-transformations in $[\G,\H]^p$ defined by
\begin{equation}
  \label{eq:Ldeftwococ}
  (L(\beta)^2_{\alpha',\alpha})_x=\beta^2_{\alpha'_x,\alpha_x}\,.
\end{equation}

\begin{thm}
  \label{thm:Lwelldef}
  The map $L$ given in \eqref{eq:Ldef} and \eqref{eq:Ldeftwococ} is
  well defined, that is,
  \begin{itemize}
    \item $L(H)$ is a $\Gray$-functor (\ref{def:grayfunct}) $L(H)\from[\G,\H]^p\to[\G,\K]^p$,
    \item $L(\beta)$ is a pseudo-transformation (\ref{def:pstransf}) $L(\beta)\from{}L(H)\to{}L(H')$,
    \item $L(B)$ is a pseudo-modification (\ref{def:psmod}) $L(B)\from{}L(\beta)\to{}L(\beta')$,
    \item $L(\Delta)$ is a perturbation (\ref{def:perturb}) $L(\Delta)\from{}L(B)\to{}L(B')$.
  \end{itemize}
\end{thm}

\begin{prf}
  The claims are proved in lemmata \ref{lem:Lwelldef0},
  \ref{lem:Lwelldef1}, \ref{lem:Lwelldef2}, and \ref{lem:Lwelldef3}.
  \qed
\end{prf}

\begin{lem}
  \label{lem:Lwelldef0}
  $L(H)\from[\G,\H]^p\to[\G,\K]^p$ is a $\Gray$-functor for $H$ a
  $\Gray$-functor, i.e.\ a 0-cell of $[[\G,\H]^p,[\G,\K]^p]^p$.
\end{lem}
\begin{prf}
  This is obvious from $H$ being a $\Gray$-functor. \qed
\end{prf}

\begin{lem}
  \label{lem:Lwelldef1}
  Let $\beta\from{}H\to{}H'$ be a pseudo-transformation, i.e.\ a
  1-cell in $[\H,\K]^p$, then $L(\beta)\from{}L(H)\to{}L(H')$, is a
  pseudo-transformation (\ref{def:pstransf}) in a canonical way
  $L(\beta)\from{}L(H)\to{}L(H')$ i.e.\ a 1-cell in
  $[[\G,\H]^p,[\G,\K]^p]^p$.
\end{lem}
\begin{prf}
  First we give the details of definition of $L(\beta)$ according to
  \eqref{eq:Ldef} and \eqref{eq:Ldeftwococ}:
  \begin{enumerate}
    \item On 0-cells $G$ of $[\G,\H]^p$ we get pseudo-transformations
    $L(\beta)_G\from{}L(H)(G)\to{}L(H')(G)$ by whiskering
    \begin{equation}
      \label{eq:Lwelldef10}    
      L(\beta)_G=\beta*_{-1}G
    \end{equation}
    by our definition \eqref{eq:Ldef} above.
    According to section \ref{sec:horiz-comp-grayc}
    the components at cells of $\G$ are $(L(\beta)_G)_{x}=\beta_{Gx}$,
    $(L(\beta)_G)(f)=\beta_{Gf}$, $(L(\beta)_G)_{\phi}=\beta_{G\phi}$
    and 2-cocycle $(L(\beta)_G)^2_{f',f}=\beta^2_{Gf',Gf}$. And by
    theorem \ref{thm:hcompthm} $\beta*_{-1}G\from{}HG\to{}H'G$ is a
    pseudo-modification, i.e.\ a 1-cell in $[\G,\K]^p$.
    \item On 1-cells definition \eqref{eq:Ldef} $\alpha\from{}G\to{}G'$ of $[\G,\H]^p$ yields
    2-cells in $[\G,\K]^p$ as follows:
    \begin{equation}
      \label{eq:Lwelldef11}
      \begin{xy}
        \xyboxmatrix{
          L(H)(G)\ar[r]^{L(\beta)_G}\ar[d]_{L(H)(\alpha)}&L(H')(G)\ar[d]^{L(H')(\alpha)}\ar@2[dl]**{}?(.5)/-1ex/;?(.5)/1ex/^{L(\beta)_\alpha}\\
          L(H)(G')\ar[r]_{L(\beta)_{G'}}&L(H')(G')\\
        }
      \end{xy}=
      \begin{xy}
        \xyboxmatrix{
          HG\ar[r]^{\beta*_{-1}G}\ar[d]_{H*_{-1}\alpha}&H'G\ar[d]^{H'*_{-1}\alpha}\ar@2[dl]**{}?/-1ex/;?/+1ex/^{\beta*_{-1}\alpha}\\
          HG'\ar[r]_{\beta*_{-1}G'}&H'G'\\
        }
      \end{xy}\,.
    \end{equation} 
    We note that by \ref{thm:hcompthm} $\beta*_{-1}\alpha$ is a
    pseudo-modification, i.e.\ a 2-cell in $[\G,\K]^p$.
    \item On 2-cells $A\from\alpha\To\alpha'$ of $[\G,\H]^p$ definition \eqref{eq:Ldef}
    gives us
    \begin{multline}
      \label{eq:Lweldef12}
      \begin{xy}
        \xyboxmatrix"A"@+.5cm{
          L(H)(G)\ar[r]^{L(\beta)_G}\ar[d]|{L(H)(\alpha)}="x"\ar@/_3pc/[d]_{L(H)(\alpha')}="y"&L(H')(G)\ar[d]^{L(H')(\alpha)}\ar@2[dl]**{}?/-1ex/;?/1ex/^{L(\beta)_\alpha}\\
          L(H)(G')\ar[r]_{L(\beta)_{G'}}&L(H')(G')\POS\ar@2"x";"y"**{}?/-1ex/;?/1ex/^{L(H)(A)}
        }
        \POS + (80,0)
        \xyboxmatrix"B"@+.5cm{
          L(H)(G)\ar[r]^{L(\beta)_G}\ar[d]_{L(H)(\alpha')}&L(H')(G)\ar[d]|{L(H')(\alpha')}="x"\ar@2[dl]**{}?/-1ex/;?/1ex/^{L(\beta)_{\alpha'}}\ar@/^3pc/[d]^{L(H')(\alpha)}="y"\\
          L(H)(G')\ar[r]_{L(\beta)_{G'}}&L(H')(G')\POS\ar@2"y";"x"**{}?/-1ex/;?/1ex/^{L(H')(A)}
        }
        \ar@3"A";"B"**{}?/-1ex/;?/1ex/^{L(\beta)_A}    
      \end{xy}\\=
      \begin{xy}
        \xyboxmatrix"A"@+.5cm{
          HG\ar[r]^{\beta*_{-1}G}\ar[d]|{H*_{-1}\alpha}="x"\ar@/_3pc/[d]_{H*_{-1}\alpha'}="y"&H'G\ar[d]^{H'*_{-1}\alpha}\ar@2[dl]**{}?/-1ex/;?/1ex/^{\beta*_{-1}\alpha}\\
          HG'\ar[r]_{\beta*_{-1}G'}&H'G'\POS\ar@2"x";"y"**{}?/-1ex/;?/1ex/^{H*_{-1}A}
        }
        \POS + (80,0)
        \xyboxmatrix"B"@+.5cm{
          HG\ar[r]^{\beta*_{-1}G}\ar[d]_{H*_{-1}\alpha'}&H'G\ar[d]|{H'*_{-1}\alpha'}="x"\ar@2[dl]**{}?/-1ex/;?/1ex/_{\beta*_{-1}{\alpha'}}\ar@/^3pc/[d]^{H'*_{-1}\alpha}="y"\\
          HG'\ar[r]_{\beta*_{-1}{G'}}&H'G'\POS\ar@2"y";"x"**{}?/-1ex/;?/1ex/^{H'*_{-1}A}
        }
        \ar@3"A";"B"**{}?/-1ex/;?/1ex/^{\beta*_{-1}A}
      \end{xy}\,.
    \end{multline}
    \item By lemma \ref{lem:Lbeta2perturb} we get a family of
    perturbations indexed by composable pairs of 1-cells
    $\xymatrix@1{G\ar[r]^{\alpha}&G'\ar[r]^{\alpha'}&G''}$ in
    $[\G,\H]^p$:
    \begin{equation}
      \label{eq:Lbetatwo}    
      \begin{xy}
        \xyboxmatrix"A"@+.5cm{
          L(H)(G)\ar[r]^{L(\beta)_G}\ar[d]_{L(H)(\alpha)}&L(H')(G)\ar[d]^{L(H')(\alpha)}\ar@2[dl]**{}?/-1ex/;?/1ex/_{L(\beta)_\alpha}\\
          L(H)(G')\ar[r]^{L(\beta)_{G'}}\ar[d]_{L(H)(\alpha')}&L(H')(G')\ar[d]^{L(H')(\alpha')}\ar@2[dl]**{}?/-1ex/;?/1ex/^{L(\beta)_{\alpha'}}\\
          L(H)(G'')\ar[r]_{L(\beta)_{G''}}&L(H')(G'')
        }
        \POS;(78,0)
        \xyboxmatrix"B"@+.5cm{
          L(H)(G)\ar[r]^{L(\beta)_G}\ar[dd]_/2ex/{L(H)(\alpha'*_0\alpha)}&L(H')(G)\ar[dd]^{L(H')(\alpha'*_0\alpha)}\ar@2[ddl]**{}?/-1ex/;?/1ex/_{L(\beta)_{\alpha'*_0\alpha}}\\
          {}&{}\\
          L(H)(G'')\ar[r]_{L(\beta)_{G''}}&L(H')(G'')
        }
        \POS\ar@3"A";"B"**{}?/-1ex/;?/1ex/^{L(\beta)^2_{\alpha',\alpha}}
      \end{xy}
    \end{equation}
    each definined on 0-cells $x$ of $\G$ as $(L(\beta)^2_{\alpha',\alpha})_x=\beta^2_{\alpha'_x,\alpha_x}$,
    i.e.\,
    \begin{equation}
      \label{eq:betatwozero}
      \begin{xy}
        \xyboxmatrix"A"@+.5cm{
          HGx\ar[r]^{\beta_{Gx}}\ar[d]_{H\alpha_x}&H'Gx\ar[d]^{H'\alpha_x}\ar@2[dl]**{}?/-1ex/;?/1ex/^{\beta_{\alpha_x}}\\
          HG'x\ar[r]^{\beta_{G'x}}\ar[d]_{H\alpha'_x}&H'G'x\ar[d]^{H'\alpha_x}\ar@2[dl]**{}?/-1ex/;?/1ex/^{\beta_{\alpha'_x}}\\
          HG''x\ar[r]_{\beta_{G''x}}&H'G''x\\
        }
        \POS;(78,0)
        \xyboxmatrix"B"@+.5cm{
          HGx\ar[r]^{\beta_{Gx}}\ar[dd]_{H(\alpha'_x\#_0\alpha_x)}&H'Gx\ar[dd]^{H'(\alpha'_x\#_0\alpha_x)}\ar@2[ddl]**{}?/-1ex/;?/1ex/^{\beta_{\alpha_x'\#_0\alpha_x}}\\
          {}&{}\\
          HG''x\ar[r]_{\beta_{G''x}}&H'G''x
        }
        \POS\ar@3"A";"B"**{}?/-1ex/;?/1ex/^{\beta^2_{\alpha'_x,\alpha_x}}
      \end{xy}\,.
    \end{equation}
  \end{enumerate}
  
  We now check the axioms for a pseudo-transformation given in
  \ref{def:pstransf}.
  \begin{enumerate}
    \item We get $L(\beta)_{\id_G}=\id_{L(\beta)_G}$ because
    $\beta*_{-1}\id_G=\id_{\beta*_{-1}G}$ by \ref{thm:pasteunit}. 
    \item For each 3-cell in $[\G,\H]^p$, i.e.\ a perturbation
    $\Gamma\from{}A\Tto{}A'\from\alpha\To\alpha'$ the following square
    commmutes:
    \begin{equation}
      \label{eq:Lweldef3cell}
      \def\outerbase{%
        \POS 0;<7.5cm,0cm>:(0,-.7)::0%
      }
      \def\subdiagram#1#2{%
        \xyboxmatrix"#1"@+.5cm{#2}%
      }
      \begin{xy}
        \outerbase
        \subdiagram{A}{
          L(H)(G) \ar[r]^{L(\beta)_G}\ar[d]|{L(H)(\alpha)}="x"\ar@/_5pc/[d]|{L(H)(\alpha')}="y"&L(H')(G)\ar[d]|{L(H')(\alpha)}\ar@2[dl]**{}?/-1ex/;?/1ex/^{L(\beta)_\alpha}\\
          L(H)(G') \ar[r]_{{L(\beta)_{G'}}}&L(H')(G')\ar@2"x";"y"**{}?/-1ex/;?/1ex/^{L(H)(A)}
        } 
        +(1,0) 
        \subdiagram{B}{
          L(H)(G)\ar[r]^{L(\beta)_G}\ar[d]|{L(H)(\alpha')}&L(H')(G)\ar[d]|{L(H')(\alpha')}="x"\ar@2[dl]**{}?/-1ex/;?/1ex/^{L(\beta)_{\alpha'}}\ar@/^5pc/[d]|{L(H')(\alpha)}="y"\\
          L(H)(G')\ar[r]_{{L(\beta)_{G'}}}&L(H')(G')\ar@2"y";"x"**{}?/-1ex/;?/1ex/^{L(H')(A)}
        }
        ,(0,1)
        \subdiagram{C}{
          L(H)(G)\ar[r]^{L(\beta)_G}\ar[d]|{L(H)(\alpha)}="x"\ar@/_5pc/[d]|{L(H)(\alpha')}="y"&L(H')(G)\ar[d]|{L(H')(\alpha)}\ar@2[dl]**{}?/-1ex/;?/1ex/^{L(\beta)_\alpha}\\
          L(H)(G')\ar[r]_{{L(\beta)_{G'}}}&L(H')(G')\ar@2"x";"y"**{}?/-1ex/;?/1ex/^{L(H)(A')}
        }
        +(1,0) 
        \subdiagram{D}{
          L(H)(G)\ar[r]^{L(\beta)_G}\ar[d]|{L(H)(\alpha')}&L(H')(G)\ar[d]|{L(H')(\alpha')}="x"\ar@2[dl]**{}?/-1ex/;?/1ex/^{L(\beta)_{\alpha'}}\ar@/^5pc/[d]|{L(H')(\alpha)}="y"\\
          L(H)(G')\ar[r]_{{L(\beta)_{G'}}}& L(H')(G')\ar@2"y";"x"**{}?/-1ex/;?/1ex/^{L(H')(A')}
        }
        \ar@3"A";"B"**{}?/-1ex/;?/1ex/^{L(\beta)_A}
        \ar@3"C";"D"**{}?/-1ex/;?/1ex/_{L(\beta)_{A'}} 
        \ar@3"A";"C"**{}?/-1ex/;?/1ex/_{({L(\beta)_{G'}}*_0L(H)(\Gamma))\\*_1L(\beta)_\alpha} 
        \ar@3"B";"D"**{}?/-1ex/;?/1ex/^{L(\beta)_{\alpha'}\\*_1(L(H')(\Gamma)*_0L(\beta)_G)}
      \end{xy}\,.
    \end{equation}
    By definition the components in \eqref{eq:Lweldef3cell} are given
    by
    \begin{equation}
      \label{eq:Lweldef3cellcheck}
      \def\outerbase{%
        \POS 0;<7.5cm,0cm>:(0,-.7)::0%
      }
      \def\subdiagram#1#2{%
        \xyboxmatrix"#1"@+.5cm{#2}%
      }
      \begin{xy}
        \outerbase
        \subdiagram{A}{
          HG\ar[r]^{\beta*_{-1}G}\ar[d]|{H*_{-1}\alpha}="x"\ar@/_5pc/[d]|{H*_{-1}\alpha'}="y"&H'G\ar[d]|{H'*_{-1}\alpha}\ar@2[dl]**{}?/-1ex/;?/1ex/^{\beta*_{-1}\alpha}\\
          HG'\ar[r]_{\beta*_{-1}G'}&H'G'\ar@2"x";"y"**{}?/-1ex/;?/1ex/^{H*_{-1}A}
        } 
        +(1,0) 
        \subdiagram{B}{
          HG\ar[r]^{\beta*_{-1}G}\ar[d]|{H*_{-1}\alpha'}&H'G\ar[d]|{H'*_{-1}\alpha'}="x"\ar@2[dl]**{}?/-1ex/;?/1ex/^{\beta*_{-1}\alpha'}\ar@/^5pc/[d]|{H'*_{-1}\alpha}="y"\\
          HG'\ar[r]_{\beta*_{-1}G'}&H'G'\ar@2"y";"x"**{}?/-1ex/;?/1ex/^{H'*_{-1}A}
        }
        ,(0,1)
        \subdiagram{C}{
          HG\ar[r]^{\beta*_{-1}G}\ar[d]|{H*_{-1}\alpha}="x"\ar@/_5pc/[d]|{H*_{-1}\alpha'}="y"&H'G\ar[d]|{H'*_{-1}\alpha}\ar@2[dl]**{}?/-1ex/;?/1ex/^{\beta*_{-1}\alpha}\\
          HG'\ar[r]_{\beta*_{-1}G'}&H'G'\ar@2"x";"y"**{}?/-1ex/;?/1ex/^{H*_{-1}A'}
        }
        +(1,0) 
        \subdiagram{D}{
          HG\ar[r]^{\beta*_{-1}G}\ar[d]|{H*_{-1}\alpha'}&H'*_{-1}G\ar[d]|{H'*_{-1}\alpha'}="x"\ar@2[dl]**{}?/-1ex/;?/1ex/^{\beta*_{-1}\alpha'}\ar@/^5pc/[d]|{H'*_{-1}\alpha}="y"\\
          HG'\ar[r]_{\beta*_{-1}G'}& H'G'\ar@2"y";"x"**{}?/-1ex/;?/1ex/^{H'*_{-1}A'}
        }
        \ar@3"A";"B"**{}?/-1ex/;?/1ex/^{\beta*_{-1}A}
        \ar@3"C";"D"**{}?/-1ex/;?/1ex/_{\beta*_{-1}A'}
        \ar@3"A";"C"**{}?/-1ex/;?/1ex/_{((\beta*_{-1}G')*_0(H*_{-1}\Gamma))\\*_1(\beta*_{-1}\alpha)}
        \ar@3"B";"D"**{}?/-1ex/;?/1ex/^{(\beta*_{-1}\alpha')\\*_1((H'*_{-1}\Gamma)*_0(\beta*_{-1}G))}
      \end{xy}\,.
    \end{equation} Since this is an equation of perturbations it is
    sufficient to compare their values at 0-cells $x$ of $\G$:
    \begin{equation}
      \label{eq:Lweldef3cellcheckx}
      \def\outerbase{%
        \POS 0;<7.5cm,0cm>:(0,-.7)::0%
      }
      \def\subdiagram#1#2{%
        \xyboxmatrix"#1"@+.5cm{#2}%
      }
      \begin{xy}
        \outerbase
        \subdiagram{A}{
          HGx\ar[r]^{\beta_{Gx}}\ar[d]|{H\alpha_x}="x"\ar@/_5pc/[d]|{H\alpha'_x}="y"&H'Gx\ar[d]|{H'\alpha_x}\ar@2[dl]**{}?/-1ex/;?/1ex/^{\beta_{\alpha_x}}\\
          HG'x\ar[r]_{\beta_{G'x}}&H'G'\ar@2"x";"y"**{}?/-1ex/;?/1ex/^{HA_x}
        } 
        +(1,0) 
        \subdiagram{B}{
          HGx\ar[r]^{\beta_{Gx}}\ar[d]|{H\alpha'x}&H'Gx\ar[d]|{H'\alpha'x}="x"\ar@2[dl]**{}?/-1ex/;?/1ex/^{\beta_{\alpha'_x}}\ar@/^5pc/[d]|{H'\alpha_x}="y"\\
          HG'x\ar[r]_{\beta_{G'x}}&H'G'x\ar@2"y";"x"**{}?/-1ex/;?/1ex/^{H'Ax}
        }
        ,(0,1)
        \subdiagram{C}{
          HGx\ar[r]^{\beta_{Gx}}\ar[d]|{H\alpha_x}="x"\ar@/_5pc/[d]|{H\alpha'_x}="y"&H'Gx\ar[d]|{H'\alpha_x}\ar@2[dl]**{}?/-1ex/;?/1ex/^{\beta_{\alpha_x}}\\
          HG'x\ar[r]_{\beta_{G'x}}&H'G'x\ar@2"x";"y"**{}?/-1ex/;?/1ex/^{HA'_x}
        }
        +(1,0) 
        \subdiagram{D}{
          HGx\ar[r]^{\beta_{Gx}}\ar[d]|{H\alpha'_x}&H'Gx\ar[d]|{H'\alpha'_x}="x"\ar@2[dl]**{}?/-1ex/;?/1ex/^{\beta_{\alpha'_x}}\ar@/^5pc/[d]|{H'\alpha_x}="y"\\
          HG'x\ar[r]_{\beta_{G'x}}& H'G'x\ar@2"y";"x"**{}?/-1ex/;?/1ex/^{H'A'_x}
        }
        \ar@3"A";"B"**{}?/-1ex/;?/1ex/^{\beta_{A_x}}
        \ar@3"C";"D"**{}?/-1ex/;?/1ex/_{\beta_{A'_x}}
        \ar@3"A";"C"**{}?/-1ex/;?/1ex/_{(\beta_{G'x}\#_0H\Gamma_x)\\\#_1\beta_{\alpha_x}}
        \ar@3"B";"D"**{}?/-1ex/;?/1ex/^{\beta_{\alpha'_x}\\\#_1(H'\Gamma_x\#_0\beta_{Gx})}
      \end{xy}\,,
    \end{equation}
    which holds by \eqref{eq:pstransf3cell}.
    \item A similar argument serves to verify equations
    \eqref{eq:pstransf2cellcomp} and \eqref{eq:pstransf2cellid}.
    \item We check that $L(\beta)^2$ as given by \eqref{eq:Lbetatwo}
    satisfies the cocycle condition \eqref{eq:pstransf2cocy}. But this
    is assured by the same condition for $\beta^2$: \eqref{eq:betatwozerococycle}
    \begin{sidewaysfigure}
      \begin{equation}
        \label{eq:betatwozerococycle}
        \begin{xy}
          \xyboxmatrix"A"@+.5cm{
            HGx\ar[r]^{\beta_{Gx}}\ar[d]_{H\alpha_x}&H'Gx\ar[d]^{H'\alpha_x}\ar@2[dl]**{}?/-1ex/;?/1ex/^{\beta_{\alpha_x}}\\
            HG'x\ar[r]^{\beta_{G'x}}\ar[d]_{H\alpha'_x}&H'G'x\ar[d]^{H'\alpha'_x}\ar@2[dl]**{}?/-1ex/;?/1ex/^{\beta_{\alpha'_x}}\\
            HG''x\ar[r]^{\beta_{G''x}}\ar[d]_{H\alpha''_x}&H'G''x\ar[d]^{H'\alpha''_x}\ar@2[dl]**{}?/-1ex/;?/1ex/^{\beta_{\alpha''_x}}\\
            HG'''x\ar[r]_{\beta_{G'''x}}&H'G'''x } \POS;(83,0)
          \xyboxmatrix"B"@+.5cm{
            HGx\ar[r]^{\beta_{Gx}}\ar[d]|/-1ex/{H(\alpha'_x\#_0\alpha_x)}&H'Gx\ar[d]^{H'(\alpha'_x\#_0\alpha_x)}\ar@2[dl]**{}?/-1ex/;?/1ex/^{\beta_{\alpha_x'\#_0\alpha_x}}\\
            HG''x\ar[r]^{\beta_{G''x}}\ar[d]_{H\alpha''_x}&H'G''x\ar[d]^{H'\alpha''_x}\ar@2[dl]**{}?/-1ex/;?/1ex/^{\beta_{\alpha_x''}}\\
            HG'''x\ar[r]_{\beta_{G'''x}}&H'G'''x } \POS;(0,-83)
          \xyboxmatrix"C"@+.5cm{
            HGx\ar[r]^{\beta_{Gx}}\ar[d]_{H\alpha_x}&H'Gx\ar[d]^{H'\alpha_x}\ar@2[dl]**{}?/-1ex/;?/1ex/^{\beta_{\alpha_x}}\\
            HG'x\ar[r]^{\beta_{G'x}}\ar[d]_{H(\alpha''_x\#_0\alpha'_x)}&H'G'x\ar[d]^{H'(\alpha''_x\#_0\alpha'_x)}\ar@2[dl]**{}?/-1ex/;?/1ex/^{\beta_{\alpha''_x\#_0\alpha'_x}}\\
            HG'''x\ar[r]_{\beta_{G'''x}}&H'G'''x\\
          } \POS;(83,-83) 
          \xyboxmatrix"D"@+.5cm{
            HGx\ar[r]^{\beta_{Gx}}\ar[dd]|/-2ex/{H(\alpha''_x\#_0\alpha'_x\#_0\alpha_x)}&H'Gx\ar[dd]|/-2ex/{H'(\alpha''_x\#_0\alpha'_x\#_0\alpha_x)}\ar@2@<1ex>[ddl]**{}?/-1ex/;?/1ex/|/3.1ex/{\beta_{\alpha''_x\#_0\alpha_x'\#_0\alpha_x}}\\
            {}&{}\\
            HG'''x\ar[r]_{\beta_{G'''x}}&H'G'''x }
          \POS\ar@3"A";"B"**{}?/-1ex/;?/1ex/^{(\beta_{\alpha''_x}\#_0H\alpha'_x\#_0H\alpha_x)\\\#_1(H'\alpha''_x\#_0\underline{\beta^2_{\alpha'_x,\alpha_x}})}
          \POS\ar@3"A";"C"**{}?/-1ex/;?/1ex/_{(\underline{\beta^2_{\alpha''_x,\alpha'_x}}\#_0H\alpha'_x)\\\#_1(H'\alpha''_x\#_0H'\alpha'_x\#_0\beta_{\alpha_x})}
          \POS\ar@3"C";"D"**{}?/-1ex/;?/1ex/_{\underline{\beta^2_{\alpha''_x\#_0\alpha'_x,\alpha_x}}}
          \POS\ar@3"B";"D"**{}?/-1ex/;?/1ex/^{\underline{\beta^2_{\alpha''_x,\alpha'_x\#_0\alpha_x}}}
        \end{xy}
      \end{equation}
    \end{sidewaysfigure}
    holds for all 0-cells $x$ in $\G$. Hence $L(\beta)^2$ is a 2-cocycle
    with components
    $L(\beta)_{\alpha',\alpha}^2\from\left(L(\beta)_{\alpha'}*_0L(G)(\alpha)\right)*_1\left(L(G')(\alpha')*_0L(\beta)_\alpha\right)\Tto{}L(\beta)_{\alpha'*_0\alpha}$.

    The cocycle $L(\beta)^2$ is also normalized, i.e.\ satisfies
    \eqref{eq:pstransf2cocnorm} because $\beta^2$ is so.

    \item Finally we check compatibility of the cocycle $L(\beta)^2$
    with whiskering of 2-cells, as shown in
    \eqref{eq:pstransf12whiskleft} and
    \eqref{eq:pstransf12whiskright}. Again, this it suffices to check
    on 0-cells $x$ of $[\G,\H]^p$ and $\beta$ being a
    pseudo-transformation.
  \end{enumerate}

  This completes the proof that $L(\beta)$ is a pseudo-transformation. \qed
\end{prf}

\begin{lem}
  \label{lem:Lwelldef2}
  For a 2-cell in $[\H,\K]^p$, i.e.\ a pseudo-modification, $L(B)$ is
  a pseudo-modification (\ref{def:psmod})
  $L(B)\from{}L(\beta)\To{}L(\beta')$
\end{lem}
\begin{prf}
  We begin by stating the defining components of $L(B)$:
  \begin{enumerate}
    \item On 0-cells $G$ of $[\G,\H]^p$
    \begin{equation}
      \label{eq:Lwelldef20}
      \def\subdiagram#1{\xyboxmatrix@+1cm{#1}}
      \begin{xy}
        \subdiagram{
          L(H)(G)\ar@/^1.5pc/[r]^{L(\beta)_G}="a"\ar@/_1.5pc/[r]_{L(\beta')_G}="b"&L(H')(G)
          \ar@2@<-1ex>"a";"b"**{}?/-1ex/;?/1ex/^{L(B)_G}
        }
      \end{xy}=
      \begin{xy}
        \subdiagram{
          H*_{-1}G\ar@/^1.5pc/[r]^{\beta*_{-1}G}="a"\ar@/_1.5pc/[r]_{\beta'*_{-1}G}="b"&H'*_{-1}G
          \ar@2@<-1ex>"a";"b"**{}?/-1ex/;?/1ex/^{B*_{-1}G}
        }
      \end{xy}\,,
    \end{equation}
    \item on 1-cells $\alpha$ of $[\G,\H]^p$
    \def\subdiagram#1#2{\xyboxmatrix"#1"@+1cm{#2}}      
    \def\outerbase{\POS 0;(75,0):(-1,0)::0,}
    \begin{multline}
      \label{eq:Lwelldef21}
      \begin{xy}
        \outerbase
        \subdiagram{A}{
          L(H)(G)\ar@/^1.5pc/[r]^{L(\beta)_G}="a"\ar@/_1.5pc/[r]_{L(\beta')_G}="b"\ar[d]_{L(H)(\alpha)}&
          L(H')(G)\ar[d]^{L(H')(\alpha)}\tara{L(\beta)_{\alpha'}}{;[d]**{}?};[dl]\\
          L(H)(G')\ar@/_1.5pc/[r]_{L(\beta')_{G'}}="b1"&
          L(H')(G')
          \ar@2"a";"b"**{}?/-1ex/;?/1ex/^{L(B)_G}
        }
        +(1,0)
        \subdiagram{B}{
          L(H)(G)\ar@/^1.5pc/[r]^{L(\beta)_G}="a"\ar[d]_{L(H)(\alpha)}&
          L(H')(G)\ar[d]^{L(H')(\alpha)}\tarb{L(\beta)_{\alpha}};{[l];[dl]**{}?}\\
          L(H)(G')\ar@/^1.5pc/[r]^{L(\beta)_{G'}}="a1"\ar@/_1.5pc/[r]_{L(\beta')_G'}="b1"&
          L(H')(G')
          \ar@2"a1";"b1"**{}?/-1ex/;?/1ex/^{L(B)_{G'}}
        }
        \ar@3"A";"B"**{}?/-1ex/;?/1ex/^{L(B)_\alpha}
      \end{xy}\\=
      \begin{xy}
        \outerbase
        \subdiagram{A}{
          HG\ar@/^1.5pc/[r]^{\beta*_{-1}G}="a"\ar@/_1.5pc/[r]_{\beta'*_{-1}G}="b"\ar[d]_{H*_{-1}\alpha}&
          H'G\ar[d]^{H'*_{-1}\alpha}\tara{\beta*_{-1}\alpha'}{;[d]**{}?};[dl]\\
          HG'\ar@/_1.5pc/[r]_{\beta'*_{-1}G'}="b1"&
          H'G'\ar@2"a";"b"**{}?/-1ex/;?/1ex/^{B*_{-1}G}
        }
        +(1,0)
        \subdiagram{B}{
          HG\ar@/^1.5pc/[r]^{\beta*_{-1}G}="a"\ar[d]_{H*_{-1}\alpha}&
          H'G\ar[d]^{H'*_{-1}\alpha}\tarb{\beta*_{-1}\alpha};{[l];[dl]**{}?}\\
          HG'\ar@/^1.5pc/[r]^{\beta*_{-1}G'}="a1"\ar@/_1.5pc/[r]_{\beta'*_{-1}G'}="b1"&
          H'G'\ar@2"a1";"b1"**{}?/-1ex/;?/1ex/^{B*_{-1}G'}
        }
        \ar@3"A";"B"**{}?/-1ex/;?/1ex/^{B*_{-1}\alpha}
      \end{xy}\,.
    \end{multline}
  \end{enumerate}
  And we verify the axioms for this to constitute a pseudo-modification:
  \begin{enumerate}
    \item $L(B)_{\id_G}=\id_{L(B)_{G}}$ by $B*_{-1}\id_G=\id_{B*_{-1}G}$ from \ref{thm:pasteunit}.
    \item Compatibility of $L(B)$ with the cocycles $L(\beta)^2$ and
    $L(\beta')^2$, i.e.\ \eqref{eq:modifexcoccomp}, is assured by $B$
    being a pseudo-modification.
    \item For 2-cells $A\from\alpha\To\alpha'\from{}G\to{}G'$ in
    $[\G,\H]^p$ \eqref{eq:modifex2cell} holds because it does so for
    $B$.\qed
  \end{enumerate}
\end{prf}

\begin{lem}
  \label{lem:Lwelldef3}
  For every 3-cell $\Delta\from{}B\Tto{}B'\from\beta\To\beta'$ of
  $[\H,\K]^p$, $L(\Delta)$ is a perturbation
  $L(\Delta)\from{}L(B)\Tto{}L(B')$.
\end{lem}
\begin{prf}
  The components of $L(\Delta)$ are given by
  $L(\Delta)_G=\Delta*_{-1}G$. But then $L(\Delta)$ satisfies the
  perturbation condition \eqref{eq:pertex1cell} because $\Delta$ does
  so for every $G$ and for every $x$.\qed
\end{prf}

  We define 
  \begin{equation}
    \label{eq:Lbetaalphax}
    (L(\beta)_\alpha)_x=
    \xymatrix@ru@+.7cm{GHx\ar[r]^{\beta_{Hx}}\ar[d]_{G\alpha_x}&G'Hx\ar[d]^{G'\alpha_x}\ar@2[dl]**{}?/-1ex/;?/1ex/^{\beta_{\alpha_x}}\\
      GH'x\ar[r]_{\beta_{H'x}}&G'H'x
    }\,,
  \end{equation}
  \begin{sidewaysfigure}\revise{replace with ref to appendix}
    \begin{center}$(L(\beta)_\alpha)_f=(\beta*_{-1}\alpha)_f$\end{center}
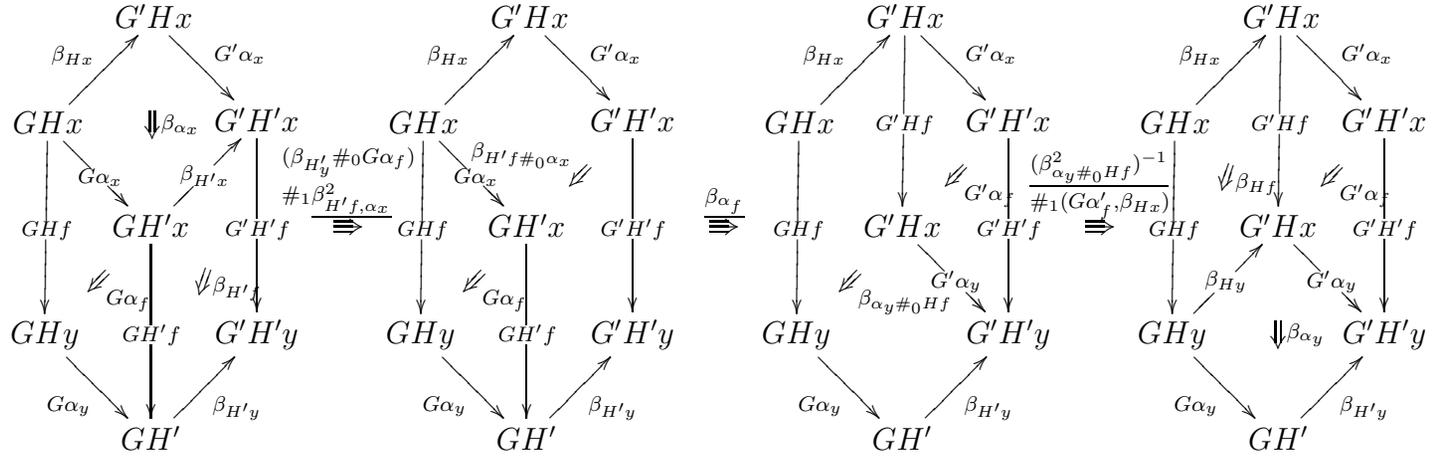
    \begin{multline}
      \label{eq:Lbetaalphaf}
      \begin{xy}
        \xyboxmatrix"A"@ru@+.3cm{
          {}&GHx\ar[r]^{\beta_{Hx}}\ar[d]|{G\alpha_x}\ar[dl]|{GHf}&G'Hx\ar[d]^{G'\alpha_x}\ar@2[dl]**{}?/-1ex/;?/1ex/^{\beta_{\alpha_x}}\\
          GHy\ar[d]_{G\alpha_y}&GH'x\ar[r]|{\beta_{H'x}}\ar[dl]|{GH'f}\ar@2[l]**{}?/-1ex/;?/+1ex/^{G\alpha_f}&G'H'x\ar[dl]|{G'H'f}\ar@2[dll]**{}?/-1ex/;?/+1ex/^{\beta_{H'f}}\\
          GH'\ar[r]_{\beta_{H'y}}&G'H'y&{}
        };(50,0)
        \xyboxmatrix"B"@ru@+.3cm{
          {}&GHx\ar[r]^{\beta_{Hx}}\ar[d]|{G\alpha_x}\ar[dl]|{GHf}&G'Hx\ar[d]^{G'\alpha_x}\\
          GHy\ar[d]_{G\alpha_y}&GH'x\ar[dl]|{GH'f}\ar@2[l]**{}?/-1ex/;?/+1ex/^{G\alpha_f}&G'H'x\ar[dl]|{G'H'f}\ar@2[l]**{}?/-1ex/;?/+1ex/_{\beta_{H'f\#_0\alpha_x}}\\
          GH'\ar[r]_{\beta_{H'y}}&G'H'y&{}
        };(100,0)
        \xyboxmatrix"C"@ru@+.3cm{
          {}&GHx\ar[r]^{\beta_{Hx}}\ar[dl]|{GHf}&G'Hx\ar[d]^{G'\alpha_x}\ar[dl]|{G'Hf}\\
          GHy\ar[d]_{G\alpha_y}&G'Hx\ar[d]|{G'\alpha_y}\ar@2[l]**{}?/-1ex/;?/1ex/^{\beta_{\alpha_y\#_0Hf}}&G'H'x\ar[dl]|{G'H'f}\ar@2[l]**{}?/-1ex/;?/+1ex/^{G'\alpha_f}\\
          GH'\ar[r]_{\beta_{H'y}}&G'H'y&{}
        };(150,0)
        \xyboxmatrix"D"@ru@+.3cm{
          {}&GHx\ar[r]^{\beta_{Hx}}\ar[dl]|{GHf}&G'Hx\ar[d]^{G'\alpha_x}\ar[dl]|{G'Hf}\ar@2[dll]**{}?/-1ex/;?/1ex/^{\beta_{Hf}}\\
          GHy\ar[d]_{G\alpha_y}\ar[r]|{\beta_{Hy}}&G'Hx\ar[d]|{G'\alpha_y}\ar@2[dl]**{}?/-1ex/;?/1ex/^{\beta_{\alpha_y}}&G'H'x\ar[dl]|{G'H'f}\ar@2[l]**{}?/-1ex/;?/+1ex/^{G'\alpha_f}\\
          GH'\ar[r]_{\beta_{H'y}}&G'H'y&{} }
        \ar@3"A";"B"**{}?/-1ex/;?/1ex/^{(\beta_{H'_y}\#_0G\alpha_f)\\\#_1\underline{\beta^2_{H'f,\alpha_x}}}
        \ar@3"B";"C"**{}?/-1ex/;?/1ex/^{\underline{\beta_{\alpha_f}}}
        \ar@3"C";"D"**{}?/-1ex/;?/1ex/^{\underline{(\beta^2_{\alpha_y\#_0Hf})^{-1}}\\\#_1(G\alpha'_f,\beta_{Hx})}
      \end{xy}\,.
    \end{multline}
  \end{sidewaysfigure}

  We pause here to note for later reference the composites of
  pseudo-transformations along the boundaries in \eqref{eq:Lwelldef11},
  i.e. $L(G')(\alpha)L(\beta)_H=(G'*_{-1}\alpha)*_0(\beta*_{-1}H)=\beta\rhc\alpha$
  and
  $L(\beta)_{H'}L(G)(\alpha)=(\beta*_{-1}H')*_0(G*_{-1}\alpha)=\beta\lhc\alpha$;
  these are given by their values on 0-, 1- and 2-cells and their
  2-cocycle data.  The data for $\beta\rhc\alpha$ on 0- and 1-cells
  are:
  \begin{equation}
    \label{eq:rhcbetaalpha01}    
    \begin{xy}
      \xyboxmatrix{
        GHx\ar[r]^{\beta_{Hx}}\ar[d]_{GHf}&G'Hx\ar[r]^{G'\alpha_x}\ar[d]|{G'Hf}\ar@2[dl]**{}?/-1ex/;?/1ex/^{\beta_{Hf}}&G'H'x\ar[d]^{G'H'f}\ar@2[dl]**{}?/-1ex/;?/1ex/^{G'(\alpha_f)}\\
        GHy\ar[r]_{\beta_{Hy}}&G'Hy\ar[r]_{G'\alpha_y}&G'H'y
      }\,,
    \end{xy}
  \end{equation}
  on 2-cells: see \eqref{eq:rhcbetaalpha2},
  \begin{sidewaysfigure}
    \begin{center}$(\beta\rhc\alpha)_{\phi}$\end{center}
    \begin{equation}
      \label{eq:rhcbetaalpha2}
      \begin{xy}
        \xyboxmatrix"A"@+.5cm{
          GHx\ar[r]^{\beta_{Hx}}\ar[d]|{GHf}="x"\ar@/_3pc/[d]|{GHf'}="y"\ar@2"x";"y"**{}?/-1ex/;?/1ex/^{GH\phi}&G'Hx\ar[r]^{G'\alpha_x}\ar[d]|{G'Hf}\ar@2[dl]**{}?/-1ex/;?/1ex/^{\beta_{Hf}}&G'H'x\ar[d]^{G'H'f}\ar@2[dl]**{}?/-1ex/;?/1ex/^{G'(\alpha_f)}\\
          GHy\ar[r]_{\beta_{Hy}}&G'Hy\ar[r]_{G'\alpha_y}&G'H'y}
        \POS;(80,0) \xyboxmatrix"B"@+.5cm{
          GHx\ar[r]^{\beta_{Hx}}\ar[d]_{GHf'}&G'Hx\ar[r]^{G'\alpha_x}\ar@/_1.5pc/[d]|{G'Hf'}="x"\ar@/^1.5pc/[d]|{G'Hf}="y"\ar@2"y";"x"**{}?/-1ex/;?/1ex/^{G'H\phi}\ar@2[dl]**{}?/1ex/;?/3ex/_{\beta_{Hf'}}&G'H'x\ar[d]^{G'H'f}\ar@2[dl]**{}?/-3ex/;?/-1ex/^{G'(\alpha_f)}\\
          GHy\ar[r]_{\beta_{Hy}}&G'Hy\ar[r]_{G'\alpha_y}&G'H'y}
        \POS;(160,0) \xyboxmatrix"C"@+.5cm{
          GHx\ar[r]^{\beta_{Hx}}\ar[d]_{GHf'}&G'Hx\ar[r]^{G'\alpha_x}\ar[d]|{G'Hf'}\ar@2[dl]**{}?/-1ex/;?/1ex/^{\beta_{Hf'}}&G'H'x\ar[d]|{G'H'f'}="x"\ar@/^3pc/[d]|{G'H'f}="y"\ar@2"y";"x"**{}?/-1ex/;?/1ex/^{G'H'\phi}\ar@2[dl]**{}?/-1ex/;?/1ex/^{G'(\alpha_f)}\\
          GHy\ar[r]_{\beta_{Hy}}&G'Hy\ar[r]_{G'\alpha_y}&G'H'y}
        \ar@3"A";"B"**{}?/-1ex/;?/1ex/|{}="here"\ar@{.}"here";p+(0,30)*!C\labelbox{(G'\alpha_y\#_0\underline{\beta_{H\phi}})\\\#_1(G'\alpha_f\#_0\beta_{Hx})}
        \ar@3"B";"C"**{}?/-1ex/;?/1ex/|{}="here"\ar@{.}"here";p+(0,30)*!C\labelbox{(G'\alpha_y\#_0\beta_{\alpha_{f'}})\\\#_1(\underline{G'\alpha_\phi}\#_0\beta_{Hf})}
      \end{xy}
    \end{equation}
  \end{sidewaysfigure}
  2-cocycle: see \eqref{eq:rhcbetaalphacoc}.
  \begin{sidewaysfigure}
    \begin{center}$(\beta\rhc\alpha)^2_{f',f}$\end{center}
    \begin{equation}
      \label{eq:rhcbetaalphacoc}
      \begin{xy}
        \xyboxmatrix"A"{
          GHx\ar[r]^{\beta_{Hx}}\ar[d]_{GHf}&G'Hx\ar[r]^{G'\alpha_x}\ar[d]|{G'Hf}\ar@2[dl]**{}?/-1ex/;?/1ex/_{\beta_{Hf}}|{}="x"&G'H'x\ar[d]^{G'H'f}\ar@2[dl]**{}?/-1ex/;?/1ex/_{G'_{\alpha_f}}\\
          GHy\ar[r]^{\beta_{Hy}}\ar[d]_{GHf'}&G'Hy\ar[r]^{G'\alpha_y}\ar[d]|{G'Hf'}\ar@2[dl]**{}?/-1ex/;?/1ex/_{\beta_{Hf'}}&G'H'y\ar[d]^{G'H'f'}\ar@2[dl]**{}?/-1ex/;?/1ex/^{G'_{\alpha_{f'}}}|{}="y"\\
          GHz\ar[r]_{\beta_{Hz}}&G'Hz\ar[r]_{G'\alpha_z}&G'H'z
          \POS\tria "x";"y"
        }
        \POS;(70,0)
        \xyboxmatrix"B"{
          GHx\ar[r]^{\beta_{Hx}}\ar[d]_{GHf}&G'Hx\ar[r]^{G'\alpha_x}\ar[d]|{G'Hf}\ar@2[dl]**{}?/-1ex/;?/1ex/_{\beta_{Hf}}|{}="x"&G'H'x\ar[d]^{G'H'f}\ar@2[dl]**{}?/-1ex/;?/1ex/_{G'_{\alpha_f}}\\
          GHy\ar[r]^{\beta_{Hy}}\ar[d]_{GHf'}&G'Hy\ar[r]^{G'\alpha_y}\ar[d]|{G'Hf'}\ar@2[dl]**{}?/-1ex/;?/1ex/_{\beta_{Hf'}}&G'H'y\ar[d]^{G'H'f'}\ar@2[dl]**{}?/-1ex/;?/1ex/^{G'_{\alpha_{f'}}}|{}="y"\\
          GHz\ar[r]_{\beta_{Hz}}&G'Hz\ar[r]_{G'\alpha_z}&G'H'z
          \POS\tria "y";"x"
        }
        \POS+(70,0)
        \xyboxmatrix"C"{
          GHx\ar[r]^{\beta_{Hx}}\ar[dd]|{GH(f'\#_0f)}&G'Hx\ar[r]^{G'\alpha_x}\ar[dd]|{G'H(f'\#_0f)}\ar@2[ddl]**{}?/-3ex/;?/-1ex/_{\beta_{H(f'\#_f)}}|{}="x"&G'H'x\ar[dd]|{G'H'(f'\#_0f)}\ar@2[ddl]**{}?/-3ex/;?/-1ex/_{G'_{\alpha_{f'\#_0f}}}\\
          {}&{}&{}\\
          GHz\ar[r]_{\beta_{Hz}}&G'Hz\ar[r]_{G'\alpha_z}&G'H'z
        }
        \ar@3"A";"B"**{}?/-1ex/;?/1ex/|{}="here"\ar@{.}"here";p+(0,30)*!C\labelbox{(G'\alpha_z\#_0\beta_{\alpha_f}\#_0GHf)\\\#_1(\underline{G'\alpha_{f'}\ten\beta_{\alpha_f}})\\\#_1(G'H'f'\#_0G'\alpha_f\#_0\beta_{Hx})}
        \ar@3"B";"C"**{}?/-1ex/;?/1ex/|{}="here"\ar@{.}"here";p+(0,30)*!C\labelbox{(G'\alpha_z\#_0\underline{\beta^2_{Hf',Hf}})\\\#_1(G'\alpha^2_{f',f}\#_0\beta_{Hx})}
      \end{xy}
    \end{equation}
  \end{sidewaysfigure}
  
  For $\beta\lhc_{-1}\alpha$ on the other hand, the data on 0- and 1-cells
  are
  \begin{equation}
    \label{eq:betalchalpha01}
    \begin{xy}
      \xyboxmatrix{
        GHx\ar[r]^{G\alpha_x}\ar[d]_{GHf}&GH'x\ar[r]^{\beta_{H'x}}\ar[d]|{GH'f}\ar@2[dl]**{}?/-1ex/;?/1ex/^{G\alpha_f}&G'H'x\ar[d]^{G'H'f}\ar@2[dl]**{}?/-1ex/;?/1ex/^{\beta_{H'f}}\\
        GHy\ar[r]_{G\alpha_y}&GH'y\ar[r]_{\beta_{H'y}}&G'H'y
      }
    \end{xy}\,,
  \end{equation}
  on 2-cells 
  \begin{sidewaysfigure}
    \begin{center}$(\beta\lhc\alpha)_\phi$\end{center}
    \begin{equation}
      \label{eq:lhcbetaalpha2}
      \begin{xy}
        \xyboxmatrix"A"@+.5cm{
          GHx\ar[r]^{G\alpha_x}\ar[d]|{GHf}="x"\ar@/_3pc/[d]|{GHf'}="y"\ar@2"x";"y"**{}?/-1ex/;?/1ex/^{GH\phi}&G'Hx\ar[r]^{\beta_{H'x}}\ar[d]|{G'Hf}\ar@2[dl]**{}?/-1ex/;?/1ex/^{G(\alpha_f)}&G'H'x\ar[d]^{G'H'f}\ar@2[dl]**{}?/-1ex/;?/1ex/^{\beta_{H'f}}\\
          GHy\ar[r]_{G\alpha_y}&G'Hy\ar[r]_{\beta_{H'y}}&G'H'y}
        \POS;(80,0) \xyboxmatrix"B"@+.5cm{
          GHx\ar[r]^{G\alpha_x}\ar[d]_{GHf'}&G'Hx\ar[r]^{\beta_{H'x}}\ar@/_1.5pc/[d]|{G'Hf'}="x"\ar@/^1.5pc/[d]|{G'Hf}="y"\ar@2"y";"x"**{}?/-1ex/;?/1ex/^{G'H\phi}\ar@2[dl]**{}?/1ex/;?/3ex/_{G(\alpha_{f'})}&G'H'x\ar[d]^{G'H'f}\ar@2[dl]**{}?/-3ex/;?/-1ex/^{\beta_{H'f}}\\
          GHy\ar[r]_{G\alpha_y}&G'Hy\ar[r]_{\beta_{H'y}}&G'H'y}
        \POS;(160,0) \xyboxmatrix"C"@+.5cm{
          GHx\ar[r]^{G\alpha_x}\ar[d]_{GHf'}&G'Hx\ar[r]^{\beta_{H'x}}\ar[d]|{G'Hf'}\ar@2[dl]**{}?/-1ex/;?/1ex/^{G\alpha_{f'}}&G'H'x\ar[d]|{G'H'f'}="x"\ar@/^3pc/[d]|{G'H'f}="y"\ar@2"y";"x"**{}?/-1ex/;?/1ex/^{G'H'\phi}\ar@2[dl]**{}?/-1ex/;?/1ex/^{\beta_{H'f'}}\\
          GHy\ar[r]_{G\alpha_y}&G'Hy\ar[r]_{\beta_{H'y}}&G'H'y}
        \ar@3"A";"B"**{}?/-1ex/;?/1ex/|{}="here"\ar@{.}"here";p+(0,30)*!C\labelbox{(\beta_{H'y}\#_0\underline{G\alpha_\phi})\\\#_1(\beta_{H'f}\#_0G\alpha_x)}
        \ar@3"B";"C"**{}?/-1ex/;?/1ex/|{}="here"\ar@{.}"here";p+(0,30)*!C\labelbox{(\beta_{H'y}\#_0G\alpha_{f'})\\\#_1(\underline{\beta_{H\phi}}\#_0G\alpha_x)}
      \end{xy}
    \end{equation}
  \end{sidewaysfigure}
  finally, for the 2-cocycle: see \eqref{eq:lhcbetaalphacoc}.
  \begin{sidewaysfigure}
    \begin{center}$(\beta\lhc\alpha)^2_{f',f}$\end{center}
    \begin{equation}
      \label{eq:lhcbetaalphacoc}
      \begin{xy}
        \xyboxmatrix"A"{
          GHx\ar[r]^{G\alpha_x}\ar[d]_{GHf}&GH'x\ar[r]^{\beta_{H'x}}\ar[d]|{GH'f}\ar@2[dl]**{}?/-1ex/;?/1ex/_{G\alpha_f}|{}="x"&G'H'x\ar[d]^{G'H'f}\ar@2[dl]**{}?/-1ex/;?/1ex/_{\beta_{H'f}}\\
          GHy\ar[r]^{G\alpha_y}\ar[d]_{GHf'}&GH'y\ar[r]^{\beta_{H'x}}\ar[d]|{GH'f'}\ar@2[dl]**{}?/-1ex/;?/1ex/_{G\alpha_f'}&G'H'y\ar[d]^{G'H'f'}\ar@2[dl]**{}?/-1ex/;?/1ex/^{\beta_{H'f'}}|{}="y"\\
          GHz\ar[r]_{G\alpha_z}&GH'z\ar[r]_{\beta_{H'z}}&G'H'z
          \POS\tria "x";"y"
        }
        \POS;(70,0)
        \xyboxmatrix"B"{
          GHx\ar[r]^{G\alpha_x}\ar[d]_{GHf}&GH'x\ar[r]^{\beta_{H'x}}\ar[d]|{GH'f}\ar@2[dl]**{}?/-1ex/;?/1ex/_{G\alpha_f}|{}="x"&G'H'x\ar[d]^{G'H'f}\ar@2[dl]**{}?/-1ex/;?/1ex/_{\beta_{H'f}}\\
          GHy\ar[r]^{G\alpha_y}\ar[d]_{GHf'}&GH'y\ar[r]^{\beta_{H'x}}\ar[d]|{GH'f'}\ar@2[dl]**{}?/-1ex/;?/1ex/_{G\alpha_f'}&G'H'y\ar[d]^{G'H'f'}\ar@2[dl]**{}?/-1ex/;?/1ex/^{\beta_{H'f'}}|{}="y"\\
          GHz\ar[r]_{G\alpha_z}&GH'z\ar[r]_{\beta_{H'z}}&G'H'z
          \POS\tria "y";"x"
        }
        \POS+(70,0)
        \xyboxmatrix"C"{
          GHx\ar[r]^{G\alpha_x}\ar[dd]|{GH(f'\#_0f)}&GH'x\ar[r]^{\beta_{H'x}}\ar[dd]|{GH'(f'\#_0f)}\ar@2[ddl]**{}?/-3ex/;?/-1ex/_{G\alpha_{f'\#_0f}}|{}="x"&G'H'x\ar[dd]|{G'H'(f'\#_0f)}\ar@2[ddl]**{}?/-3ex/;?/-1ex/_{\beta_{H'(f'\#_0f)}}\\
          {}&{}&{}\\
          GHz\ar[r]_{G\alpha_z}&GH'z\ar[r]_{\beta_{H'z}}&G'H'z
        }
        \ar@3"A";"B"**{}?/-1ex/;?/1ex/|{}="here"\ar@{.}"here";p+(0,30)*!C\labelbox{(\beta_{H'z}\#_0G\alpha_{f'}\#_0GHf)\\\#_1(\underline{\beta_{H'f'}\ten{}G\alpha_f})\\\#_1(G'H'f'\#_0\beta_{H'f}\#_0G\alpha_x)}
        \ar@3"B";"C"**{}?/-1ex/;?/1ex/|{}="here"\ar@{.}"here";p+(0,30)*!C\labelbox{(\beta_{H'z}\#_0\underline{G\alpha^2_{f',f}})\\\#_1(\underline{\beta^2_{H'f',H'f}}\#_0G\alpha_x)}
      \end{xy}
    \end{equation}
  \end{sidewaysfigure}

  In order to make sure that $L(\beta)_\alpha$ as shown in
  \eqref{eq:Lwelldef11} is indeed a pseudo-modification we need to check
  the following conditions on the 3-cell components of $L(\beta)$,
  $\beta\rhc\alpha$ and $\beta\lhc\alpha$, namely that
  \begin{enumerate}
    \item \revise{the $(L(\beta)_\alpha)_f$ are compatible with the
    $(\beta\rhc\alpha)^2_{f',f}$ and $(\beta\lhc\alpha)^2_{f',f}$,}
    \item \revise{for 2-cells $\phi\from{}f\To{}f'\from{}x\to{}y$ the 3-cells
    $(L(\beta)_\alpha)_f$, $(L(\beta)_\alpha)_{f'}$,
    $(\beta\rhc\alpha)_\phi$ and $(\beta\lhc\alpha)_\phi$ are compatible.}
  \end{enumerate}

  \revise{On 3-cells… i.e.\ $L(\beta)_\Gamma$}

\begin{lem}
  \label{lem:Lbeta2perturb}
  For a pair of composable 1-cells
  $\xymatrix@1{G\ar[r]^\alpha&G'\ar[r]^{\alpha'}&G''}$ in $[\G,\H]^p$
  $L(\beta)^2_{\alpha',\alpha}$ is indeed a perturbation.
\end{lem}
\begin{prf}
  We need to verify that $L(\beta)^2_{\alpha',\alpha}$ with
  $(L(\beta)^2_{\alpha',\alpha})_x=\beta^2_{\alpha'_x,\alpha_x}$
  satisfies \eqref{eq:pertex1cell}. That is, we need to verify that
  for every $x$ we have a commuting diagram:
  \begin{equation}
    \label{eq:Lbeta2perturbstate}
    \def\subdiagram#1#2{\xyboxmatrix"#1"@+2cm{#2}}
    \def\aofflabel#1#2{\offlabel{#2};p+#1}
    \def\frontlabel#1{\aofflabel{<0cm,-1cm>}{#1}}
    \def\backlabel#1{\aofflabel{<0cm,1cm>}{#1}}
    \begin{xy}
      \save (90,0):(0,-1):: 
      ,(0,0)="A11" 
      ,(1,0)="A21" 
      ,(0,1)="A12"
      ,(1,1)="A22"
      \restore 
      \POS,"A11"
      \subdiagram{A11}{
        L(H)(G)x\ar@/^2pc/[r]^{(L(H')(\alpha'*_0\alpha)*_0L(\beta)_G)_x}="a"\ar@/_2pc/[r]_{(L(\beta)_{G''}*_0L(H)(\alpha'*_0\alpha))_x}="b"\ar[d]|{L(H)(G)f}\tarb{}"a";"b"\aofflabel{<-3.1cm,1.1cm>}{(L(\beta)_{\alpha'}*_0L(H)(\alpha))\\*_1(L(H')(\alpha')*_0L(\beta)_\alpha)_x}&L(H')(G'')x\ar[d]|{L(H')(G'')f}\tarb{}[dl]\frontlabel{(L(\beta)_{G''}*_0L(H)(\alpha'*_0\alpha))_f}\\
        L(H)(G)y\ar@/_2pc/[r]_{(L(\beta)_{G''}*_0L(H)(\alpha'*_0\alpha))_y}&L(H')(G'')y
      }
      ,"A12"
      \subdiagram{A12}{
        L(H)(G)x\ar@/^2pc/[r]^{(L(H')(\alpha'*_0\alpha)*_0L(\beta)_G)_x}\ar[d]|{L(H)(G)f}&L(H')(G'')x\ar[d]|{L(H')(G'')f}\tara{}[dl]\backlabel{(L(H')(\alpha'*_0\alpha)*_0L(\beta)_G)_f}\\
        L(H)(G)x\ar@/^2pc/[r]^{(L(H')(\alpha'*_0\alpha)*_0L(\beta)_G)_y}="a1"\ar@/_2pc/[r]_{(L(\beta)_{G''}*_0L(H)(\alpha'*_0\alpha))_y}="b1"&L(H')(G'')y\tarb{}"a1";"b1"\aofflabel{<-3.1cm,-1.1cm>}{(L(\beta)_{\alpha'}*_0L(H)(\alpha))\\*_1(L(H')(\alpha')*_0L(\beta)_\alpha)_y}
      }   
      ,"A21"
      \subdiagram{A21}{
        L(H)(G)x\ar@/^2pc/[r]^{(L(H')(\alpha'*_0\alpha)*_0L(\beta)_G)_x}="a"\ar@/_2pc/[r]_{(L(\beta)_{G''}*_0L(H)(\alpha'*_0\alpha))_x}="b"\ar[d]|{L(H)(G)f}="a2"\tarb{}"a";"b"\aofflabel{<3.1cm,1.1cm>}{(L(\beta)_{\alpha'*_0\alpha})_x}&L(H')(G'')x\ar[d]|{L(H')(G'')f}\tarb{}[dl]\frontlabel{(L(\beta)_{G''}*_0L(H)(\alpha'*_0\alpha))_f}\\
        L(H)(G)y\ar@/_2pc/[r]_{(L(\beta)_{G''}*_0L(H)(\alpha'*_0\alpha))_y}="b1"&L(H')(G'')y
      }
      ,"A22"
      \subdiagram{A22}{
        L(H)(G)x\ar@/^2pc/[r]^{(L(H')(\alpha'*_0\alpha)*_0L(\beta)_G)_x}\ar[d]|{L(H)(G)f}&L(H')(G'')x\ar[d]|{L(H')(G'')f}\tara{}[dl]\backlabel{(L(H')(\alpha'*_0\alpha)*_0L(\beta)_G)_f}\\
        L(H)(G)y\ar@/^2pc/[r]^{(L(H')(\alpha'*_0\alpha)*_0L(\beta)_G)_y}="a1"\ar@/_2pc/[r]_{(L(\beta)_{G''}*_0L(H)(\alpha'*_0\alpha))_y}="b1"&L(H')(G'')y\tarb{}"a1";"b1"\aofflabel{<3.1cm,-1.1cm>}{(L(\beta)_{\alpha'*_0\alpha})_y}
      }
      \ttara{}"A11";"A21"\aofflabel{<0cm,3cm>}{(L(\beta)_{G''}*_0L(H)(\alpha'*_0\alpha))_f\\\#_1(L(H')(G'')f\#_0\underline{(L(\beta)^2_{\alpha',\alpha})_x})}
      \ttara{\underline{(L(\beta)_{\alpha'*_0\alpha})_f}}"A21";"A22"
      \ttarb{((L(\beta)_{\alpha'}*_0L(H)(\alpha))\\*_1(L(H')(\alpha')*_0L(\beta)_\alpha))_f}"A11";"A12"
      \ttarb{}"A12";"A22"\aofflabel{<0cm,-3cm>}{(\underline{(L(\beta)^2_{\alpha',\alpha})_y}\#_0L(H)(G)f)\\\#_1(L(H')(\alpha'*_0\alpha)*_0L(\beta)_G)_f}      
    \end{xy}\,.
  \end{equation}
  This commutes by \eqref{eq:Lbeta2perturbdetail}.
  \begin{sidewaysfigure}
    \begin{equation}
      \label{eq:Lbeta2perturbdetail}
      \def\outerbase{%
        \POS 0;<3.5cm,0cm>:(0,-1)::0%
      }
      \def\internalbase{\POS 0;<.1cm,0cm>:0,}
      \def\subdiagram#1#2{%
        \POS*!C++\xybox{%
          \internalbase%
          \vertices%
          \outerframe%
          #2%
        }="#1",%
      }
      \tbdef{(10,5)}{(8,-5)}{(0,-10)}
      \def\drawvertex#1{%
        \POS*{}="#1"%
      }
      \def\drawvertexat#1#2{%
        \t(#2)%
        \drawvertex{#1}%
      }
      \def\vertices{%
        \drawvertexat{a}{0,0,0}\drawvertexat{b}{1,0,0}%
        \drawvertexat{c}{0,1,0}\drawvertexat{d}{1,1,0}%
        \drawvertexat{e}{0,2,0}\drawvertexat{f}{1,2,0}%
        \drawvertexat{a'}{0,0,1}\drawvertexat{b'}{1,0,1}%
        \drawvertexat{c'}{0,1,1}\drawvertexat{d'}{1,1,1}%
        \drawvertexat{e'}{0,2,1}\drawvertexat{f'}{1,2,1}%
      }
      \def\outerframe{%
        \ar@{-}"a";"a'"%
        \ar@{-}"a";"b"%
        \ar@{-}"a'";"c'"%
        \ar@{-}"b";"d"%
        \ar@{-}"c'";"e'"%
        \ar@{-}"d";"f"%
        \ar@{-}"e'";"f'"%
        \ar@{-}"f";"f'"%
      }
      \begin{xy}
        \outerbase
        \subdiagram{A0}{
          \ar@{-}"a";"c"
          \ar@{-}"c";"d"
          \ar@{-}"c";"c'"
          \ar@{-}"c";"e"
          \ar@{-}"e";"e'"
          \ar@{-}"e";"f"
          \pla{\beta_{\alpha_x}}"b";"c"
          \pla{\beta_{\alpha'_x}}"d";"e"
          \pla{H\alpha_f}"c";"a'"
          \pla{H\alpha'_f}"e";"c'"
          \pla{\beta_{G''f}}"f";"e'"
        }
        +(1,0)
        \subdiagram{B0}{
          \ar@{-}"a";"c"
          \ar@{-}"c";"d"
          \ar@{-}"c";"c'"
          \ar@{-}"c";"e"
          \ar@{-}"e";"e'"
          \pla{\beta_{\alpha_x}}"b";"c"
          \pla{\beta_{G''f\#_0\alpha'_x}}"e";"f"
          \pla{H\alpha_f}"c";"a'"
          \pla{H\alpha'_f}"e";"c'"
        }
        +(1,0)
        \def\subdiagramC#1#2#3{%
          \subdiagram{#1}{%
            \ar@{-}"a";"c"%
            \ar@{-}"c";"c'"%
            \ar@{-}"c";"d"%
            \ar@{-}"d";"d'"%
            \ar@{-}"d'";"f'"%
            \pla{\beta_{\alpha_x}}"b";"c"\POS="X"%
            \pla{H\alpha_f}"c";"a'"%
            \pla{H'\alpha'_f}"f";"d'"\POS="Y"%
            \pla{\beta_{\alpha'_y\#_0G'f}}"d'";"c'"%
            \tria"#2";"#3"%
          }%
        }%
        \subdiagramC{C0}{X}{Y}
        +(1,0)
        \def\subdiagramA#1#2#3#4#5{%
          \subdiagram{#1}{%
            \ar@{-}"a";"c"%
            \ar@{-}"c";"c'"%
            \ar@{-}"c";"d"%
            \ar@{-}"d";"d'"%
            \ar@{-}"d'";"c'"%
            \ar@{-}"d'";"f'"%
            \pla{\beta_{\alpha_x}}"b";"c"\POS="X"%
            \pla{H\alpha_f}"c";"a'"\POS="X'"%
            \pla{H'\alpha'_f}"f";"d'"\POS="Y"%
            \pla{\beta_{G'f}}"d";"c'"%
            \pla{\beta_{\alpha'_y}}"f'";"c'"\POS="Y'"%
            \tria"#2";"#3"%
            \tria"#4";"#5"%
          }%
        }%
        \subdiagramA{D0}{X}{Y}{Y'}{X'}
        +(1,0)
        \subdiagramA{E0}{X}{Y}{X'}{Y'}
        (0,1)
        \subdiagram{A1}{
          \ar@{-}"a";"c"
          \ar@{-}"c";"c'"
          \ar@{-}"c";"e"
          \ar@{-}"e";"e'"
          \ar@{-}"e";"f"
          \pla{\beta_{\alpha'_x\#_0\alpha_x}}"d";"c"
          \pla{H\alpha_f}"c";"a'"
          \pla{H\alpha'_f}"e";"c'"
          \pla{\beta_{G''f}}"f";"e'"
        }
        (2,1)
        \subdiagramC{C1}{Y}{X}
        +(1,0)
        \subdiagramA{D1}{Y}{X}{Y'}{X'}      
        +(1,0)
        \subdiagramA{E1}{Y}{X}{X'}{Y'}
        (2,2)
        \subdiagram{C2}{
          \ar@{-}"a";"c"
          \ar@{-}"c";"c'"
          \ar@{-}"d";"d'"
          \ar@{-}"d'";"f'"
          \pla{H\alpha_f}"c";"a'"
          \pla{H'\alpha'_f}"f";"d'"
          \pla{\beta_{\alpha'_y\#_0G'f\#_0\alpha_x}}"d";"c"
        }
        +(1,0)
        \def\subdiagramB#1#2#3{%
          \subdiagram{#1}{%
            \ar@{-}"a";"c"%
            \ar@{-}"c";"c'"%
            \ar@{-}"d";"d'"%
            \ar@{-}"d'";"f'"%
            \ar@{-}"c'";"d'"%
            \pla{H\alpha_f}"c";"a'"\POS="X"%
            \pla{H'\alpha'_f}"f";"d'"%
            \pla{\beta_{G'f\#_0\alpha_x}}"d";"c"%
            \pla{\beta_{\alpha'_y}}"d'";"e'"\POS="Y"%
            \tria"#2";"#3"%
          }%
        }%
        \subdiagramB{D2}{Y}{X}
        +(1,0)
        \subdiagramB{E2}{X}{Y}
        (0,3)
        \subdiagram{A3}{
          \ar@{-}"a";"c"
          \ar@{-}"c";"c'"
          \ar@{-}"c";"e"
          \ar@{-}"e";"e'"
          \pla{\beta_{G''f\#_0\alpha'_x\#_0\alpha_x}}"d";"c"
          \pla{H\alpha_f}"c";"a'"
          \pla{H\alpha'_f}"e";"c'"
        }
        (4,3)
        \subdiagram{E3}{
          \ar@{-}"b";"b'"
          \ar@{-}"b'";"d'"
          \ar@{-}"c'";"d'"
          \ar@{-}"d";"d'"
          \ar@{-}"d'";"f'"
          \pla{H'\alpha_f}"d";"b'"
          \pla{H'\alpha'_f}"f";"d'"
          \pla{\beta_{\alpha_y\#_0\beta_{Gf}}}"b'";"a'"
          \pla{\beta_{\alpha'_y}}"f'";"c'"
        }
        (1,4)
        \subdiagram{B4}{
          \ar@{-}"b";"b'"
          \ar@{-}"b'";"d'"
          \ar@{-}"d";"d'"
          \ar@{-}"d'";"f'"
          \pla{H'\alpha_f}"d";"b'"
          \pla{H'\alpha'_f}"f";"d'"
          \pla{\beta_{\alpha'_y\#_0\alpha_y\#_0\beta_{Gf}}}"d'";"c'"
        }
        (3,4)
        \subdiagram{D4}{
          \ar@{-}"b";"b'"
          \ar@{-}"b'";"a'"
          \ar@{-}"b'";"d'"
          \ar@{-}"d";"d'"
          \ar@{-}"d'";"f'"
          \pla{\beta_{Gf}}"b";"a'"
          \pla{H'\alpha_f}"d";"b'"
          \pla{H'\alpha'_f}"f";"d'"
          \pla{\beta_{\alpha'_y\#_0\alpha_y}}"d'";"c'"
        }      
        (4,4)
        \subdiagram{E4}{
          \ar@{-}"b";"b'"
          \ar@{-}"b'";"a'"
          \ar@{-}"b'";"d'"
          \ar@{-}"d";"d'"
          \ar@{-}"d'";"c'"
          \ar@{-}"d'";"f'"
          \pla{\beta_{Gf}}"b";"a'"
          \pla{H'\alpha_f}"d";"b'"
          \pla{H'\alpha'_f}"f";"d'"
          \pla{\beta_{\alpha_y}}"d'";"a'"
          \pla{\beta_{\alpha'_y}}"f'";"c'"
        }      
        \ttara{}"A0";"A1"      
        \ttara{}"A0";"B0"
        \ttara{}"A1";"A3"
        \ttara{}"A3";"B4"
        \ttara{}"A3";"C2"
        \ttara{}"A3";"B0"
        \ttara{}"B0";"C0"
        \ttara{}"C0";"C1"
        \ttara{}"C0";"C1"      
        \ttara{}"C0";"D0"
        \ttara{}"C1";"C2"
        \ttara{}"C1";"D1"
        \ttara{}"C2";"B4"           
        \ttara{}"B4";"D4"
        \ttara{}"D0";"D1"
        \ttara{}"D0";"E0"
        \ttara{}"D1";"D2"
        \ttara{}"D1";"E1"
        \ttara{}"D2";"C2"
        \ttara{}"E0";"E1"
        \ttara{}"E1";"E2"
        \ttara{}"E2";"E3"
        \ttara{}"E2";"D2"      
        \ttara{}"E3";"B4"
        \ttara{}"E3";"E4"
        \ttara{}"E4";"D4"
        \pla{\txt{\eqref{eq:pstransf2cocy}}}"C2";"D1"
        \pla{\txt{\eqref{eq:pstransf2cocy}}}"A1";"B0"
        \pla{\txt{\eqref{eq:pstransf2cocy}}}"E3";"D4"
        \pla{\txt{\eqref{eq:pstransf12whiskleft}}}"A3";"C1"
        \pla{\txt{\eqref{eq:pstransf12whiskright}}}"B4";"D2"
        \pla{\txt{\eqref{eq:pstransf2cellcomp}}}"C2";{"A3";"B4"**{}?}
      \end{xy}
    \end{equation}
    \begin{center}
      Verification of \eqref{eq:Lbeta2perturbstate}; along the top and
      right is
      $((L(\beta)_{\alpha'}*_0L(H)(\alpha))*_1(L(H')(\alpha')*_0L(\beta)_\alpha))_f$,
      the three arrows along the lower left make up
      $(L(\beta)_{\alpha'*_0\alpha})_f$.
    \end{center}
  \end{sidewaysfigure}
  \qed
\end{prf}

\begin{lem}
  \label{lem:betacomp}
  For pseudo-transformations and $\Gray$-functors
  \begin{equation*}
    \begin{xy}
      \xyboxmatrix{
        \H\ar@/^2pc/[r]^{H}="x"\ar[r]|{H'}="y"\ar@/_2pc/[r]_{H''}="z"\tara{\beta}"x";"y"\tara{\beta'}"y";"z"&\K
      }
    \end{xy}
  \end{equation*}
  we have
  $L(\beta')*_{0}L(\beta)=L(\beta'*_{0}\beta)$, diagrammatically
  \begin{equation}
    \label{eq:betacompeq}    
    \begin{xy}
      \xyboxmatrix@+1cm{
        [\G,\H]^p\ar@/^2pc/[r]^{L(H)}="x"\ar[r]|{L(H')}="y"\ar@/_2pc/[r]_{L(H'')}="z"\ar@2"x";"y"**{}?/-1ex/;?/1ex/^{L(\beta)}\ar@2"y";"z"**{}?/-1ex/;?/1ex/^{L(\beta')}&[\G,\K]^p
      }
    \end{xy}=
    \begin{xy}
      \xyboxmatrix@+1cm{
        [\G,\H]^p\ar@/^1.5pc/[r]^{L(H)}="x"\ar@/_1.5pc/[r]_{L(H'')}="y"\ar@<-2ex>@2"x";"y"**{}?/-1ex/;?/1ex/^{L(\beta'*_0\beta)}&[\G,\K]^p
      }
    \end{xy}\,.
  \end{equation}
\end{lem}
\begin{prf}
  In order for the two pseudo-transformations to be equal we need
  their respective defining data to coincide:
  \begin{itemize}
    \item by lemma \ref{lem:betacomp0} for all $G\from\G\to\H$
    \begin{equation}
      \label{eq:betacompeq1}
      (L(\beta')*_0L(\beta))_G=L(\beta')_G*_0L(\beta)_G=L(\beta'*_0\beta)_G\,,
    \end{equation}    
    \item by lemma \ref{lem:betacomp1} for all $\alpha\from{}G\to{}G'$
    \begin{equation}
      \label{eq:betacompeq2}
      \begin{xy}
        \xyboxmatrix{
          L(H)(G)\ar[r]^{L(\beta)_G}\ar[d]_{L(H)(\alpha)}&L(H')(G)\ar[r]^{L(\beta')_G}\ar[d]|{L(H')(\alpha)}\ar@2[dl]**{}?(.5)/-1ex/;?(.5)/1ex/^{L(\beta)_\alpha}&L(H'')(G)\ar[d]^{L(H'')(\alpha)}\ar@2[dl]**{}?(.5)/-1ex/;?(.5)/1ex/^{L(\beta')_\alpha}\\
          L(H)(G')\ar[r]_{L(\beta)_{G'}}&L(H')(G')\ar[r]_{L(\beta')_{G'}}&L(H'')(G')\\
        }
      \end{xy}=
      \begin{xy}
        \xyboxmatrix{
          L(H)(G)\ar[rr]^{L(\beta'*_0\beta)_G}\ar[d]_{L(H)(\alpha)}&{}&L(H'')(G)\ar[d]^{L(H')(\alpha)}\ar@2[dll]**{}?(.5)/-1ex/;?(.5)/1ex/^{L(\beta'*_0\beta)_\alpha}\\
          L(H)(G')\ar[rr]_{L(\beta'*_0\beta)_{G'}}&{}&L(H'')(G')\\
        }
      \end{xy}\,,
    \end{equation}\revise{ fix notation
    \item by lemma \ref{lem:betacomp3} for all
    $A\from\alpha\To\alpha'$ the equation of 3-cells
    \eqref{eq:betacompcond3} holds,}
    \begin{sidewaysfigure}
      \begin{multline}
        \label{eq:betacompcond3}
        \begin{xy}
          \xyboxmatrix"A"@+.1cm{
            L(G)(H)\ar[r]^{L(\beta)_H}\ar[d]|{L(G)(\alpha)}="x"\ar@/_5pc/[d]_{L(G)(\alpha')}="y"\ar@2"x";"y"**{}?/-1ex/;?/1ex/^{L(G)(A)}&L(G')(H)\ar[r]^{L(\beta')_H}\ar[d]|{L(G')(\alpha)}\ar@2[dl]**{}?(.5)/-1ex/;?(.5)/1ex/^{L(\beta)_\alpha}&L(G'')(H)\ar[d]^{L(G'')(\alpha)}\ar@2[dl]**{}?(.5)/-1ex/;?(.5)/1ex/^{L(\beta')_\alpha}\\
            L(G)(H')\ar[r]_{L(\beta)_{H'}}&L(G')(H')\ar[r]_{L(\beta')_{H'}}&L(G'')(H')
          };(62.5,45) 
          \xyboxmatrix"B"@+.1cm{
            L(G)(H)\ar[r]^{L(\beta)_H}\ar[d]_{L(G)(\alpha')}&L(G')(H)\ar[r]^{L(\beta')_H}\ar@/_2.5pc/[d]|{L(G')(\alpha')}="x"\ar@/^2.5pc/[d]|{L(G')(\alpha)}="y"\ar@2"y";"x"**{}?/-1ex/;?/1ex/^{L(G')(A)}\ar@2[dl]**{}?/1.5ex/;?/3.5ex/_{L(\beta)_{\alpha'}}&L(G'')(H)\ar[d]^{L(G'')(\alpha)}\ar@2[dl]**{}?/-3.5ex/;?/-1.5ex/^{L(\beta')_\alpha}\\
            L(G)(H')\ar[r]_{L(\beta)_{H'}}&L(G')(H')\ar[r]_{L(\beta')_{H'}}&L(G'')(H')
          };(125,0)
          \xyboxmatrix"C"@+.1cm{
            L(G)(H)\ar[d]_{L(G)(\alpha')}\ar[r]^{L(\beta)_H}&L(G')(H)\ar[r]^{L(\beta')_H}\ar[d]|{L(G')(\alpha')}\ar@2[dl]**{}?(.5)/-1ex/;?(.5)/1ex/^{L(\beta)_{\alpha'}}&L(G'')(H)\ar[d]|{L(G'')(\alpha')}="x"\ar@/^5pc/[d]^{L(G'')(\alpha)}="y"\ar@2"y";"x"**{}?/-1ex/;?/1ex/^{L(G'')(A)}\ar@2[dl]**{}?(.5)/-1ex/;?(.5)/1ex/^{L(\beta')_{\alpha'}}\\
            L(G)(H')\ar[r]_{L(\beta)_{H'}}&L(G')(H')\ar[r]_{L(\beta')_{H'}}&L(G'')(H')
          }
          \ar@3"A";"B"**{}?/-1ex/;?/1ex/^{(L(\beta')_{H'}*_0\underline{L(\beta)_A})\\*_1(L(\beta')_\alpha)*_0L(\beta)_H}
          \ar@3"B";"C"**{}?/-1ex/;?/1ex/^{(L(\beta'_{H'})*_0L(\beta)_{\alpha'})\\*_1(\underline{L(\beta')_A}*_0L(\beta)_H)}
        \end{xy}\\=
        \begin{xy}
          \xyboxmatrix"A"@+.1cm{
            L(G)(H)\ar[rr]^{L(\beta'*_0\beta)_H}\ar[d]|{L(G)(\alpha)}="x"\ar@/_5pc/[d]_{L(G)(\alpha')}="y"\ar@2"x";"y"**{}?/-1ex/;?/1ex/^{L(G)(A)}&{}&L(G'')(H)\ar[d]|{L(G'')(\alpha)}\ar@2[dll]**{}?(.5)/-1ex/;?(.5)/1ex/^{L(\beta'*_0\beta)_{\alpha}}\\
            L(G)(H')\ar[rr]_{L(\beta'*_0\beta)_{H'}}&{}&L(G'')(H') 
          };(110,0) 
          \xyboxmatrix"B"@+.1cm{
            L(G)(H)\ar[rr]^{L(\beta'*_0\beta)_H}\ar[d]_{L(G)(\alpha')}&{}&L(G'')(H)\ar[d]|{L(G'')(\alpha')}="x"\ar@/^+5pc/[d]^{L(G'')(\alpha)}="y"\ar@2[dll]**{}?(.5)/-1ex/;?(.5)/1ex/^{L(\beta'*_0\beta)_{\alpha'}}\ar@2"y";"x"**{}?/-1ex/;?/1ex/^{L(G'')(A)}\\
            L(G)(H')\ar[rr]_{L(\beta'*_0\beta)_{H'}}&{}&L(G'')(H') 
          }
          \ar@3"A";"B"**{}?/-1ex/;?/1ex/^{\underline{L(\beta'*_0\beta)_A}}
        \end{xy}
      \end{multline}
    \end{sidewaysfigure}
    \item by lemma \ref{lem:betacomppairs}, finally, equation \eqref{eq:betacompcond2coc} holds for all composable pairs
    $\alpha\from{}H\To{}H'$, $\alpha'\from{}H'\To{}H''$.
    \begin{sidewaysfigure}
      \begin{multline}
        \label{eq:betacompcond2coc}
        \begin{xy}
          \xyboxmatrix"A"@+.1cm{
            L(G)(H)\ar[rr]^{(L(\beta')*_0L(\beta))_H}\ar[d]_{L(G)(\alpha)}&{}&L(G'')(H)\ar[d]^{L(G'')(\alpha)}\ar@2[dll]**{}?/-1ex/;?/1ex/_{(L(\beta')*_0L(\beta))_{\alpha}}\\
            L(G)(H')\ar[rr]^{(L(\beta')*_0L(\beta))_{H'}}\ar[d]_{L(G)(\alpha')}&{}&L(G)(H'')\ar[d]^{L(G'')(\alpha')}\ar@2[dll]**{}?/-1ex/;?/1ex/^{(L(\beta')*_0L(\beta))_{\alpha'}}\\
            L(G)(H'')\ar[rr]_{(L(\beta')*_0L(\beta))_{H''}}&{}&L(G'')(H'')
          };(85,0)
          \xyboxmatrix"B"@+.1cm{
            L(G)(H)\ar[r]^{L(\beta)_H}\ar[d]_{L(G)(\alpha)}&L(G')(H)\ar[r]^{L(\beta')_H}\ar[d]|{L(G')(\alpha)}\ar@2[dl]**{}?/-1ex/;?/1ex/_{L(\beta)_{\alpha}}|{}="x"&L(G'')(H)\ar[d]^{L(G'')(\alpha)}\ar@2[dl]**{}?/-1ex/;?/1ex/^{L(\beta')_{\alpha}}\\
            L(G)(H')\ar[r]^{L(\beta)_{H'}}\ar[d]_{L(G)(\alpha')}&L(G')(H')\ar[r]^{L(\beta')_{H'}}\ar[d]|{L(G')(\alpha')}\ar@2[dl]**{}?/-1ex/;?/1ex/^{L(\beta)_{\alpha'}}&L(G)(H'')\ar[d]^{L(G'')(\alpha')}\ar@2[dl]**{}?/-1ex/;?/1ex/^{L(\beta')_{\alpha'}}|{}="y"\\
            L(G)(H'')\ar[r]_{L(\beta)_{H''}}&L(G')(H'')\ar[r]_{L(\beta')_{H''}}&L(G'')(H'')
            \tria"y";"x"
          };(170,0)
          \xyboxmatrix"C"@+.1cm{
            L(G)(H)\ar[r]^{L(\beta)_H}\ar[dd]|{L(G)(\alpha'*_0\alpha)}&L(G')(H)\ar[r]^{L(\beta')_H}\ar[dd]|{L(G')(\alpha'*_0\alpha)}\ar@2[ddl]**{}?/-1ex/;?/1ex/^{L(\beta)_{\alpha'*_0\alpha}}&L(G'')(H)\ar[dd]|{L(G'')(\alpha'*_0\alpha)}\ar@2[ddl]**{}?/-1ex/;?/1ex/^{L(\beta')_{\alpha'*_0\alpha}}\\
            {}&{}&{}\\
            L(G)(H'')\ar[r]_{L(\beta)_{H''}}&L(G')(H'')\ar[r]_{L(\beta')_{H''}}&L(G'')(H'')
          }
          \ar@3"A";"B"**{}?<>(.5)/-1ex/;?<>(.5)/1ex/|{}="here"\ar@{.}"here";p+(0,30)*!C\labelbox{(L(\beta')_{H''}*_0L(\beta)_{\alpha'}*_0L(G)(\alpha))\\*_1(\underline{L(\beta')_{\alpha'}\ten{}L(\beta)_{\alpha}})\\*_1(L(G'')(\alpha')*_0L(\beta')(\alpha)*_0L(\beta)_H)}
          \ar@3"B";"C"**{}?/-1ex/;?/1ex/|{}="here"\ar@{.}"here";p+(0,30)*!C\labelbox{(L(\beta')_{H''}*_0\underline{L(\beta)^2_{\alpha',\alpha}})\\*_1(\underline{L(\beta')^2_{\alpha',\alpha}}*_0L(\beta)_H)}
        \end{xy}\\=
        \begin{xy}
          \xyboxmatrix"A"@+.1cm{
            L(G)(H)\ar[rr]^{L(\beta'*_0\beta)_H}\ar[d]_{L(G)(\alpha)}&{}&L(G'')(H)\ar[d]^{L(G'')(\alpha)}\ar@2[dll]**{}?/-1ex/;?/1ex/^{L(\beta'*_0\beta)_{\alpha}}\\
            L(G)(H')\ar[rr]^{L(\beta'*_0\beta)_{H'}}\ar[d]_{L(G)(\alpha')}&{}&L(G)(H'')\ar[d]^{L(G'')(\alpha')}\ar@2[dll]**{}?/-1ex/;?/1ex/^{L(\beta'*_0\beta)_{\alpha'}}\\
            L(G)(H'')\ar[rr]_{L(\beta'*_0\beta)_{H''}}&{}&L(G'')(H'')
          };(80,0)
          \xyboxmatrix"B"@+.1cm{
            L(G)(H)\ar[rr]^{L(\beta'*_0\beta)_H}\ar[dd]|{L(G)(\alpha'*_0\alpha)}&{}&L(G'')(H)\ar[dd]|{L(G'')(\alpha'*_0\alpha)}\ar@2[ddll]**{}?/-1ex/;?/1ex/^{L(\beta'*_0\beta)_{\alpha'*_0\alpha}}\\
            {}&{}&{}\\
            L(G)(H'')\ar[rr]_{L(\beta'*_0\beta)_{H''}}&{}&L(G'')(H'')
          }
          \ar@3"A";"B"**{}?/-1ex/;?/1ex/^{\underline{L(\beta'*_0\beta)^2_{\alpha',\alpha}}}
        \end{xy}
      \end{multline}
    \end{sidewaysfigure}
    \qed
  \end{itemize}
\end{prf}

In the sequel we prove the lemmas that constitute the substance of the
above proof.

\begin{lem}
  \label{lem:betacomp0}
  $L(\beta')_H*_0L(\beta)_H=L(\beta'*_0\beta)_H$ for all
  $H\from\G\to\H$.
\end{lem}
\begin{prf}
  On 0-cells, which are $\Gray$-functors $H\from\G\to\H$: we have a pseudo-transformation
  $(L(\beta')*_0L(\beta))_H\from{}GH\to{}G''H$ given on 0-cells $x$ of
  $\G$ by
  \begin{equation}
    \label{eq:LbetaLbetaHx}
    ((L(\beta')*_0L(\beta))_H)_x=((\beta'*_{-1}H)*_0(\beta*_{-1}H))_x=\beta'_{Hx}\#_0\beta_{Hx}\,,
  \end{equation}
  on 1-cells $f\from{}x\to{}y$
  \begin{equation}
    \label{eq:LbetaLbetaHf}
    ((L(\beta')*_0L(\beta))_H)_f=((\beta'*_{-1}H)*_0(\beta*_{-1}H))_f=
    \begin{xy}
      \xyboxmatrix{
        GHx\ar[r]^{\beta_{Hx}}\ar[d]_{GHf}&G'Hx\ar[r]^{\beta'_{Hx}}\ar[d]|{G'Hf}\ar@2[dl]**{}?/-1ex/;?/1ex/^{\beta_{Hf}}&G''Hx\ar[d]^{G''Hf}\ar@2[dl]**{}?/-1ex/;?/1ex/^{\beta'_{Hf}}\\
        GHy\ar[r]_{\beta_{Hy}}&G'Hy\ar[r]_{\beta'_{Hy}}&G''Hy
      }
    \end{xy}\,,
  \end{equation}
  on 2-cells by the diagram \eqref{eq:LbetaLbetaHphi}.
  \begin{sidewaysfigure}
    \begin{center}$((L(\beta')*_0L(\beta))_H)_\phi=(\beta'*_0\beta)_{H\phi}$\end{center}
    \begin{equation}
      \label{eq:LbetaLbetaHphi}
      \begin{xy}
        \xyboxmatrix"A"{
          GHx\ar[r]^{\beta_{Hx}}\ar[d]|{GHf}="x"\ar@/_3pc/[d]|{GHf'}="y"\ar@2"x";"y"**{}?/-1ex/;?/1ex/^{GH\phi}&G'Hx\ar[r]^{\beta'_{Hx}}\ar[d]|{G'Hf}\ar@2[dl]**{}?/-1ex/;?/1ex/^{\beta_{Hf}}&G''Hx\ar[d]^{G''Hf}\ar@2[dl]**{}?/-1ex/;?/1ex/^{\beta'_{Hf}}\\
          GHy\ar[r]_{\beta_{Hy}}&G'Hy\ar[r]_{\beta'_{Hy}}&G''Hy
        };(80,0)
        \xyboxmatrix"B"{
          GHx\ar[r]^{\beta_{Hx}}\ar[d]_{GHf'}&G'Hx\ar[r]^{\beta'_{Hx}}\ar@/_1.45pc/[d]|{G'Hf'}="x"\ar@/^1.45pc/[d]|{G'Hf}="y"\ar@2[dl]**{}?/1ex/;?/3ex/_{\beta_{Hf'}}\ar@2"y";"x"**{}?/-1ex/;?/1ex/^{G'H\phi}&G''Hx\ar[d]^{G''Hf}\ar@2[dl]**{}?/-3ex/;?/-1ex/^{\beta'_{Hf}}\\
          GHy\ar[r]_{\beta_{Hy}}&G'Hy\ar[r]_{\beta'_{Hy}}&G''Hy
        };(160,0)
        \xyboxmatrix"C"{
          GHx\ar[r]^{\beta_{Hx}}\ar[d]_{GHf'}&G'Hx\ar[r]^{\beta'_{Hx}}\ar[d]|{G'Hf'}\ar@2[dl]**{}?/-1ex/;?/1ex/^{\beta_{Hf'}}&G''Hx\ar[d]|{G''Hf'}="x"\ar@/^4pc/[d]|{G''Hf}="y"\ar@2[dl]**{}?/-1ex/;?/1ex/^{\beta'_{Hf'}}\ar@2"y";"x"**{}?/-1ex/;?/1ex/^{G''H\phi}\\
          GHy\ar[r]_{\beta_{Hy}}&G'Hy\ar[r]_{\beta'_{Hy}}&G''Hy
        }
        \ar@3"A";"B"**{}?/-1ex/;?/1ex/|{}="here"\ar@{.}"here";p+(0,30)*!C\labelbox{(\beta'_{Hy}\#_0\underline{\beta_{H\phi}})\\\#_1(\beta'_{Hf}\#_0\beta_{Hx})}
        \ar@3"B";"C"**{}?/-1ex/;?/1ex/|{}="here"\ar@{.}"here";p+(0,30)*!C\labelbox{(\beta'_{Hy}\#_0\beta_{Hf'})\\\#_1(\underline{\beta'_{H\phi}}\#_0'\beta_{Hx})}
      \end{xy}
    \end{equation}
  \end{sidewaysfigure} 
  On composable 1-cells $f',f$ from $\G$ we get the 2-cocycle shown in \eqref{eq:LbetaLbetaHtwo}.
  \begin{sidewaysfigure}
    \begin{center}
      $(L(\beta')*_0L(\beta))^2_{f',f}=(\beta'*_0\beta)^2_{Hf',Hf}$
        \revise{tensor goes the wrong way?}
        \revise{should be $((L(\beta')*_0L(\beta))_H)^2_{f',f}$}
    \end{center}
    \begin{equation}
      \label{eq:LbetaLbetaHtwo}
      \begin{xy}
        \xyboxmatrix"A"{
          GHx\ar[r]^{\beta_{Hx}}\ar[d]_{GHf}&G'Hx\ar[r]^{\beta'_{Hx}}\ar[d]|{G'Hf}\ar@2[dl]**{}?/-1ex/;?/1ex/_{\beta_{Hf}}="x"&G''Hx\ar[d]^{G''Hf}\ar@2[dl]**{}?/-1ex/;?/1ex/_{\beta'_{Hf}}\\
          GHy\ar[r]^{\beta_{Hy}}\ar[d]_{GHf'}&G'Hy\ar[r]^{\beta'_{Hy}}\ar[d]|{G'Hf'}\ar@2[dl]**{}?/-1ex/;?/1ex/^{\beta_{Hf'}}&G''Hy\ar[d]^{G''Hf'}\ar@2[dl]**{}?/-1ex/;?/1ex/^{\beta'_{Hf'}}="y"\\
          GHz\ar[r]_{\beta_{Hz}}&G'Hz\ar[r]_{\beta'_{Hz}}&G''Hz
          \tria"y";"x"
        };(80,0)
        \xyboxmatrix"B"{
          GHx\ar[r]^{\beta_{Hx}}\ar[d]_{GHf}&G'Hx\ar[r]^{\beta'_{Hx}}\ar[d]|{G'Hf}\ar@2[dl]**{}?/-1ex/;?/1ex/_{\beta_{Hf}}="x"&G''Hx\ar[d]^{G''Hf}\ar@2[dl]**{}?/-1ex/;?/1ex/_{\beta'_{Hf}}\\
          GHy\ar[r]^{\beta_{Hy}}\ar[d]_{GHf'}&G'Hy\ar[r]^{\beta'_{Hy}}\ar[d]|{G'Hf'}\ar@2[dl]**{}?/-1ex/;?/1ex/^{\beta_{Hf'}}&G''Hy\ar[d]^{G''Hf'}\ar@2[dl]**{}?/-1ex/;?/1ex/^{\beta'_{Hf'}}="y"\\
          GHz\ar[r]_{\beta_{Hz}}&G'Hz\ar[r]_{\beta'_{Hz}}&G''Hz
          \tria"x";"y"
        };(160,0)
        \xyboxmatrix"C"{
          GHx\ar[r]^{\beta_{Hx}}\ar[dd]|{GH(f'\#_0f)}&G'Hx\ar[r]^{\beta'_{Hx}}\ar[dd]|{G'H(f'\#_0f)}\ar@2[ddl]**{}?/-3.5ex/;?/-1.5ex/_{\beta_{H(f'\#_0f)}}="x"&G''Hx\ar[dd]|{G''H(f'\#_0f)}\ar@2[ddl]**{}?/-3.5ex/;?/-1.5ex/_{\beta'_{H(f'\#_0f)}}\\
          {}&{}&{}\\
          GHz\ar[r]_{\beta_{Hz}}&G'Hz\ar[r]_{\beta'_{Hz}}&G''Hz
        }
        \ar@3"A";"B"**{}?/-1ex/;?/1ex/|{}="here"\ar@{.}"here";p+(0,30)*!C\labelbox{(\beta'_{Hz}\#_0\beta_{Hf'}\#_0GHf)\\\#_1(\underline{\beta'_{Hf'}\ten\beta_{Hf}})\\\#_1(G''Hf'\#_0\beta'_{Hf}\#_0\beta_{Hx})}
        \ar@3"B";"C"**{}?/-1ex/;?/1ex/|{}="here"\ar@{.}"here";p+(0,30)*!C\labelbox{(\beta'_{Hz}\#_0\underline{\beta^2_{Hf',Hf}})\\\#_1(\underline{\beta'^2_{Hf',Hf}}\#_0\beta_{Hx})}
      \end{xy}
    \end{equation}
  \end{sidewaysfigure}

  \revise{Now let us consider the right hand side of
    \eqref{eq:betacompeq1}}, that is, we need to give
  $L(\beta'*_0\beta)_H$; this is simply substituting images under $H$
  into the composite spelled out in section
  \ref{sec:comp-pseu-transf}, which is given on 0-cells by
  \begin{equation}
    \label{eq:LbetabetaHx}
    (L(\beta'*_0\beta)_H)_x=((\beta'*_0\beta)*_{-1}H)_x=(\beta'*_0\beta)_{Hx}=\beta'_{Hx}\#_0\beta_{Hx}\,,
  \end{equation}
  on 1-cells $f\from{}x\to{}y$ by
  \begin{multline}
    \label{eq:LbetabetaHf}
    (L(\beta'*_0\beta)_H)_f=((\beta'*_0\beta)*_{-1}H)_f\\=(\beta'*_0\beta)_{Hf}=
    \begin{xy}
      \xyboxmatrix{
        GHx\ar[r]^{\beta_{Hx}}\ar[d]_{GHf}&G'Hx\ar[r]^{\beta'_{Hx}}\ar[d]|{G'Hf}\ar@2[dl]**{}?/-1ex/;?/1ex/^{\beta_{Hf}}&G''Hx\ar[d]^{G''Hf}\ar@2[dl]**{}?/-1ex/;?/1ex/^{\beta'_{Hf}}\\
        GHy\ar[r]_{\beta_{Hy}}&G'Hy\ar[r]_{\beta'_{Hy}}&G''Hy
      }
    \end{xy}\,,
  \end{multline}
  on 2-cells $\phi\from{}f\To{}f'$ by
  \begin{equation*}
    (L(\beta'*_0\beta)_H)_\phi=((\beta'*_0\beta)*_{-1}H)_\phi=(\beta'*_0\beta)_{H\phi}
  \end{equation*}
  with details given in \eqref{eq:LbetabetaHphi},
  \begin{sidewaysfigure}
    \begin{center}$(L(\beta'*_0\beta)_H)_\phi=((\beta'*_0\beta)*_{-1}H)_\phi=(\beta'*_0\beta)_{H\phi}$\end{center}
    \begin{equation}
      \label{eq:LbetabetaHphi}
      \begin{xy}
        \xyboxmatrix"A"{
          GHx\ar[r]^{\beta_{Hx}}\ar[d]|{GHf}="x"\ar@/_3pc/[d]_{GHf'}="y"\ar@2"x";"y"**{}?/-1ex/;?/1ex/^{GH\phi}&G'Hx\ar[r]^{\beta'_{Hx}}\ar[d]|{G'Hf}\ar@2[dl]**{}?/-1ex/;?/1ex/^{\beta_{Hf}}&G''Hx\ar[d]|{G''Hf}\ar@2[dl]**{}?/-1ex/;?/1ex/^{\beta'_{Hf}}\\
          GHy\ar[r]_{\beta_{Hy}}&G'Hy\ar[r]_{\beta'_{Hy}}&G''Hy} ;(70,0)
        \xyboxmatrix"B"{
          GHx\ar[r]^{\beta_{Hx}}\ar[d]_{GHf'}&G'Hx\ar[r]^{\beta'_{Hx}}\ar@/^1.5pc/[d]|{G'Hf}="x"\ar@/_1.5pc/[d]|{G'Hf'}="y"\ar@2"x";"y"**{}?/-1ex/;?/1ex/^{G'H\phi}\ar@2[dl]**{}?/1ex/;?/3ex/_{\beta_{f'}}&G''Hx\ar[d]^{G''Hf}\ar@2[dl]**{}?/-3ex/;?/-1ex/^{\beta'_{Hf}}\\
          GHy\ar[r]_{\beta_{Hy}}&G'Hy\ar[r]_{\beta'_{Hy}}&G''Hy} ;(140,0)
        \xyboxmatrix"C"{
          GHx\ar[r]^{\beta_{Hx}}\ar[d]|{GHf'}&G'Hx\ar[r]^{\beta'_{Hx}}\ar[d]|{G'Hf'}\ar@2[dl]**{}?/-1ex/;?/1ex/^{\beta_{f'}}&G''Hx\ar@/^3pc/[d]|{G''Hf}="x"\ar[d]|{G''Hf'}="y"\ar@2"x";"y"**{}?/-1ex/;?/1ex/^{G''H\phi}\ar@2[dl]**{}?/-1ex/;?/1ex/^{\beta'_{f'}}\\
          GHy\ar[r]_{\beta_{Hy}}&G'Hy\ar[r]_{\beta'_{Hy}}&G''Hy}
        \ar@3"A";"B"**{}?/-1ex/;?/1ex/|{}="here"\ar@{.}"here";p+(0,25)*!C\labelbox{(\beta'_{Hy}\#_0\underline{\beta_{H\phi}})\\\#_1(\beta'_{Hf}\#_0\beta_{Hx})}
        \ar@3"B";"C"**{}?/-1ex/;?/1ex/|{}="here"\ar@{.}"here";p+(0,25)*!C\labelbox{(\beta'_{Hy}\#_0\beta_{Hf'})\\\#_1(\underline{\beta'_{H\phi}}\#_0\beta_{Hx})}
      \end{xy}
    \end{equation}
  \end{sidewaysfigure}

  \revise{For composable pairs of 1-cells $f',f$ we first consider the
    2-cocycle associated to the pseudo-transformation $\beta'*_0\beta\from{}G\to{}G''$:
  }
\qed  
\end{prf}

\begin{lem}
  \label{lem:betacomp1}
  For pseudo-transformations $\alpha\from{}G\To{}G'$
  \eqref{eq:betacompeq2} holds.
\end{lem}
\begin{prf} \revise{According to theorem \ref{thm:betabetaalphapaste} } 
  On 1-cells, which are pseudo-transformations $\alpha\from{}G\To{}G'$
  we get a pseudo-modification
  \begin{multline}
    \label{eq:LbetaLbetaalpha}    
    (L(\beta')*_0L(\beta))_\alpha=(L(\beta')_{G'}*_{-1}L(\beta)_\alpha)*_0(L(\beta')_\alpha*_{-1}L(\beta)_G)\\
    =\begin{xy}
      \xyboxmatrix@+.5cm@ur{
        HG\ar[r]^{\beta*_{-1}G}\ar[d]_{H*_{-1}\alpha}&H'G\ar[r]^{\beta'*_{-1}G}\ar[d]|{H'*_{-1}\alpha}\ar@<-1ex>@2[dl]**{}?/-1ex/;?/1ex/^{\beta*_{-1}\alpha}&H''G\ar[d]^{H''*_{-1}\alpha}\ar@<-1ex>@2[dl]**{}?/-1ex/;?/1ex/^{\beta'*_{-1}\alpha}\\
        HG'\ar[r]_{\beta*_{-1}G'}&H'G'\ar[r]_{\beta'*_{-1}G'}&H''G'
      }
    \end{xy}=
    \begin{xy}
      \xyboxmatrix@+.5cm@ur{
        HG\ar[r]^{\beta*_{-1}G}\ar[d]_{H*_{-1}\alpha}&H'G\ar[r]^{\beta'*_{-1}G}&H''G\ar[d]^{H''*_{-1}\alpha}\ar@<-1ex>@2[dll]**{}?/-1ex/;?/1ex/^(-1){(\beta'*_0\beta)*_{-1}\alpha}\\
        HG'\ar[r]_{\beta*_{-1}G'}&H'G'\ar[r]_{\beta'*_{-1}G'}&H''G'
      }
    \end{xy}\\
    =L(\beta'*_{-1}\beta)_\alpha\,.
  \end{multline}
  \qed
\end{prf}

\begin{lem}
  \label{lem:betacomp2}

\end{lem}

\begin{lem}
  \label{lem:betacomppairs}

\end{lem}

  On 3-cells, i.e.\ pseudo-modifications
  $A\from\alpha{}\To\alpha'\from{}H\to{}H'$, we get a perturbation
  $(L(\beta)*_{-1}L(\beta))_A\from{}(L(\beta)*_{-1}L(\beta))_{\alpha}\Tto{}(L(\beta)*_{-1}L(\beta))_{\alpha'}$,
  which on 0-cells of $\G$ is given by
  \begin{multline}
    \label{eq:LbetaLbetaAx}
    ((L(\beta')*_0L(\beta))_A)_x\\=
    \begin{xy}
      \xyboxmatrix"A"@ur@+.5cm{
        GHx\ar[r]^{\beta_{Hx}}\ar[d]|{G\alpha_x}="x"\ar@/_3pc/[d]|{G\alpha'_x}="y"\ar@<-1ex>@2"x";"y"**{}?/-1ex/;?/1ex/^{GA_x}&G'Hx\ar[r]^{\beta'_{Hx}}\ar[d]|{G'\alpha_x}\ar@2[dl]**{}?/-1ex/;?/1ex/^{\beta_{\alpha_x}}&G''Hx\ar[d]^{G''\alpha_x}\ar@2[dl]**{}?/-1ex/;?/1ex/^{\beta'_{\alpha_x}}\\
        GH'x\ar[r]_{\beta_{H'x}}&G'H'x\ar[r]_{\beta'_{H'x}}&G''H'x
      };(60,0)
      \xyboxmatrix"B"@ur@+.5cm{
        GHx\ar[r]^{\beta_{Hx}}\ar[d]_{G\alpha'_x}&G'Hx\ar[r]^{\beta'_{Hx}}\ar@/^1.5pc/[d]|{G'\alpha_x}="x"\ar@/_1.5pc/[d]|{G'\alpha'_x}="y"\ar@<-1ex>@2"x";"y"**{}?/-1ex/;?/1ex/^{G'A_x}\ar@2[dl]**{}?/1ex/;?/3ex/^{\beta_{\alpha'_x}}&G''Hx\ar[d]^{G''\alpha_x}\ar@2[dl]**{}?/-3ex/;?/-1ex/^{\beta'_{\alpha_x}}\\
        GH'x\ar[r]_{\beta_{H'x}}&G'H'x\ar[r]_{\beta'_{H'x}}&G''H'x
      };(120,0)
      \xyboxmatrix"C"@ur@+.5cm{
        GHx\ar[r]^{\beta_{Hx}}\ar[d]_{G\alpha'_x}&G'Hx\ar[r]^{\beta'_{Hx}}\ar[d]|{G'\alpha_x}\ar@2[dl]**{}?/-1ex/;?/1ex/^{\beta_{\alpha'_x}}&G''Hx\ar@/^3pc/[d]|{G''\alpha_x}="x"\ar[d]|{G''\alpha'_x}="y"\ar@<-1ex>@2"x";"y"**{}?/-1ex/;?/1ex/^{G''A_x}\ar@2[dl]**{}?/-1ex/;?/1ex/^{\beta'_{\alpha'_x}}\\
        GH'x\ar[r]_{\beta_{H'x}}&G'H'x\ar[r]_{\beta'_{H'x}}&G''H'x
      };(70,0)
      \ar@3"A";"B"**{}?/-1ex/;?/1ex/|{}="here"\ar@{.}"here";p+(0,35)*!C\labelbox{(\beta'_{H'x}\#_0\underline{\beta_{(A_x)}})\\\#_1(\beta'_{\alpha_x}\#_0\beta_{Hx})}
      \ar@3"B";"C"**{}?/-1ex/;?/1ex/|{}="here"\ar@{.}"here";p+(0,35)*!C\labelbox{(\beta'_{H'_x}\#_0\beta_{\alpha'_x})\\\#_1(\underline{\beta'_{A_x}}\#_0\beta_{Hx})}
    \end{xy}\,.
  \end{multline}
  Finally, for the pseudo-transformation
  $L(\beta')*_0L(\beta)$ we need to give the 2-cocycle
  $(L(\beta')*_{0}L(\beta))^2_{\alpha',\alpha}$ for pairs of
  composable 1-cells from $[\G,\H]^p$, which means that for evey
  composable pair $\alpha,\alpha'$ of pseudo-transformations we need
  to give a perturbation
  \begin{equation}
    \label{eq:LbetaLbetatwo}
    \begin{xy}
      \xyboxmatrix"A"@+.5cm{
        GH\ar[d]_{G*_{-1}\alpha}\ar[r]^{(\beta'*_0\beta)_{H}}&G''H\ar[d]^{G''*_{-1}\alpha}\ar@2[dl]**{}?/-1ex/;?/1ex/^{(\beta'*_0\beta)_\alpha}\\
        GH'\ar[d]_{G*_{-1}\alpha'}\ar[r]_{(\beta'*_0\beta)_{H'}}&G''H'\ar[d]^{G''*_{-1}\alpha'}\ar@2[dl]**{}?/-1ex/;?/1ex/^{(\beta'*_0\beta)_{\alpha'}}\\
        GH''\ar[r]_{(\beta'*_0\beta)_{H''}}&G''H'' 
      }
      ;(75,0)
      \xyboxmatrix"B"@+.5cm{
        GH\ar[dd]|/3ex/{G*_{-1}(\alpha'*_0\alpha)}\ar[r]^{(\beta'*_0\beta)_{H}}&G''H\ar[dd]^{G''*_{-1}(\alpha'*_0\alpha)}\ar@2@<-3ex>[ddl]**{}?/-1ex/;?/1ex/^/2ex/{(\beta'*_0\beta)_{\alpha'*_0\alpha}}\\
        {}&{}\\
        GH''\ar[r]_{(\beta'*_0\beta)_{H''}}&G''H'' 
      }
      \ar@3"A";"B"**{}?/-1ex/;?/1ex/^{(\beta'*_0\beta)^2_{\alpha',\alpha}}
    \end{xy}\,.
  \end{equation}
  This perturbation takes values on 0-cells of $\G$:
  \begin{multline}
    \label{eq:LbetaLbetatwox}
    ((\beta'*_0\beta)^2_{\alpha',\alpha})_{x}\\
    =\begin{xy}
      \xyboxmatrix"A"@+.5cm@ur{
        GHx\ar[d]_{G\alpha_x}\ar[r]^{\beta'_{Hx}\#_0\beta_{Hx}}&G''Hx\ar[d]^{G''\alpha_x}\ar@2@<-3ex>[dl]**{}?/-1ex/;?/1ex/^{(\beta'*_0\beta)_{\alpha_x}}\\
        GH'x\ar[d]_{G\alpha'_x}\ar[r]|{\beta'_{H'x}\#_0\beta_{H'x}}&G''H'x\ar[d]^{G''\alpha'_x}\ar@2@<-3ex>[dl]**{}?/-1ex/;?/1ex/^{(\beta'*_0\beta)_{\alpha'_x}}\\
        GH''x\ar[r]_{\beta'_{H''x}\#_0\beta_{H''x}}&G''H''x
      }
      ;(75,0)
      \xyboxmatrix"B"@+.5cm@ur{
        GHx\ar[dd]_{G(\alpha'_x\#_0\alpha_x)}\ar[r]^{\beta'_{Hx}\#_0\beta_{Hx}}&G''Hx\ar[dd]^{G''(\alpha'_x\#_0\alpha_x)}\ar@2[ddl]**{}?/-1ex/;?/1ex/|/-2ex/{(\beta'*_0\beta)_{(\alpha'*_0\alpha)_x}}\\
        {}&{}\\
        GH''x\ar[r]_{\beta'_{H''x}\#_0\beta_{H''x}}&G''H''x
      }
      \ar@3"A";"B"**{}?/-1ex/;?/1ex/^{((\beta'*_0\beta)^2_{\alpha',\alpha})_x}
    \end{xy}\\=
    \begin{xy}
      \xyboxmatrix"A"@+.5cm@ur{
        GHx\ar[r]^{\beta_{Hx}}\ar[d]_{G\alpha_x}&G'Hx\ar[r]^{\beta'_{Hx}}\ar[d]|{G'\alpha_x}\ar@2[dl]**{}?/-1ex/;?/1ex/_{\beta_{\alpha_x}}|{}="x"&G''Hx\ar[d]^{G''\alpha_x}\ar@2[dl]**{}?/-1ex/;?/1ex/^{\beta'_{\alpha_x}}\\
        GH'x\ar[r]|{\beta_{H'x}}\ar[d]_{G\alpha'_x}&G'H'x\ar[r]|{\beta'_{H'x}}\ar[d]|{G'\alpha'_x}\ar@2[dl]**{}?/-1ex/;?/1ex/^{\beta_{\alpha'_x}}&G''H'x\ar[d]^{G''\alpha'_x}\ar@2[dl]**{}?/-1ex/;?/1ex/^{\beta'_{\alpha'_x}}|{}="y"\\
        GH''x\ar[r]_{\beta_{H''x}}&G'H''x\ar[r]_{\beta'_{H''x}}&G''H''x
        \tria"x";"y"
      };(80,0)
      \xyboxmatrix"B"@+.5cm@ur{
        GHx\ar[r]^{\beta_{Hx}}\ar[d]_{G\alpha_x}&G'Hx\ar[r]^{\beta'_{Hx}}\ar[d]|{G'\alpha_x}\ar@2[dl]**{}?/-1ex/;?/1ex/_{\beta_{\alpha_x}}|{}="x"&G''Hx\ar[d]^{G''\alpha_x}\ar@2[dl]**{}?/-1ex/;?/1ex/^{\beta'_{\alpha_x}}\\
        GH'x\ar[r]|{\beta_{H'x}}\ar[d]_{G\alpha'_x}&G'H'x\ar[r]|{\beta'_{H'x}}\ar[d]|{G'\alpha'_x}\ar@2[dl]**{}?/-1ex/;?/1ex/^{\beta_{\alpha'_x}}&G''H'x\ar[d]^{G''\alpha'_x}\ar@2[dl]**{}?/-1ex/;?/1ex/^{\beta'_{\alpha'_x}}|{}="y"\\
        GH''x\ar[r]_{\beta_{H''x}}&G'H''x\ar[r]_{\beta'_{H''x}}&G''H''x
        \tria"y";"x"
      }
      \ar@3"A";"B"**{}?/-1ex/;?/1ex/|{}="here"\ar@{.}"here";p+(0,30)*!C\labelbox{
        (\beta'_{H''x}\#_0\beta_{\alpha'_x}\#_0G\alpha_x)\\
        \#_1(\beta'_{\alpha'_x}\ten\beta_{\alpha_x})\\
        \#_1(G''\alpha'_x\#_0\beta_{H_x)}
      }
    \end{xy}\,.
  \end{multline}
\revise{  And on 1-cells the pseudo-modifications on the left and right hand
  sides of \eqref{eq:LbetaLbetatwo} are given by
  \begin{equation}
    \label{eq:LbetaLbetaleftf}
    \begin{xy}
      \xyboxmatrix"A"@+.5cm@ur{
        {}&GHx\ar[r]^{\beta_{Hx}}\ar[d]|{G\alpha_x}\ar[dl]_{GHf}&G'Hx\ar[r]^{\beta'_{Hx}}\ar[d]|{G'\alpha_x}\ar@2[dl]**{}?/-1ex/;?/1ex/_{\beta_{\alpha_x}}|{}="x"&G''Hx\ar[d]^{G''\alpha_x}\ar@2[dl]**{}?/-1ex/;?/1ex/^{\beta'_{\alpha_x}}\\
        GHy\ar[d]_{G\alpha_y}&GH'x\ar[r]|{\beta_{H'x}}\ar[d]|{G\alpha'_x}\ar[dl]|{GH'f}\ar@2[l]**{}?/-1ex/;?/1ex/^{G\alpha_f}&G'H'x\ar[r]|{\beta'_{H'x}}\ar[d]|{G'\alpha'_x}\ar@2[dl]**{}?/-1ex/;?/1ex/_{\beta_{\alpha'_x}}|{}="z"&G''H'x\ar[d]^{G''\alpha'_x}\ar@2[dl]**{}?/-1ex/;?/1ex/^{\beta'_{\alpha'_x}}|{}="y"\\
        GH'y\ar[d]_{G\alpha'_y}&GH''x\ar[r]|{\beta_{H''x}}\ar[dl]|{GH''f}\ar@2[l]**{}?/-1ex/;?/1ex/^{G\alpha'_f}&G'H''x\ar[r]|{\beta'_{H''x}}\ar[dl]|{G'H''f}\ar@2@<-1ex>[dll]**{}?/-1ex/;?/1ex/^{\beta_{H''f}}&G''H''x\ar[dl]^{G''H''f}\ar@2@<-1ex>[dll]**{}?/-1ex/;?/1ex/^{\beta'_{H''f}}|{}="w"\\
        GH''y\ar[r]_{\beta_{H''_y}}&G'H''y\ar[r]_{\beta'_{H''_y}}&GH''y&{}
        \tria"x";"y"
        \tria"z";"w"
      }
    \end{xy}
  \end{equation}
}

  \revise{ \qed}

A natural transformation with components
\begin{equation}
  \label{eq:2}
  i_\G\from[\1,\G]^p\to\G\,.
\end{equation}

An extranatural transformation with components
\begin{equation}
  \label{eq:3}
  j_\G\from\1\to[\G,\G]^p\,.
\end{equation}

\section{Functors and Transformations of $\Gray$-Categories}
\label{sec:funct-transf-gray}

We collect some defintions from the literature. 

\begin{defn}[e.g.\ {\citep{gray1974}, \citep{crans}}]
  \label{def:grayfunct}
  A \defterm{$\Gray$-category} is an enriched category, enriched over the
  category of 2-catgories with the $\Gray$-tensor product. A
  \defterm{$\Gray$-functor} is correspondingly an enriched functor.
\end{defn}

\begin{rem}
  \label{rem:grayfunct}
  A $\Gray$-functor $G\from\G\to\H$ maps $i$-cells to $i$-cells and
  preserves units, compositions, whiskers and tensors.
  
  In more detail, this means that the following equations hold.
  \begin{enumerate}
    \item Units: for $x,f,\phi,\Gamma$ 0-,1-,2-,3-cells respectively
    \begin{align}
      \label{eq:grayfunctunit}
        G\id_x&=\id_{Gx}\\
        G\id_f&=\id_{Gf}\\
        G\id_\phi&=\id_{G\phi}\\
        G\id_\Gamma&=\id_{G\Gamma}\,.
    \end{align}
    \item Whiskers: for suitably incident cells $f,g,\phi,\psi,\Gamma,\Delta$
    \begin{align}
      \label{eq:grayfunctwhisk}
      G(g\#_0\phi)&=Gg\#_0G\phi\\
      G(\psi\#_0f)&=G\psi\#_0Gf\\
      G(g\#_0\Gamma)&=Gg\#_0G\Gamma\\
      G(\Delta\#_0f)&=G\Delta\#_0Gf\\
      G(\phi\#_1\Delta)&=G\phi\#_0G\Delta\\
      G(\Gamma\#_1\psi)&=G\Gamma\#_1G\psi\,.
    \end{align}
    \item Composites: for suitably incident cells $f,g,\phi,\psi,\Gamma,\Delta$
    \begin{align}
      \label{eq:grayfunctcomp}
      G(g\#_0f)&=Gg\#_0Gf\\
      G(\psi\#_1\phi)&=G\psi\#_1G\phi\\
      G(\Delta\#_2\Gamma)&=G\Delta\#_2G\Gamma\,.
    \end{align}
    \item Tensors: for 2-cells incident on a 0-cell
    \begin{align}
      \label{eq::grayfunctten}
      G(\psi\ten\phi)&=G\psi\ten{}G\phi\,.
    \end{align}
  \end{enumerate}  
\end{rem}

The following notions were introduced in \citep{gohla2014} with some
more generality; we give only the ones used in this paper.

\begin{defn}
  \label{def:pstransf}
  For $\Gray$-functors $F,G\from\G\to\H$ a
  \defterm{pseudo-transformation} $\alpha\from{}F\to{}G$ is given by
  the following data: 
\begin{enumerate}
  \item for each 0-cell $x$ of $\G$ a 1-cell
    $\alpha_x\from{}Fx\to{}Gx$,
  \item for each 1-cell $f\from{}x\to{}y$ of $\G$ an invertible 2-cell
    \begin{equation}
      \label{eq:pstransf1cell}
      \begin{xy}
        \xyboxmatrix{
          Fx\ar[r]^{\alpha_x}\ar[d]_{Ff}&
          Gx\ar[d]^{Gf}\ar@2[dl]**{}?/-1ex/;?/1ex/^{\alpha_f}\\
          Fy\ar[r]_{\alpha_y}&Gy 
        }
      \end{xy}
    \end{equation}
  \item for each 2-cell $\phi\from{}f\to{}f'$ of $\G$ an invertible 3-cell of $\H$
    \begin{equation}
      \label{eq:pstransf2cell}
      \begin{matrix}
        \begin{xy}
          \xyboxmatrix"A"{Fx\ar[r]^{\alpha_x}\ar[d]|{Ff}="x"\ar@/_3pc/[d]_{Ff'}="y"&
            Gx\ar[d]^{Gf}\ar@2[dl]**{}?/-1ex/;?/1ex/^{\alpha_f}\\Fy
            \ar[r]_{\alpha_y}&Gy\POS\ar@2"x";"y"**{}?/-1ex/;?/1ex/^{F\phi}
          } \POS + (50,0)
          \xyboxmatrix"B"{Fx\ar[r]^{\alpha_x}\ar[d]_{Ff'}&
            Gx\ar[d]|{Gf'}="x"\ar@2[dl]**{}?/-1ex/;?/1ex/^{\alpha_{f'}}
            \ar@/^3pc/[d]^{Gf}="y"\\Fy\ar[r]_{\alpha_y}&Gy
            \POS\ar@2"y";"x"**{}?/-1ex/;?/1ex/^{G\phi}
          } \ar@3"A";"B"**{}?/-1ex/;?/1ex/^{\alpha_\phi} 
        \end{xy}
      \end{matrix}
    \end{equation}
  \item for each pair of composable 1-cells $f\from{}x\to{}y$,
    $f'\from{}y\to{}z$ an invertible 3-cell
    \begin{equation}
      \label{eq:pstransf2coch}
      \begin{xy}
        \xyboxmatrix"A"{
          Fx\ar[r]^{\alpha_x}\ar[d]_{Ff}&
          Gx\ar[d]^{Gf}\ar@2[dl]**{}?/-1ex/;?/1ex/^{\alpha_f}\\
          Fy\ar[r]_{\alpha_y}\ar[d]_{Ff'}&Gy\ar@2[dl]**{}?/-1ex/;?/1ex/^{\alpha_{f'}}\ar[d]^{Gf'}\\
          Fz\ar[r]_{\alpha_z}&Gz
        }\POS + (50,0)
        \xyboxmatrix"B"{
          Fx\ar[r]^{\alpha_x}\ar[dd]_{F(f'\#_0f)}&
          Gx\ar[dd]^{G(f'\#_0f)}="x"\ar@<+1.5ex>@2[ddl]**{}?/-1ex/;?/1ex/_{\alpha_{f'\#_0f}}\\
          {}&{}\\
          Fz\ar[r]_{\alpha_z}&Gz
        }
        \ar@3"A";"B"**{}?/-1ex/;?/1ex/^{\alpha^2_{f',f}}
      \end{xy}\,.
    \end{equation}
  \end{enumerate}
  These must satisfy the following conditions: 
  \begin{enumerate}
    \item On identities of 0-cells:
    \begin{equation}
      \label{eq:pstransf0id}
      \alpha_{\id_x}=\id_{\alpha_x}
    \end{equation}
    \item for each 3-cell
    $\Gamma\from{}\phi\to{}\phi'$ the square of 3-cells in $\H$
    \begin{equation}
      \label{eq:pstransf3cell}
      \begin{xy}
        \xyboxmatrix"A"{
          Fx \ar[r]^{\alpha_x}\ar[d]|{Ff}="x"\ar@/_3pc/[d]|{Ff'}="y"&Gx\ar[d]^{Gf}\ar@2[dl]**{}?/-1ex/;?/1ex/^{\alpha_f}\\
          Fy \ar[r]_{{\alpha_y}}&Gy\ar@2"x";"y"**{}?/-1ex/;?/1ex/^{F\phi}
        } 
        \POS + (50,0) 
        \xyboxmatrix"B"{
          Fx\ar[r]^{\alpha_x}\ar[d]_{Ff'}&Gx\ar[d]|{Gf'}="x"\ar@2[dl]**{}?/-1ex/;?/1ex/^{\alpha_{f'}}\ar@/^3pc/[d]|{Gf}="y"\\
          Fy\ar[r]_{{\alpha_y}}&Gy\ar@2"y";"x"**{}?/-1ex/;?/1ex/^{G\phi}
        }
        \POS + (-50,-50) 
        \xyboxmatrix"C"{
          Fx\ar[r]^{\alpha_x}\ar[d]|{Ff}="x"\ar@/_3pc/[d]|{Ff'}="y"&Gx\ar[d]^{Gf}\ar@2[dl]**{}?/-1ex/;?/1ex/^{\alpha_f}\\
          Fy\ar[r]_{{\alpha_y}}&Gy\ar@2"x";"y"**{}?/-1ex/;?/1ex/^{F\phi'}
        }
        \POS + (50,0) 
        \xyboxmatrix"D"{
          Fx\ar[r]^{\alpha_x}\ar[d]_{Ff'}&Gx\ar[d]|{Gf'}="x"\ar@2[dl]**{}?/-1ex/;?/1ex/^{\alpha_{f'}}\ar@/^3pc/[d]|{Gf}="y"\\
          Fy\ar[r]_{{\alpha_y}}& Gy\ar@2"y";"x"**{}?/-1ex/;?/1ex/^{G\phi'}
        }
        \ar@3"A";"B"**{}?/-1ex/;?/1ex/^{\alpha_\phi}
        \ar@3"C";"D"**{}?/-1ex/;?/1ex/_{\alpha_{\phi'}} 
        \ar@3"A";"C"**{}?/-1ex/;?/1ex/_{({\alpha_y}\#_0F\Gamma)\#_1\alpha_f} 
        \ar@3"B";"D"**{}?/-1ex/;?/1ex/^{\alpha_{f'}\#_1(G\Gamma\#_0\alpha_x)}
      \end{xy}
    \end{equation}commutes.
    \item  For every pair $\phi\from{}f\To{}f', \phi'\from{}f'\To{}f''$:
    \begin{equation}
      \label{eq:pstransf2cellcomp}
      \begin{xy}
        \xyboxmatrix"A"{
          Fx\ar[r]^{\alpha_x}\ar[d]|{Ff}="x"\ar@/_2.5pc/[d]|{Ff'}="y"\ar@/_5pc/[d]|{Ff''}="z"&Gx\ar[d]^{Gf}\ar@2[dl]**{}?/-1ex/;?/1ex/^{\alpha_f}\\
          Fy\ar[r]_{\alpha_y}&Gy\ar@2"x";"y"**{}?/-1ex/;?/1ex/^{F\phi}\ar@2"y";"z"**{}?/-1ex/;?/1ex/^{F\phi'} 
        } 
        \POS + (55, 0)
        \xyboxmatrix"B"{
          Fx\ar[r]^{\alpha_x}\ar[d]|{Ff'}="x"\ar@/_2.5pc/[d]|{Ff''}="y"&Gx\ar[d]|{Gf'}="w"\ar@/^2.5pc/[d]|{Gf}="z"\ar@2[dl]**{}?/-1ex/;?/1ex/^{\alpha_{f'}}\\
          Fx\ar[r]_{\alpha_{y}}&Gy\ar@2"x";"y"**{}?/-1ex/;?/1ex/^{F\phi'}\ar@2"z";"w"**{}?/-1ex/;?/1ex/^{G\phi} 
        }
        \POS + (55, 0)
        \xyboxmatrix"C"{
          Fx\ar[r]^{\alpha_x}\ar[d]_{Ff''}&Gx\ar[d]|{Gf''}="z"\ar@2[dl]**{}?/-1ex/;?/1ex/^{\alpha_{f''}}\ar@/^2.5pc/[d]|{Gf'}="y"\ar@/^5pc/[d]|{Gf}="x"\\
          Fy\ar[r]_{\alpha_y}&Gy\ar@2"x";"y"**{}?/-1ex/;?/1ex/^{G\phi}\ar@2"y";"z"**{}?/-1ex/;?/1ex/^{G\phi'} 
        }
        \ar@3"A";"B"**{}?/-1ex/;?/1ex/|{}="x"\save\POS"x"+(0,13)*!C\labelbox{(\alpha_y\#_0F\phi')\\\#_1\alpha_\phi}\ar@{.}"x";c\restore
        \ar@3"B";"C"**{}?/-1ex/;?/1ex/|{}="x"\save\POS"x"+(0,13)*!C\labelbox{\alpha_{\phi'}\\\#_1(G\phi\#_0\alpha_x)}\ar@{.}"x";c\restore
        \ar@3@/_5pc/"A";"C"_{\alpha_{\phi'\#_1\phi}}
      \end{xy}\,,
    \end{equation}
    and for identity 2-cells $\id_f\from{}f\To{}f$ we have an identity
    3-cell
    \begin{equation}
      \label{eq:pstransf2cellid}
      \alpha_{\id_f}=\id_{\alpha_f}\,.
    \end{equation}
    \item The family of 3-cells has to satisfy a kind of cocycle
    condition: For a composable triple $f,f',f''$ of 1-cells
    $\alpha^2$ has to satisfy equation (\ref{eq:pstransf2cocy}).
    \begin{sidewaysfigure}
      \begin{equation}
        \label{eq:pstransf2cocy}
        \begin{xy}
          \xyboxmatrix"A"@+.5cm{
            Fx\ar[r]^{\alpha_x}\ar[d]_{Ff}&Gx\ar[d]^{Gf}\ar@2[dl]**{}?/-1ex/;?/1ex/^{\alpha_f}\\
            Fy\ar[r]^{\alpha_y}\ar[d]_{Ff'}&Gy\ar[d]^{Gf'}\ar@2[dl]**{}?/-1ex/;?/1ex/^{\alpha_{f'}}\\
            Fz\ar[r]^{\alpha_z}\ar[d]_{Ff''}&Gz\ar[d]^{Gf''}\ar@2[dl]**{}?/-1ex/;?/1ex/^{\alpha_{f''}}\\
            Fw\ar[r]_{\alpha_w}&Gw
          }
          \POS;(83,0)
          \xyboxmatrix"B"@+.5cm{
            Fx\ar[r]^{\alpha_x}\ar[d]|/-1ex/{F(f'\#_0f)}&Gx\ar[d]^{G(f'\#_0f)}\ar@2[dl]**{}?/-1ex/;?/1ex/^{\alpha_{f'\#_0f}}\\
            Fz\ar[r]^{\alpha_z}\ar[d]_{Ff''}&Gz\ar[d]^{Gf''}\ar@2[dl]**{}?/-1ex/;?/1ex/^{\alpha_{f''}}\\
            Fw\ar[r]_{\alpha_w}&Gw
          }
          \POS;(0,-83)
          \xyboxmatrix"C"@+.5cm{
            Fx\ar[r]^{\alpha_x}\ar[d]_{Ff}&Gx\ar[d]^{Gf}\ar@2[dl]**{}?/-1ex/;?/1ex/^{\alpha_f}\\
            Fy\ar[r]^{\alpha_y}\ar[d]_{F(f''\#_0f')}&Gy\ar[d]^{G(f''\#_0f')}\ar@2[dl]**{}?/-1ex/;?/1ex/^{\alpha_{f''\#_0f'}}\\
            Fw\ar[r]_{\alpha_w}&Gw\\
          }
          \POS;(83,-83)
          \xyboxmatrix"D"@+.5cm{
            Fx\ar[r]^{\alpha_x}\ar[dd]|/-2ex/{F(f''\#_0f'\#_0f)}&Gx\ar[dd]|/-2ex/{G(f''\#_0f'\#_0f)}\ar@2@<1ex>[ddl]**{}?/-1ex/;?/1ex/|/3ex/{\alpha_{f''\#_0f'\#_0f}}\\
            {}&{}\\
            Fw\ar[r]_{\alpha_w}&Gw
          }
          \POS\ar@3"A";"B"**{}?/-1ex/;?/1ex/^{(\alpha_{f''}\#_0F(f'\#_0f))\\\#_1(Gf''\#_0\underline{\alpha^2_{f',f}})}
          \POS\ar@3"A";"C"**{}?/-1ex/;?/1ex/_{(\underline{\alpha^2_{f'',f'}}\#_0Ff)\\\#_1(G(f''\#_0f')\#_0\alpha_{f})}
          \POS\ar@3"C";"D"**{}?/-1ex/;?/1ex/_{\underline{\alpha^2_{f''\#_0f',f}}}
          \POS\ar@3"B";"D"**{}?/-1ex/;?/1ex/^{\underline{\alpha^2_{f'',f'\#_0f}}}
        \end{xy}
      \end{equation}
    \end{sidewaysfigure}
    furthermore, $\alpha^2$ has to satisfy the normalization
    condition:
    \begin{equation}
      \label{eq:pstransf2cocnorm}
      \alpha^2_{f',f}=
      \begin{cases}
        \id_{\alpha_{f}}&\text{if }f'=\id_y\\
        \id_{\alpha_{f'}}&\text{if }f=\id_x\\
      \end{cases}
    \end{equation}
    \item The family of 3-cells $\alpha^2$ has to be compatible with left and
    right whiskering according to (\ref{eq:pstransf12whiskleft}) and
    (\ref{eq:pstransf12whiskright}).
    \begin{sidewaysfigure}
      \begin{equation} 
        \label{eq:pstransf12whiskleft}
        \begin{xy}
          (60,0) : (0,0)
          \xyboxmatrix"A"{
            Fx\ar[r]^{\alpha_x}\ar[d]_{Ff}& Gx\ar[d]^{Gf}\ar@2[dl]**{}?/-1ex/;?/1ex/^{\alpha_f}\\
            Fy\ar[r]|{\alpha_y}\ar[d]|{Fg}="x"\ar@/_2.5pc/[d]|{Fg'}="y"&Gy\ar[d]^{Gg'}\ar@2[dl]**{}?/-1ex/;?/1ex/^{\alpha_g}\\
            Fz\ar[r]_{\alpha_z}&Gz\ar@2"x";"y"**{}?/-1ex/;?/1ex/^{F\gamma}
          }
          +(1,0)
          \xyboxmatrix"B"{
            Fx\ar[r]^{\alpha_x}\ar[d]_{Ff}&Gx\ar[d]^{Gf}\ar@2[dl]**{}?/-1ex/;?/1ex/^{\alpha_f}="X"\\
            Fy\ar[r]|{\alpha_y}\ar[d]_{Fg'}&Gy\ar[d]|{Gg''}="y"\ar@2[dl]**{}?/-1ex/;?/1ex/^{\alpha_{g'}}\ar@/^2.5pc/[d]|{Gg'}="x"\\
            Fz\ar[r]_{\alpha_z}&Gz\ar@2"x";"y"**{}?/-1ex/;?/1ex/^{G\gamma}="Y"
            \tria"X";"Y"
          }
          +(1,0)
          \xyboxmatrix"C"{
            Fx\ar[r]^{\alpha_x}\ar[d]_{Ff}&Gx\ar[d]^{Gf}\ar@2[dl]**{}?/-1ex/;?/1ex/^{\alpha_f}="X"\\
            Fy\ar[r]|{\alpha_y}\ar[d]_{Fg'}&Gy\ar[d]|{Gg''}="y"\ar@2[dl]**{}?/-1ex/;?/1ex/^{\alpha_{g'}}\ar@/^2.5pc/[d]|{Gg'}="x"\\
            Fz\ar[r]_{\alpha_z}&Gz\ar@2"x";"y"**{}?/-1ex/;?/1ex/^{G\gamma}="Y"
            \tria"Y";"X"
          }
          +(1,0)
          \xyboxmatrix"D"{
            Fx\ar[r]^{\alpha_x}\ar[dd]_{F(g'\#_0f)}&Gx\ar[dd]|(.6){G(g'\#_0f)}="y"\ar@/^3pc/[dd]^{G(g\#_0f)}="x"\ar@2[ddl]**{}?/-1ex/;?/1ex/_(-1){\alpha_{g'\#_0f}}\ar@2"x";"y"**{}?/-1ex/;?/1ex/_{G(\gamma\#_0f)}\\
            &\\
            Fz\ar[r]_{\alpha_z}&Gz
          }
          ,(0,-1)
          \xyboxmatrix"E"{
            Fx\ar[r]^{\alpha_x}\ar[d]_{Ff}\ar@/_4pc/[dd]|{F(g'\#_0f)}="x"&Gx\ar[d]^{Gf}\ar@2[dl]**{}?/-1ex/;?/1ex/^{\alpha_f}\\
            Fy\ar[r]|{\alpha_y}\ar[d]_{Fg}&Gy\ar[d]^{Gg'}\ar@2[dl]**{}?/-1ex/;?/1ex/^ {\alpha_g}\\
            Fz\ar[r]_{\alpha_z}&Gz\ar@2"2,1";"x"**{}?/-1ex/;?/1ex/^{F(\gamma\#_0f)}
          }
          ;p+(3,0)="G"**{}?(.5)
          \xyboxmatrix"F"{
            Fx\ar[r]^{\alpha_x}\ar[dd]|{F(g'\#_0f)}="x"\ar@/_4pc/[dd]_{F(g'\#_0f)}="y"\ar@2"x";"y"**{}?/-1ex/;?/1ex/_{F(\gamma\#_0f)}&Gx\ar[dd]|(.6){G(g\#_0f)}="w"\ar@2[ddl]**{}?/-1ex/;?/1ex/_(-1){\alpha_{g\#_0f}}\\
            &\\
            Fz\ar[r]_{\alpha_z}&Gz
          }
          \ar@3"A";"B"**{}?/-1ex/;?/1ex/^{}="here"\ar@{.}"here";p+(0,.4)*\labelbox{
            (\underline{\alpha_\gamma}\#_0Ff)\\
            \#_1(Gg'\#_0\alpha_{f})
          }
          \ar@3"B";"C"**{}?/-1ex/;?/1ex/^{}="here"\ar@{.}"here";p+(0,.4)*\labelbox{
            (\alpha_{g'}\#_0Ff)\\
            \#_1(\underline{G\gamma\ten\alpha_f})
          }
          \ar@3"C";"D"**{}?/-1ex/;?/1ex/^{}="here"\ar@{.}"here";p+(0,.4)*\labelbox{
            \underline{\alpha^2_{g',f}}\\
            \#_1(G(\gamma\#_0f)\#_0\alpha_x)
          }
          \ar@3"E";"F"**{}?/-1ex/;?/1ex/^{}="here"\ar@{.}"here";p+(0,-.4)*\labelbox{
            (\alpha_z\#_0F(\gamma\#_0f))\\
            \#_1\underline{\alpha^2_{g,f}}
          }
          \ar@3"F";"D"**{}?/-1ex/;?/1ex/^{}="here"\ar@{.}"here";p+(0,-.4)*\labelbox{
            \underline{\alpha_{\gamma\#_0f}}            
          }
          \ar@{=}"A";"E"**{}?/-1ex/;?/1ex/
        \end{xy}
      \end{equation}
      \begin{center}
        Compatibility of the cocycle $\alpha^2$ with left whiskers $\gamma\#_0f$.
      \end{center}
    \end{sidewaysfigure}

    \begin{sidewaysfigure}
      \begin{equation} 
        \label{eq:pstransf12whiskright}
        \begin{xy}
          (60,0) : (0,0)
          \xyboxmatrix"A"{
            Fx\ar[r]^{\alpha_x}\ar[d]_{Ff}="x"\ar@/_2.5pc/[d]|{Ff'}="y"\ar@2"x";"y"**{}?/-1ex/;?/1ex/^{F\delta}="X"&Gx\ar[d]^{Gf}\ar@2[dl]**{}?/-1ex/;?/1ex/^{\alpha_f}\\
            Fy\ar[r]|{\alpha_y}\ar[d]_{Fg}&Gy\ar[d]^{Gg}\ar@2[dl]**{}?/-1ex/;?/1ex/^{\alpha_g}="Y"\\
            Fz\ar[r]_{\alpha_z}&Gz
            \tria"Y";"X"
          }
          +(1,0)
          \xyboxmatrix"B"{
            Fx\ar[r]^{\alpha_x}\ar[d]_{Ff}="x"\ar@/_2.5pc/[d]|{Ff'}="y"\ar@2"x";"y"**{}?/-1ex/;?/1ex/^{F\delta}="X"&Gx\ar[d]^{Gf}\ar@2[dl]**{}?/-1ex/;?/1ex/^{\alpha_f}\\
            Fy\ar[r]|{\alpha_y}\ar[d]|{Fg}&Gy\ar[d]^{Gg}\ar@2[dl]**{}?/-1ex/;?/1ex/^{\alpha_g}="Y"\\
            Fz\ar[r]_{\alpha_z}&Gz
            \tria"X";"Y"
          }
          +(1,0)
          \xyboxmatrix"C"{
            Fx\ar[r]^{\alpha_x}\ar[d]_{Ff'}&Gx\ar[d]^{Gf'}="y"\ar@2[dl]**{}?/-1ex/;?/1ex/^{\alpha_{f'}}\ar@/^2.5pc/[d]|{Gf}="x"\ar@2"x";"y"**{} ?/-1ex/;?/1ex/^{G\delta}\\
            Fy\ar[r]|{\alpha_y}\ar[d]_{Fg}&Gy\ar[d]^{Gg}\ar@2[dl]**{}?/-1ex/;?/1ex/^{\alpha_{g}}\\
            Fz\ar[r]_{\alpha_z}&Gz
          }
          +(1,0)
          \xyboxmatrix"D"{
            Fx\ar[r]^{\alpha_x}\ar[dd]_{F(g\#_0f')}&Gx\ar[dd]|(.6){G(g\#_0f')}="x"\ar@/^2.5pc/[dd]^{G(g\#_0f)}="y"\ar@2[ddl]**{}?/-1ex/;?/1ex/_(-1){\alpha_{g\#_0f'}}\ar@2"x";"y"**{}?/-1ex/;?/1ex/^{G(g\#_0\delta)}
            \\
            {}&\\
            Fz\ar[r]_{\alpha_z}&Gz
          }
          ,(0,-1)
          \xyboxmatrix"E"{
            Fx\ar[r]^{\alpha_x}\ar[d]_{Ff}\ar@/_4pc/[dd]|{F(g\#_0f)}="y"&Gx\ar[d]^{Gf}\ar@2[dl]**{}?/-1ex/;?/1ex/^{\alpha_f}\\
            Fy\ar[r]|{\alpha_y}\ar[d]_{Fg}="x"
            &Gy\ar[d]^{Gg}\ar@2[dl]**{}?/-1ex/;?/1ex/^{\alpha_g}\\
            Fz\ar[r]_{\alpha_z}&Gz\ar@2"2,1";"y" **{}?/-1ex/;?/1ex/_{F(g\#_0\delta)}
          }
          ;p+(3,0)="G"**{}?(.5)
          \xyboxmatrix"F"{
            Fx\ar[r]^{\alpha_x}\ar[dd]|{F(g\#_0f)}="x"\ar@/_4pc/[dd]|{F(g\#_0f')}="y"&Gx\ar[dd]^{G(g\#_0f)}="w"\ar@2[ddl]**{}?/-1ex/;?/1ex/_(-1){\alpha_{g\#_0f}}\\
            &\\
            Fz\ar[r]_{\alpha_z}&Gz\POS\ar@2"x";"y" **{} ?<;?>_{F(g\#_0\delta)}="Y"
          }
          ,"G"
          \ar@3"A";"B"**{}?/-1ex/;?/1ex/^{}="here"\ar@{.}"here";p+(0,.4)*\labelbox{
            (\alpha_g\ten{}F\delta)^{-1}\\
            \#_1(\underline{Gg\#_0\alpha_f})
          }
          \ar@3"B";"C"**{}?/-1ex/;?/1ex/^{}="here"\ar@{.}"here";p+(0,.4)*\labelbox{
            (\alpha_g\#_0Ff')\\
            \#_1(Gg\#_0\underline{\alpha_\delta})
          }
          \ar@3"C";"D"**{}?/-1ex/;?/1ex/^{}="here"\ar@{.}"here";p+(0,.4)*\labelbox{
            \underline{\alpha^2_{g,f'}}\\
            \#_1(G(g\#_0\delta)\#_0\alpha_x)
          }
          \ar@3"E";"F"**{}?/-1ex/;?/1ex/^{}="here"\ar@{.}"here";p+(0,-.4)*\labelbox{
            (\alpha_z\#_0F(g\#_0\delta))\\
            \#_1(\underline{\alpha^2_{g,f}})
          }
          \ar@3"F";"D"**{}?/-1ex/;?/1ex/^{}="here"\ar@{.}"here";p+(0,-.4)*\labelbox{
            \underline{\alpha_{g\#_0\delta}}\\
          }
          \ar@{=}"A";"E"**{}?/-1ex/;?/1ex/
        \end{xy}
      \end{equation}
      \begin{center}
        Compatibility of the cocycle $\alpha^2$ with right whiskers $g\#_0\delta$.
      \end{center}
    \end{sidewaysfigure}
  \end{enumerate}
\end{defn}

\begin{defn}
  \label{def:psmod}
  A \defterm{pseudo-modification} $A\from\alpha\To\alpha'$ between
  pseudo-modifications in the sense of definition \ref{def:pstransf}
  is given by the following data:
  \begin{enumerate}
    \item For every 0-cell $x$ in $\G$ a 2-cell 
    \begin{equation}
      \label{eq:modifex0cell}
      \begin{xy}
        \xyboxmatrix{
          Fx\ar@/^1.5pc/[r]^{\alpha_x}="a"\ar@/_1.5pc/[r]_{\beta_x}="b"&Gx
          \ar@2"a";"b"**{}?/-1ex/;?/1ex/^{A_x}
        }
      \end{xy}
    \end{equation}
    \item For every 1-cell $f\from{}x\to{}y$ a 3-cell in $\H$
    \begin{equation}
      \label{eq:modifex1cell}
      \begin{xy}
        \xyboxmatrix"A"{
          Fx\ar@/^1.5pc/[r]^{\alpha_x}="a"\ar@/_1.5pc/[r]_{\beta_x}="b"\ar[d]_{Ff}&
          Gx\ar[d]^{Gf}\ar@<+1em>@2[dl]**{}?(.5);?(.9)|{\beta_f}\\
          Fy\ar@/_1.5pc/[r]_{\beta_y}="b1"&
          Gy
          \ar@2"a";"b"**{}?/-1ex/;?/1ex/^{A_x}
        }
        +(40,0)
        \xyboxmatrix"B"{
          Fx\ar@/^1.5pc/[r]^{\alpha_x}="a"\ar[d]_{Ff}&
          Gx\ar[d]^{Gf}\ar@<-1em>@2[dl]**{}?(.1);?(.5)|{\alpha_f}\\
          Fy\ar@/^1.5pc/[r]^{\alpha_y}="a1"\ar@/_1.5pc/[r]_{\beta_y}="b1"&
          Gy
          \ar@2"a1";"b1"**{}?/-1ex/;?/1ex/^{A_y}
        }
        \ar@3"A";"B"**{}?/-1ex/;?/1ex/^{A_f}
      \end{xy}
    \end{equation}
  \end{enumerate}
  such that the following conditions hold:
  \begin{enumerate}
    \item Units are preserved:
    \begin{equation}
      \label{eq:modifex0cellunit}
      A_{\id_{x}}=\id_{A_x}
    \end{equation}
    \item Compatibility with the cocycles of $\alpha$ and $\beta$
    according to \eqref{eq:modifexcoccomp}
    \begin{sidewaysfigure}
      \begin{equation}
        \label{eq:modifexcoccomp}
        \begin{xy}
          \save (55,0):(0,-1):: 
          ,(0,0)="A11" 
          ,(1,0)="A21" 
          ,(2,0)="A31" 
          ,(3,0)="A41" 
          ,(0,1)="A12"
          ,(3,1)="A32" 
          ,"A12";"A32"**{};?(.5)="A22" 
          ,(0,2)="A13" 
          ,(1,2)="A23" 
          ,(2,2)="A33"
          ,(3,2)="A43"
          ,(0,3)="A14" 
          ,(1,3)="A24" 
          ,(2,3)="A34" 
          ,(3,3)="A44" 
          \restore 
          \POS,"A11"
          \xyboxmatrix"A11"{
            Fx\ar@/^1.5pc/[r]^{\alpha_x}="a"\ar@/_1.5pc/[r]|{\beta_x}="b"\ar[d]_{Ff}&Gx\ar[d]^{Gf}\ar@<+1em>@2[dl]**{}?(.5);?(.9)|{\beta_f}\\
            Fy\ar@/_1.5pc/[r]|{\beta_y}="b1"\ar[d]_{Ff'}&Gy\ar[d]^{Gf'}\ar@<+1em>@2[dl]**{}?(.5);?(.9)|{\beta_{f'}}\\
            Fz\ar@/_1.5pc/[r]_{\beta_z}="b2"&Gz\ar@2"a";"b"**{}?(.2);?(.8)|{A_x}
          },"A12"
          \xyboxmatrix"A12"{
            Fx\ar@/^1.5pc/[r]^{\alpha_x}="a"\ar@/_1.5pc/[r]|{\beta_x}="b"\ar[dd]|{F(f'\#_0f)}&Gx\ar[dd]|{G(f'\#_0f)}="z"\ar@<+1em>@2[ddl]**{}?(.5);?(.9)|{\beta_{f'\#_0f}}\\
            {}&{}\\
            Fz\ar@/_1.5pc/[r]_{\beta_z}="b2"&Gz\ar@2"a";"b"**{}?(.2);?(.8)|{A_x}="K"
          }
          ,"A22"
          ,"A21"
          \xyboxmatrix"A21"{
            Fx\ar@/^1.5pc/[r]^{\alpha_x}="a"\ar[d]_{Ff}&Gx\ar[d]^{Gf}\ar@<-1em>@2[dl]**{}?(.1);?(.5)|{\alpha_f}\\
            Fy\ar@/_1.5pc/[r]|{\beta_y}="b1"\ar@/^1.5pc/[r]|{\alpha_y}="a1"\ar[d]_{Ff'}&Gy\ar[d]^{Gf'}\ar@<+1em>@2[dl]**{}?(.5);?(.9)|{\beta_{f'}}\\
            Fz\ar@/_1.5pc/[r]_{\beta_z}&Gz\ar@2"a1";"b1"**{}?(.2);?(.8)|{A_y}
          }
          ,"A31"
          \xyboxmatrix"A31"{
            Fx\ar@/^1.5pc/[r]^{\alpha_x}="a"\ar[d]_{Ff}&Gx\ar[d]^{Gf}\ar@<-1em>@2[dl]**{}?(.1);?(.5)|{\alpha_f}\\
            Fy\ar@/^1.5pc/[r]|{\alpha_y}="a1"\ar[d]_{Ff'}&Gy\ar[d]^{Gf'}\ar@<-1em>@2[dl]**{}?(.1);?(.5)|{\alpha_{f'}}\\
            Fz\ar@/^1.5pc/[r]|{\alpha_z}="a2"\ar@/_1.5pc/[r]_{\beta_z}="b2"&Gz\ar@2"a2";"b2"**{}?(.2);?(.8)|{A_z}="Y"
          }
          ,"A41"
          ,"A32"
          \xyboxmatrix"A32"{
            Fx\ar@/^1.5pc/[r]^{\alpha_x}="a"\ar[dd]_{F(f'\#_0f)}&Gx\ar[dd]|{G(f'\#_0f)}="z"\ar@<-1em>@2[ddl]**{}?(.1);?(.5)|{\alpha_{f'\#_0f}}\\
            {}&{}\\
            Fz\ar@/^1.5pc/[r]|{\alpha_z}="a2"\ar@/_1.5pc/[r]_{\beta_z}="b2"&Gz&\ar@2"a2";"b2"**{}?(.2);?(.8)|{A_z}
          }
          \ar@3"A11";"A21"**{}?/-1ex/;?/1ex/^{(\beta_{f'}\#_0Ff)\\\#_1(Gf'\#_0\underline{A_f})}
          \ar@3"A21";"A31"**{}?/-1ex/;?/1ex/^{(\underline{A_{f'}}\#_0Ff)\\\#_1(Gf'\#_0\alpha_f)}
          \ar@3"A31";"A32"**{}?/-1ex/;?/1ex/^{(A_z\#_0F(f'\#_0f))\\\#_1(\underline{\alpha^2_{f',f}})}
          \ar@3"A11";"A12"**{}?/-1ex/;?/1ex/_{(\underline{\beta^2_{f',f}})\\\#_1(G(f'\#_0f)\#_0A_x)}
          \ar@3"A12";"A32"**{}?/-1ex/;?/1ex/_{\underline{A_{f'\#_0f}}}
        \end{xy}
      \end{equation}
      \begin{center}
        Compatibility of the modification $A$ with the cocycles of
        $\alpha$ and $\beta$
      \end{center}
    \end{sidewaysfigure}
    \item For 2-cells $g\from{}f\To{}f'$ in $\G$ the images under $F$
    and $G$ as well the data of $A$, $\alpha$ and $\beta$ are
    compatible as shown in (\ref{eq:modifex2cell})
    \begin{sidewaysfigure}
      \begin{equation}
        \label{eq:modifex2cell}
        \begin{xy}
          \save (55,0):(0,-1):: 
          ,(0,0)="A11" 
          ,(1,0)="A21" 
          ,(2,0)="A31" 
          ,(3,0)="A41" 
          ,(0,1)="A12"
          ,(1,1)="A22"
          ,(2,1)="A32" 
          ,(3,1)="A42" 
          ,(0,2)="A13" 
          ,(1,2)="A23" 
          ,(2,2)="A33"
          ,(3,2)="A43"
          ,(0,3)="A14" 
          ,(1,3)="A24" 
          ,(2,3)="A34" 
          ,(3,3)="A44" 
          \restore 
          \POS,"A11"
          \xyboxmatrix"A11"{
            Fx\ar@/^1.5pc/[r]^{\alpha_x}="a"\ar@/_1.5pc/[r]_{\beta_x}="b"
            \ar[d]|{Ff}="a2"\ar@/_3pc/[d]_{Ff'}="b2"&
            Gx\ar[d]^{Gf}\ar@<+1em>@2[dl]**{}?(.5);?(.9)|{\beta_f}\\
            Fy\ar@/_1.5pc/[r]_{\beta_y}="b1"&
            Gy
            \ar@2"a";"b"**{}?(.2);?(.8)|{A_x}="Y"
            \ar@2"a2";"b2"**{}?(.2);?(.8)|{Fg}="X"
          }
          ,"A21"
          \xyboxmatrix"A21"{
            Fx\ar@/^1.5pc/[r]^{\alpha_x}="a"
            \ar[d]|{Ff}="a2"\ar@/_3pc/[d]_{Ff'}="b2"&
            Gx\ar[d]^{Gf}\ar@<-1em>@2[dl]**{}?(.1);?(.5)|{\alpha_f}\\
            Fy\ar@/^1.5pc/[r]^{\alpha_y}="a1"\ar@/_1.5pc/[r]_{\beta_y}="b1"&
            Gy
            \ar@2"a1";"b1"**{}?(.2);?(.8)|{A_y}="Y"
            \ar@2"a2";"b2"**{}?(.2);?(.8)|{Fg}="X"
            \tria"Y";"X"
          }
          ,"A31"
          \xyboxmatrix"A31"{
            Fx\ar@/^1.5pc/[r]^{\alpha_x}="a"
            \ar[d]|{Ff}="a2"\ar@/_3pc/[d]_{Ff'}="b2"&
            Gx\ar[d]^{Gf}\ar@<-1em>@2[dl]**{}?(.1);?(.5)|{\alpha_f}\\
            Fy\ar@/^1.5pc/[r]^{\alpha_y}="a1"\ar@/_1.5pc/[r]_{\beta_y}="b1"&
            Gy
            \ar@2"a1";"b1"**{}?(.2);?(.8)|{A_y}="Y"
            \ar@2"a2";"b2"**{}?(.2);?(.8)|{Fg}="X"
            \tria"X";"Y"
          }
          ,"A41"
          \xyboxmatrix"A41"{
            Fx\ar@/^1.5pc/[r]^{\alpha_x}="a"\ar[d]_{Ff'}&
            Gx\ar[d]|{Gf'}="a2"\ar@<-1em>@2[dl]**{}?(.1);?(.5)|{\alpha_{f'}}
            \ar@/^3pc/[d]^{Gf}="b2"\\
            Fy\ar@/^1.5pc/[r]^{\alpha_y}="a1"\ar@/_1.5pc/[r]_{\beta_y}="b1"&
            Gy
            \ar@2"a1";"b1"**{}?(.2);?(.8)|{A_y}="Y"
            \ar@2"b2";"a2"**{}?(.2);?(.8)|{Gg}="X"
          },"A12"
          \xyboxmatrix"A12"{
            Fx\ar@/^1.5pc/[r]^{\alpha_x}="a"\ar@/_1.5pc/[r]_{\beta_x}="b"
            \ar[d]|{Ff}="a2"\ar@/_3pc/[d]_{Ff'}="b2"&
            Gx\ar[d]^{Gf}\ar@<+1em>@2[dl]**{}?(.5);?(.9)|{\beta_{f'}}\\
            Fy\ar@/_1.5pc/[r]_{\beta_y}="b1"&
            Gy
            \ar@2"a";"b"**{}?(.2);?(.8)|{A_x}="Y"
            \ar@2"a2";"b2"**{}?(.2);?(.8)|{Fg}="X"
          }
          ,"A22"
          \xyboxmatrix"A22"{
            Fx\ar@/^1.5pc/[r]^{\alpha_x}="a"\ar@/_1.5pc/[r]_{\beta_y}="b"
            \ar[d]_{Ff'}&
            Gx\ar[d]|{Gf'}="a2"\ar@/^3pc/[d]^{Gf}="b2"
            \ar@<+1em>@2[dl]**{}?(.6);?(.9)|{\beta_f}\\
            Fy\ar@/_1.5pc/[r]_{\beta_y}="a1"&
            Gy
            \ar@2"a";"b"**{}?(.2);?(.8)|{A_x}="Y"
            \ar@2"b2";"a2"**{}?(.2);?(.8)|{Gg}="X"
            \tria"Y";"X"
          }
          ,"A32"
          \xyboxmatrix"A32"{
            Fx\ar@/^1.5pc/[r]^{\alpha_x}="a"\ar@/_1.5pc/[r]_{\beta_y}="b"
            \ar[d]_{Ff'}&
            Gx\ar[d]|{Gf'}="a2"\ar@/^3pc/[d]^{Gf}="b2"
            \ar@<+1em>@2[dl]**{}?(.6);?(.9)|{\beta_{f'}}\\
            Fy\ar@/_1.5pc/[r]_{\beta_y}="a1"&
            Gy
            \ar@2"a";"b"**{}?(.2);?(.8)|{A_x}="Y"
            \ar@2"b2";"a2"**{}?(.2);?(.8)|{Gg}="X"
            \tria"X";"Y"
          }
          ,"A42"
          \xyboxmatrix"A42"{
            Fx\ar@/^1.5pc/[r]^{\alpha_x}="a"\ar[d]_{Ff'}&
            Gx\ar[d]|{Gf'}="a2"\ar@<-1em>@2[dl]**{}?(.1);?(.5)|{\alpha_{f'}}
            \ar@/^3pc/[d]^{Gf}="b2"\\
            Fy\ar@/^1.5pc/[r]^{\alpha_y}="a1"\ar@/_1.5pc/[r]_{\beta_y}="b1"&
            Gy
            \ar@2"a1";"b1"**{}?(.2);?(.8)|{A_y}="Y"
            \ar@2"b2";"a2"**{}?(.2);?(.8)|{Gg}="X"
          }
          \ar@{=}"A11";"A12"**{}?/-1ex/;?/1ex/
          \ar@{=}"A41";"A42"**{}?/-1ex/;?/1ex/
          \ar@3"A11";"A21"**{}?/-1ex/;?/1ex/^{(\beta_y\#_0Fg)\\\#_1\underline{A_f}}
          \ar@3"A21";"A31"**{}?/-1ex/;?/1ex/^{\overline{A_y\ten{}Fg}\\\#_1\alpha_f}
          \ar@3"A31";"A41"**{}?/-1ex/;?/1ex/^{(A_y\#_0Ff')\\\#_1\underline{\alpha_{g}}}
          \ar@3"A12";"A22"**{}?/-1ex/;?/1ex/_{\underline{\beta_g}\\\#_1(Gf\#_0A_x)}
          \ar@3"A22";"A32"**{}?/-1ex/;?/1ex/_{\beta_{f'}\\\#_1(Gg\ten{}A_x)}
          \ar@3"A32";"A42"**{}?/-1ex/;?/1ex/_{\underline{A_{f'}}\\\#_1(Gg\#_0\alpha_x)}
        \end{xy}
      \end{equation}
      \begin{center}
        Compatibility of 2-cells with $A$, $\alpha$ and $\beta$
      \end{center}
    \end{sidewaysfigure}
  \end{enumerate}
\end{defn}

\begin{defn}
  \label{def:perturb}
  A \defterm{perturbation} $\Gamma\from{}A\Tto{}A'$ between
  pseudo-modifications in the sense definition \ref{def:psmod} is
  given by a family 3-cells indexed by 0-cells $x$ of $\G$:
  \begin{equation}
    \label{eq:pertex0cell}
    \begin{xy}
      \save (40,0):(0,-1):: 
      ,(0,0)="A11" 
      ,(1,0)="A21" 
      \restore
      \POS,"A11"
      \xyboxmatrix"A11"{
        Fx\ar@/^1.5pc/[r]^{\alpha_x}="a"
        \ar@/_1.5pc/[r]_{\beta_x}="b"& Gx
        \ar@2"a";"b"**{}?(.2);?(.8)|{A_x}
      }
      ,"A21"
      \xyboxmatrix"A21"{
        Fx\ar@/^1.5pc/[r]^{\alpha_x}="a"
        \ar@/_1.5pc/[r]_{\beta_x}="b"& Gx
        \ar@2"a";"b"**{}?(.2);?(.8)|{B_x}
      }
      \ar@3"A11";"A21"**{}?<;?>^{\Gamma_x}
    \end{xy}
  \end{equation}
  such that
  \begin{equation}
    \label{eq:pertex1cell}
    \begin{xy}
          \save (45,0):(0,-1):: 
          ,(0,0)="A11" 
          ,(1,0)="A21" 
          ,(0,1)="A12"
          ,(1,1)="A22"
          \restore 
          \POS,"A11"
          \xyboxmatrix"A11"{
            Fx\ar@/^1.5pc/[r]^{\alpha_x}="a"\ar@/_1.5pc/[r]_{\beta_x}="b"
            \ar[d]_{Ff}="a2"&
            Gx\ar[d]^{Gf}\ar@<+1em>@2[dl]**{}?(.5);?(.9)|{\beta_f}\\
            Fy\ar@/_1.5pc/[r]_{\beta_y}="b1"&
            Gy
            \ar@2"a";"b"**{}?(.2);?(.8)|{A_x}="Y"
          }
          ,"A12"
          \xyboxmatrix"A12"{
            Fx\ar@/^1.5pc/[r]^{\alpha_x}="a"
            \ar[d]_{Ff}&
            Gx\ar[d]^{Gf}\ar@<-1em>@2[dl]**{}?(.1);?(.5)|{\alpha_f}\\
            Fy\ar@/^1.5pc/[r]^{\alpha_y}="a1"\ar@/_1.5pc/[r]_{\beta_y}="b1"&
            Gy
            \ar@2"a1";"b1"**{}?(.2);?(.8)|{A_y}="Y"
          }   
          ,"A21"
          \xyboxmatrix"A21"{
            Fx\ar@/^1.5pc/[r]^{\alpha_x}="a"\ar@/_1.5pc/[r]_{\beta_x}="b"
            \ar[d]_{Ff}="a2"&
            Gx\ar[d]^{Gf}\ar@<+1em>@2[dl]**{}?(.5);?(.9)|{\beta_f}\\
            Fy\ar@/_1.5pc/[r]_{\beta_y}="b1"&
            Gy
            \ar@2"a";"b"**{}?(.2);?(.8)|{B_x}="Y"
          }
          ,"A22"
          \xyboxmatrix"A22"{
            Fx\ar@/^1.5pc/[r]^{\alpha_x}="a"
            \ar[d]_{Ff}&
            Gx\ar[d]^{Gf}\ar@<-1em>@2[dl]**{}?(.1);?(.5)|{\alpha_f}\\
            Fy\ar@/^1.5pc/[r]^{\alpha_y}="a1"\ar@/_1.5pc/[r]_{\beta_y}="b1"&
            Gy
            \ar@2"a1";"b1"**{}?(.2);?(.8)|{B_y}="Y"
          }
          \ar@3"A11";"A21"**{}?<;?>^{\beta_f\\\#_1(Gf\#_0\Gamma_x)}
          \ar@3"A21";"A22"**{}?<;?>^{B_f}
          \ar@3"A11";"A12"**{}?<;?>_{A_f}
          \ar@3"A12";"A22"**{}?<;?>_{(\Gamma_y\#_0Ff)\\\#_1\alpha_f}      
    \end{xy}
  \end{equation}
  commutes.
\end{defn}

\section{Composites in ${[\G,\H]^p}$}
\label{sec:comp-pseu-transf}

\subsection{Composites of Pseudo-Transformations }
\label{sec:comp-pseudo-transf}

For reference we the data for the composite
$\xymatrix@1{G\ar[r]^{\beta}&G'\ar[r]^{\beta}&G''}$ of
pseudo-transformations. This also apprears \citep[][Appendix
C]{gohla2014}, albeit for pseudo-functors. 
\begin{enumerate}
  \item On 0-cells:
  \begin{equation}
    \label{eq:betabetax}
    (\beta'*_0\beta)_x=\beta'_x\#_0\beta_x\,.
  \end{equation}
  \item On 1-cells:
  \begin{equation}
    \label{eq:betabetaf}
    (\beta'*_0\beta)_f=
    \begin{xy}
      \xyboxmatrix{
        Gx\ar[r]^{\beta_x}\ar[d]_{Gf}&G'x\ar[r]^{\beta'_x}\ar[d]|{G'f}\ar@2[dl]**{}?/-1ex/;?/1ex/^{\beta_f}&G''x\ar[d]^{G''f}\ar@2[dl]**{}?/-1ex/;?/1ex/^{\beta'_f}\\
        Gy\ar[r]_{\beta_y}&G'y\ar[r]_{\beta'_y}&G''
      }
    \end{xy}\,.
  \end{equation}
  \item On 2-cells: See \eqref{eq:betabetaphi}.
  \begin{sidewaysfigure}
    \begin{center}$(\beta'*_0\beta)_\phi$\end{center}
    \begin{equation}
      \label{eq:betabetaphi}
      \begin{xy}
        \xyboxmatrix"A"{
          Gx\ar[r]^{\beta_x}\ar[d]|{Gf}="x"\ar@/_3pc/[d]_{Gf'}="y"\ar@2"x";"y"**{}?/-1ex/;?/1ex/^{G\phi}&G'x\ar[r]^{\beta'_x}\ar[d]|{G'f}\ar@2[dl]**{}?/-1ex/;?/1ex/^{\beta_f}&G''x\ar[d]^{G''f}\ar@2[dl]**{}?/-1ex/;?/1ex/^{\beta'_f}\\
          Gy\ar[r]_{\beta_y}&G'y\ar[r]_{\beta'_y}&G''y
        }
        ;(70,0)
        \xyboxmatrix"B"{
          Gx\ar[r]^{\beta_x}\ar[d]_{Gf'}&G'x\ar[r]^{\beta'_x}\ar@/^1.5pc/[d]|{G'f}="x"\ar@/_1.5pc/[d]|{G'f'}="y"\ar@2"x";"y"**{}?/-1ex/;?/1ex/^{G'\phi}\ar@2[dl]**{}?/1ex/;?/3ex/_{\beta_{f'}}&G''x\ar[d]^{G''f}\ar@2[dl]**{}?/-3ex/;?/-1ex/^{\beta'_f}\\
          Gy\ar[r]_{\beta_y}&G'y\ar[r]_{\beta'_y}&G''y
        }
        ;(140,0)
        \xyboxmatrix"C"{
          Gx\ar[r]^{\beta_x}\ar[d]_{Gf'}&G'x\ar[r]^{\beta'_x}\ar[d]|{G'f'}\ar@2[dl]**{}?/-1ex/;?/1ex/^{\beta_{f'}}&G''x\ar@/^3pc/[d]|{G''f}="x"\ar[d]|{G''f'}="y"\ar@2"x";"y"**{}?/-1ex/;?/1ex/^{G''_\phi}\ar@2[dl]**{}?/-1ex/;?/1ex/^{\beta'_{f'}}\\
          Gy\ar[r]_{\beta_y}&G'y\ar[r]_{\beta'_y}&G''y
        }
        \ar@3"A";"B"**{}?/-1ex/;?/1ex/|{}="here"\ar@{.}"here";p+(0,25)*!C\labelbox{(\beta'_y\#_0\underline{\beta_\phi})\\\#_1(\beta'_f\#_0\beta_x)}
        \ar@3"B";"C"**{}?/-1ex/;?/1ex/|{}="here"\ar@{.}"here";p+(0,25)*!C\labelbox{(\beta'_y\#_0\beta_{f'})\\\#_1(\underline{\beta'_\phi}\#_0\beta_x)}
      \end{xy}\,.
    \end{equation}
  \end{sidewaysfigure}
  \item On composable pairs: see \eqref{eq:betabetatwo}.
  \begin{sidewaysfigure}
    \begin{center}$(\beta'*_0\beta)^2_{f',f}$\end{center}
    \begin{equation}
      \label{eq:betabetatwo}      
      \begin{xy}
        \xyboxmatrix"A"{
          Gx\ar[r]^{\beta_x}\ar[d]_{Gf}&G'x\ar[r]^{\beta'_x}\ar[d]|{G'f}\ar@2[dl]**{}?/-1ex/;?/1ex/_{\beta_f}|{}="x"&G''x\ar[d]^{G''f}\ar@2[dl]**{}?/-1ex/;?/1ex/^{\beta'_f}\\
          Gy\ar[r]|{\beta_y}\ar[d]_{Gf'}&G'y\ar[r]|{\beta'_y}\ar[d]|{G'f'}\ar@2[dl]**{}?/-1ex/;?/1ex/^{\beta_{f'}}&G''y\ar[d]^{G''f'}\ar@2[dl]**{}?/-1ex/;?/1ex/^{\beta'_{f'}}|{}="y"\\
          Gz\ar[r]_{\beta_z}&G'z\ar[r]_{\beta'_z}&G''z \tria"x";"y"
        };(70,0) \xyboxmatrix"B"{
          Gx\ar[r]^{\beta_x}\ar[d]_{Gf}&G'x\ar[r]^{\beta'_x}\ar[d]|{G'f}\ar@2[dl]**{}?/-1ex/;?/1ex/_{\beta_f}|{}="x"&G''x\ar[d]^{G''f}\ar@2[dl]**{}?/-1ex/;?/1ex/^{\beta'_f}\\
          Gy\ar[r]|{\beta_y}\ar[d]_{Gf'}&G'y\ar[r]|{\beta'_y}\ar[d]|{G'f'}\ar@2[dl]**{}?/-1ex/;?/1ex/^{\beta_{f'}}&G''y\ar[d]^{G''f'}\ar@2[dl]**{}?/-1ex/;?/1ex/^{\beta'_{f'}}|{}="y"\\
          Gz\ar[r]_{\beta_z}&G'z\ar[r]_{\beta'_z}&G''z \tria"y";"x"
        };(140,0) \xyboxmatrix"C"{
          Gx\ar[r]^{\beta_x}\ar[dd]_{G(f'\#_0f)}&G'x\ar[r]^{\beta'_x}\ar[dd]|{G'(f'\#_0f)}\ar@2[ddl]**{}?/-3ex/;?/-1ex/_{\beta_{f'\#_0f}}|{}="x"&G''x\ar[dd]^{G''(f'\#_0f)}\ar@2[ddl]**{}?/1ex/;?/3ex/^{\beta'_{f'\#_0f}}\\
          {}&{}&{}\\
          Gz\ar[r]_{\beta_z}&G'z\ar[r]_{\beta'_z}&G''z }
        \ar@3"A";"B"**{}?/-1ex/;?/1ex/|{}="here"\ar@{.}"here";p+(0,25)*!C\labelbox{(\beta'_z\#_0\beta_{f'}\#_0Gf)\\\#_1(\underline{\beta'_{f'}\ten\beta_f})\\\#_1(Gf'\#_0\beta'_f\#_0\beta_x)}
        \ar@3"B";"C"**{}?/-1ex/;?/1ex/|{}="here"\ar@{.}"here";p+(0,25)*!C\labelbox{(\beta'_z\#_0\underline{\beta^2_{f',f}})\\\#_1(\underline{\beta'^2_{f',f}}\#_0\beta_x)}
      \end{xy}
    \end{equation}
  \end{sidewaysfigure}
\end{enumerate}

\subsection{Composites of Pseudo-Modifications}
\label{sec:comp-pseudo-transf-1}

The composite $A'*_1A$ of pseudo-modifications 
\begin{equation*}
  \begin{xy}
    \xyboxmatrix{
      G\ar@/^2pc/[r]^{\alpha}|{}="x"\ar[r]|{\alpha'}="y"\ar@/_2pc/[r]_{\alpha''}|{}="z"\tara{A}"x";"y"\tara{A'}"y";"z"&G'
    }
  \end{xy}
\end{equation*}
is given by
\begin{equation}
  \label{eq:psmodcomp0}
  (A'*_1A)_x=A'_x\#_1A_x
\end{equation}
\begin{equation}
  \label{eq:psmodcomp1}
  (A'*_1A)_f=
  \begin{xy}
    0;(60,0):0
    \xyboxmatrix"A"@+1cm{
      Gx\ar@/^1.5pc/[r]^{\alpha_x}="a"\ar@/_1.5pc/[r]_{\alpha''_x}="b"\ar[d]_{Gf}&
      G'x\ar[d]^{G'f}\tara{\alpha''_f}{;[d]**{}?};[dl]\\
      Gy\ar@/_1.5pc/[r]_{\alpha''_y}="b1"&
      G'y
      \ar@<-2ex>@2"a";"b"**{}?/-1ex/;?/1ex/^{(A'*_1A)_x}
    }
    +(1,0)
    \xyboxmatrix"B"@+1cm{
      Gx\ar@/^2pc/[r]^{\alpha_x}="a"\ar[r]|{\alpha'_x}="b"\ar[d]_{Gf}\tara{A_x}"a";"b"&
      G'x\ar[d]^{G'f}\tara{\alpha'_f};[dl]\\
      Gy\ar[r]|{\alpha'_y}="a1"\ar@/_2pc/[r]_{\alpha''_y}="b1"\tara{A'_y}"a1";"b1"&
      G'y
    }
    +(1,0)
    \xyboxmatrix"C"@+1cm{
      Gx\ar@/^1.5pc/[r]^{\alpha_x}\ar[d]_{Gf}&
      G'x\ar[d]^{G'f}\tara{\alpha_f};{[dl];[l]**{}?}\\
      Gy\ar@/^1.5pc/[r]^{\alpha_y}="a"\ar@/_1.5pc/[r]_{\alpha''_y}="b"&
      G'y
      \ar@<-2ex>@2"a";"b"**{}?/-1ex/;?/1ex/^{(A'*_1A)_y}
    }
    \ar@3"A";"B"**{}?/-1ex/;?/1ex/^{\underline{A'_f}\#_1\\(G'f\#_0A_x)}
    \ar@3"B";"C"**{}?/-1ex/;?/1ex/^{(A'_y\#_0Gf)\\\#_1\underline{A_f}}
  \end{xy}\,.
\end{equation}

\subsection{Whiskers of Pseudo-Modifications}
\label{sec:whisk-pseudo-modif}

For a pseudo-transformation and a pseudo-modification
\begin{equation*}
  \begin{xy}
    \xyboxmatrix{
      G\ar@/^1.5pc/[r]^{\alpha}|{}="x"\ar@/_1.5pc/[r]_{\alpha'}|{}="y"\tara{A}"x";"y"&G'\ar[r]^{\beta}&G''
    }
  \end{xy}
\end{equation*}
we define the right whisker $\beta*_0A$ as follows: 
\begin{equation}
  \label{eq:psmodwhiskrightx}
  (\beta*_0A)_x=\beta_x\#_0A_x
\end{equation}
and \eqref{eq:psmodwhiskrightf}.
\begin{sidewaysfigure}
  \begin{multline}
    \label{eq:psmodwhiskrightf}
    (\beta*_0A)_f=\\
    \begin{xy}
      0;(85,0):0
      \xyboxmatrix"A"@+.5cm{
        Gx\ar@/^1.5pc/[r]^{\alpha_x}="x"\ar@/_1.5pc/[r]_{\alpha'_x}="y"\ar[d]_{Gf}\POS{\tara{A_x}"x";"y"}="X"&G'x\ar[r]^{\beta_x}\ar[d]|{G'f}\tara{\alpha'_f}{[d]**{}?};[dl]&G''x\ar[d]^{G''f}\POS{\tara{\beta_f}[dl]}="Y"\\
        Gy\ar@/_1.5pc/[r]_{\alpha'_y}&G'y\ar[r]_{\beta_y}&G''y \tria"X";"Y"
      }+(1,0)
      \xyboxmatrix"B"@+.5cm{
        Gx\ar@/^1.5pc/[r]^{\alpha_x}="x"\ar@/_1.5pc/[r]_{\alpha'_x}="y"\ar[d]_{Gf}\POS{\tara{A_x}"x";"y"}="X"&G'x\ar[r]^{\beta_x}\ar[d]|{G'f}\tara{\alpha'_f}{[d]**{}?};[dl]&G''x\ar[d]^{G''f}\POS{\tara{\beta_f}[dl]}="Y"\\
        Gy\ar@/_1.5pc/[r]_{\alpha'_y}&G'y\ar[r]_{\beta_y}&G''y \tria"Y";"X"
      }+(1,0)
      \xyboxmatrix"C"@+.5cm{
        Gx\ar@/^1.5pc/[r]^{\alpha_x}\ar[d]_{Gf}&G'x\ar[r]^{\beta_x}\ar[d]|{G'f}\tara{\alpha_f};{[dl];[l]**{}?}&G''x\ar[d]^{G''f}\tara{\beta_f}[dl]\\
        Gy\ar@/^1.5pc/[r]^{\alpha_y}="x"\ar@/_1.5pc/[r]_{\alpha'_y}="y"\tara{A_y}"x";"y"&G'y\ar[r]_{\beta_y}&G''y     
      }
      \ttara{(\beta_y\#_0\alpha'_f)\\\#_1(\underline{\beta_f\ten{}A_x})^{-1}}"A";"B"
      \ttara{(\underline{\beta_y\#_0A_f})\\\#_1(\beta_f\#_0\alpha_x)}"B";"C"
    \end{xy}
  \end{multline}
\end{sidewaysfigure}

For a pseudo-transformation and a pseudo-modification
\begin{equation*}
  \begin{xy}
    \xyboxmatrix{
      G\ar[r]^{\alpha}&G'\ar@/^1.5pc/[r]^{\beta}|{}="x"\ar@/_1.5pc/[r]_{\beta'}|{}="y"\tara{B}"x";"y"&G''
    }
  \end{xy}
\end{equation*}
we define the right whisker $B*_0\alpha$ as follows: 
\begin{equation}
  \label{eq:psmodwhiskleftx}
  (B*_0\alpha)_x=B_x\#_0\alpha_x
\end{equation}
and \eqref{eq:psmodwhiskleftf}.
\begin{sidewaysfigure}
  \begin{multline}
    \label{eq:psmodwhiskleftf}
    (B*_0\alpha)_f=\\
    \begin{xy}
      0;(85,0):0
      \xyboxmatrix"A"@+.5cm{
        Gx\ar[r]^{\alpha_x}\ar[d]_{Gf}&G'x\ar@/^1.5pc/[r]^{\beta_x}="x"\ar@/_1.5pc/[r]_{\beta'_x}="y"\ar[d]|{G'f}\POS{\tara{B_x}"x";"y"}\tara{\alpha_f}[dl]&G''x\ar[d]^{G''f}\POS{\tara{\beta'_f}{;[d]**{}?};[dl]}\\
        Gy\ar[r]_{\alpha_y}&G'y\ar@/_1.5pc/[r]_{\beta'_y}&G''y
      }+(1,0)
      \xyboxmatrix"B"@+.5cm{
        Gx\ar[r]^{\alpha_x}\ar[d]_{Gf}&G'x\ar@/^1.5pc/[r]^{\beta_x}\ar[d]|{G'f}\POS{\tara{\alpha_f}[dl]}="X"&G''x\ar[d]^{G''f}\tarb{\beta_f};{[dl];[l]**{}?}\\
        Gy\ar[r]_{\alpha_y}&G'y\ar@/^1.5pc/[r]^{\beta_y}="x"\ar@/_1.5pc/[r]_{\beta'_y}="y"\POS{\tara{B_y}"x";"y"}="Y"&G''y
        \tria"Y";"X"
      }+(1,0)
      \xyboxmatrix"C"@+.5cm{
        Gx\ar[r]^{\alpha_x}\ar[d]_{Gf}&G'x\ar@/^1.5pc/[r]^{\beta_x}\ar[d]|{G'f}\POS{\tara{\alpha_f}[dl]}="X"&G''x\ar[d]^{G''f}\tarb{\beta_f};{[dl];[l]**{}?}\\
        Gy\ar[r]_{\alpha_y}&G'y\ar@/^1.5pc/[r]^{\beta_y}="x"\ar@/_1.5pc/[r]_{\beta'_y}="y"\POS{\tara{B_y}"x";"y"}="Y"&G''y
        \tria"X";"Y"
      }
      \ttara{(\beta'_y\#_0\alpha_f)\\\#_1(\underline{B_f\#_0\alpha_x})}"A";"B"
      \ttara{(\underline{B_y\ten\alpha_f})^{-1}\\\#_1(\beta'_f\#_0\alpha_x)}"B";"C"
    \end{xy}
  \end{multline}
\end{sidewaysfigure}

\section{Horizontal Composition in $\Gray\Cat_{\fQ^1}$}
\label{sec:horiz-comp-grayc}

We define the composition $\_*_{-1}\_$ of $\Gray$-functors, pseudo-transformation,
pseudo-modification, and perturbations along a $\Gray$-category:

\begin{center}
  \begin{tabular}{l|llll}
    ${*_{-1}}$&$G$&$\alpha$&$A$&$\Gamma$\\\hline
    $H$&$H*_{-1}G$&$H*_{-1}\alpha$&$H*_{-1}A$&$H*_{-1}\Gamma$\\
    $\beta$&$\beta*_{-1}G$&$\beta*_{-1}\alpha$&$\beta*_{-1}A$&$\beta*_{-1}\Gamma$\\
    $B$&$B*_{-1}G$&$B*_{-1}\alpha$&$B*_{-1}A$&$B*_{-1}\Gamma$\\
    $\Delta$&$\Delta*_{-1}G$&$\Delta*_{-1}\alpha$&$\Delta*_{-1}A$&$\Delta*_{-1}\Gamma$\\
  \end{tabular}
\end{center}

\revise{On a 0-cell $x$} each of these is defined as
\begin{center}
  \begin{tabular}{c|cccc}
    ${*_{-1}}$&$G$&$\alpha\from{}G\to{}G'$&$A$&$\Gamma$\\\hline
    $H$&$HGx$&$H\alpha_x$&$HA_x$&$H\Gamma_x$\\
    $\beta\from{}H\to{}H'$&$\beta_{Hx}$&
    $\xymatrix@ru@+.7cm{HGx\ar[r]^{
        \beta_{Gx}}\ar[d]_{H\alpha_x}&H'Gx\ar[d]^{H'\alpha_x}\ar@2[dl]**{}?/-1ex/;?/1ex/^{\beta_{\alpha_x}}\\
      HG'x\ar[r]_{\beta_{H'x}}&H'G'x
    }$&$\beta_{A_x}$&$\id$\\
    $B$&$B_{Gx}$&$B_{\alpha_x}$&$\id$&$\id$\\
    $\Delta$&$\Delta_{Gx}$&$\id$&$\id$&$\id$\\
  \end{tabular}
\end{center}

\revise{On a 1-cell $f$} each of these is defined as
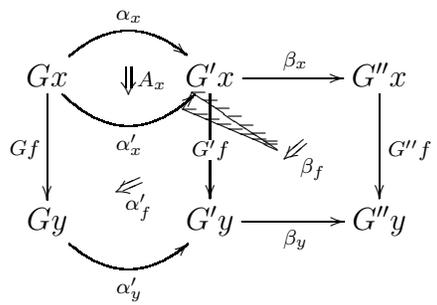
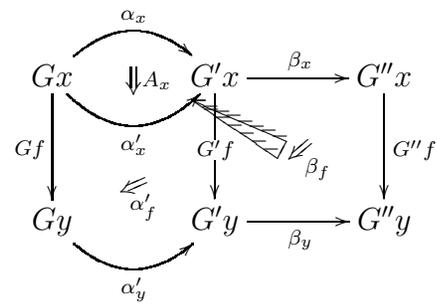
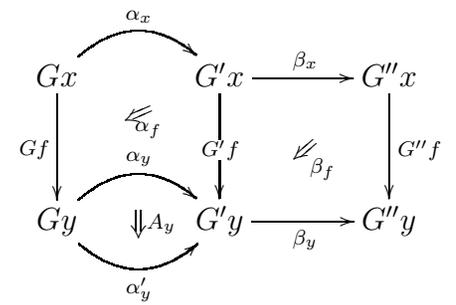
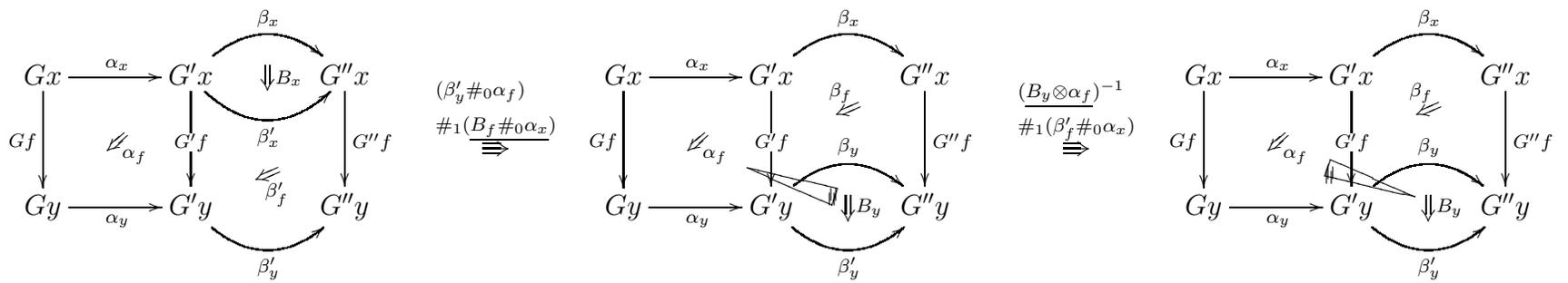
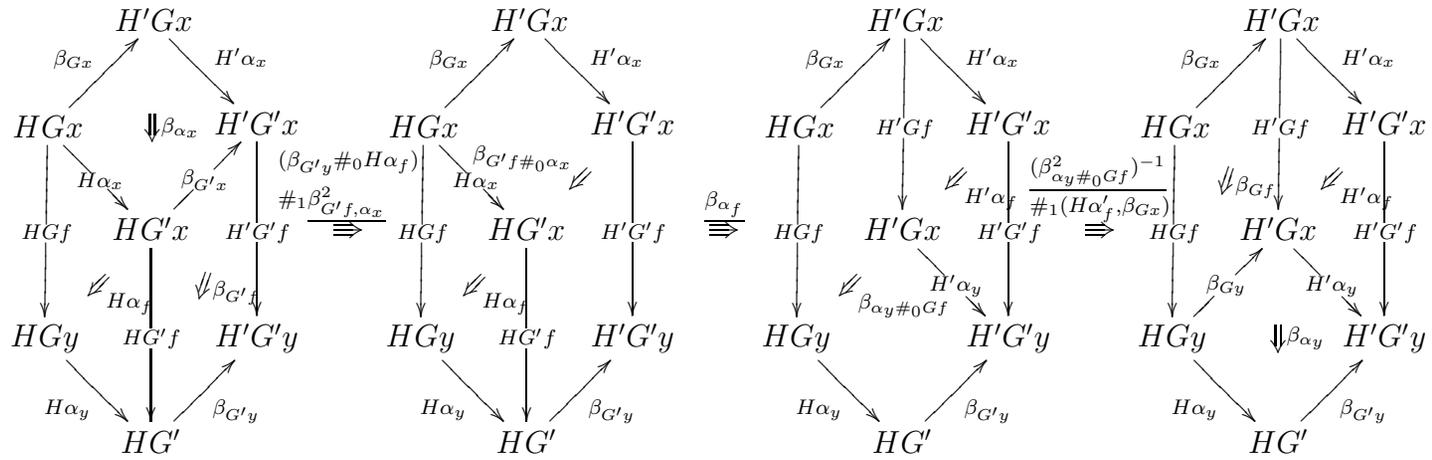
\begin{center}
  \begin{tabular}{c|cccc}
    ${*_{-1}}$&$G$&$\alpha\from{}G\to{}G'$&$A$&$\Gamma$\\\hline
    $H$&$HGf$&$H\alpha_f$&$HA_f$&$\id$\\
    $\beta\from{}H\to{}H'$&$\beta_{Gf}$&$(\beta*_{-1}\alpha)_f$ see \eqref{eq:betaalphaf}&$\id$&$\id$\\
    $B$&$B_{Gf}$&$\id$&$\id$&$\id$\\
    $\Delta$&$\id$&$\id$&$\id$&$\id$\\
  \end{tabular}
\end{center}

\begin{sidewaysfigure}
  \begin{center}$(\beta*_{-1}\alpha)_f$\end{center}
  \begin{multline}
    \label{eq:betaalphaf}
    \begin{xy}
      \xyboxmatrix"A"@ru@+.3cm{
        {}&HGx\ar[r]^{\beta_{Gx}}\ar[d]|{H\alpha_x}\ar[dl]|{HGf}&H'Gx\ar[d]^{H'\alpha_x}\ar@2[dl]**{}?/-1ex/;?/1ex/^{\beta_{\alpha_x}}\\
        HGy\ar[d]_{H\alpha_y}&HG'x\ar[r]|{\beta_{G'x}}\ar[dl]|{HG'f}\ar@2[l]**{}?/-1ex/;?/+1ex/^{H\alpha_f}&H'G'x\ar[dl]|{H'G'f}\ar@2[dll]**{}?/-1ex/;?/+1ex/^{\beta_{G'f}}\\
        HG'\ar[r]_{\beta_{G'y}}&H'G'y&{}
      };(50,0)
      \xyboxmatrix"B"@ru@+.3cm{
        {}&HGx\ar[r]^{\beta_{Gx}}\ar[d]|{H\alpha_x}\ar[dl]|{HGf}&H'Gx\ar[d]^{H'\alpha_x}\\
        HGy\ar[d]_{H\alpha_y}&HG'x\ar[dl]|{HG'f}\ar@2[l]**{}?/-1ex/;?/+1ex/^{H\alpha_f}&H'G'x\ar[dl]|{H'G'f}\ar@2[l]**{}?/-1ex/;?/+1ex/_{\beta_{G'f\#_0\alpha_x}}\\
        HG'\ar[r]_{\beta_{G'y}}&H'G'y&{}
      };(100,0)
      \xyboxmatrix"C"@ru@+.3cm{
        {}&HGx\ar[r]^{\beta_{Gx}}\ar[dl]|{HGf}&H'Gx\ar[d]^{H'\alpha_x}\ar[dl]|{H'Gf}\\
        HGy\ar[d]_{H\alpha_y}&H'Gx\ar[d]|{H'\alpha_y}\ar@2[l]**{}?/-1ex/;?/1ex/^{\beta_{\alpha_y\#_0Gf}}&H'G'x\ar[dl]|{H'G'f}\ar@2[l]**{}?/-1ex/;?/+1ex/^{H'\alpha_f}\\
        HG'\ar[r]_{\beta_{G'y}}&H'G'y&{}
      };(150,0)
      \xyboxmatrix"D"@ru@+.3cm{
        {}&HGx\ar[r]^{\beta_{Gx}}\ar[dl]|{HGf}&H'Gx\ar[d]^{H'\alpha_x}\ar[dl]|{H'Gf}\ar@2[dll]**{}?/-1ex/;?/1ex/^{\beta_{Gf}}\\
        HGy\ar[d]_{H\alpha_y}\ar[r]|{\beta_{Gy}}&H'Gx\ar[d]|{H'\alpha_y}\ar@2[dl]**{}?/-1ex/;?/1ex/^{\beta_{\alpha_y}}&H'G'x\ar[dl]|{H'G'f}\ar@2[l]**{}?/-1ex/;?/+1ex/^{H'\alpha_f}\\
        HG'\ar[r]_{\beta_{G'y}}&H'G'y&{} }
      \ar@3"A";"B"**{}?/-1ex/;?/1ex/^{(\beta_{G'y}\#_0H\alpha_f)\\\#_1\underline{\beta^2_{G'f,\alpha_x}}}
      \ar@3"B";"C"**{}?/-1ex/;?/1ex/^{\underline{\beta_{\alpha_f}}}
      \ar@3"C";"D"**{}?/-1ex/;?/1ex/^{\underline{(\beta^2_{\alpha_y\#_0Gf})^{-1}}\\\#_1(H\alpha'_f,\beta_{Gx})}
    \end{xy}\,.
  \end{multline}
\end{sidewaysfigure}

\revise{On a 2-cell $\phi$} each of these is defined as
\begin{center}
  \begin{tabular}{l|llll}
    ${*_{-1}}$&$G$&$\alpha$&$A$&$\Gamma$\\\hline
    $H$&$HG\phi$&$H\alpha_\phi$&$\id$&$\id$\\
    $\beta$&$\beta_{G\phi}$&$\id$&$\id$&$\id$\\
    $B$&$\id$&$\id$&$\id$&$\id$\\
    $\Delta$&$\id$&$\id$&$\id$&$\id$\\
  \end{tabular}
\end{center}

\revise{On a 3-cell $\Sigma$} each of these is defined as
\begin{center}
  \begin{tabular}{l|llll}
    ${*_{-1}}$&$G$&$\alpha$&$A$&$\Gamma$\\\hline
    $H$&$HG\Sigma$&$\id$&$\id$&$\id$\\
    $\beta$&$\id$&$\id$&$\id$&$\id$\\
    $B$&$\id$&$\id$&$\id$&$\id$\\
    $\Delta$&$\id$&$\id$&$\id$&$\id$\\
  \end{tabular}
\end{center}

\revise{On a pair of composable 1-cells $f,f'$} each of these is defined as
\begin{center}
  \begin{tabular}{c|cccc}
    ${*_{-1}}$&$G$&$\alpha\from{}G\to{}G'$&$A$&$\Gamma$\\\hline
    $H$&$HG(f'\#_0f)$&$H(\alpha^2_{f',f})$&\revise{$HA_f$?}&$\id$\\
    $\beta\from{}H\to{}H'$&$\beta^2_{Gf',Gf}$&$\id$&$\id$&$\id$\\
    $B$&\revise{$B_{Gf}$?}&$\id$&$\id$&$\id$\\
    $\Delta$&$\id$&$\id$&$\id$&$\id$\\
  \end{tabular}
\end{center}

We also need to define the non-dimension raising horizontal composites:
\begin{center}
  \begin{tabular}{l|llll}
    ${\lhc_{-1}}$&$G$&$\alpha\from{}G\to{}G'$&$A\from\alpha\To\alpha'$&$\Gamma\from{}A\Tto{}A'$\\\hline
    $H$&$H*_{-1}G$&$H*_{-1}\alpha$&$H*_{-1}A$&$H*_{-1}\Gamma$\\
    $\beta\from{}H\to{}H'$&$\beta*_{-1}G$&$(\beta*_{-1}G')*_0(H*_{-1}\alpha)$&$(\beta*_{-1}G')*_0(H*_{-1}A)$&$(\beta*_{-1}G')*_0(H*_{-1}\Gamma)$\\
    $B\from\beta\To\beta'$&$B*_{-1}G$&$(B*_{-1}G')*_0(H*_{-1}\alpha)$&$(B*_{-1}G')*_0(H*_{-1}A)$&$(B*_{-1}G')*_0(H*_{-1}\Gamma)$\\
    $\Delta\from{}B\Tto{}B'$&$\Delta*_{-1}G$&$(\Delta*_{-1}G')*_0(H*_{-1}\alpha)$&$(\Delta*_{-1}G')*_0(H*_{-1}A)$&$(\Delta*_{-1}G')*_0(H*_{-1}\Gamma)$\\
  \end{tabular}
\end{center}

\begin{center}
  \begin{tabular}{l|llll}
    ${\rhc_{-1}}$&$G$&$\alpha\from{}G\to{}G'$&$A\from\alpha\To\alpha'$&$\Gamma\from{}A\Tto{}A'$\\\hline
    $H$&$H*_{-1}G$&$H*_{-1}\alpha$&$H*_{-1}A$&$H*_{-1}\Gamma$\\
    $\beta\from{}H\to{}H'$&$\beta*_{-1}G$&$(H'*_{-1}\alpha)*_0(\beta*_{-1}G)$&$(H'*_{-1}A)*_0(\beta*_{-1}G)$&$(H'*_{-1}\Gamma)*_0(\beta*_{-1}G)$\\
    $B\from\beta\To\beta'$&$B*_{-1}G$&$(H'*_{-1}\alpha)*_0(B*_{-1}G)$&$(H'*_{-1}A)*_0(B*_{-1}G)$&$(H'*_{-1}\Gamma)*_0(B*_{-1}G)$\\
    $\Delta\from{}B\Tto{}B'$&$\Delta*_{-1}G$&$(H'*_{-1}\alpha)*_0(\Delta*_{-1}G)$&$(H'*_{-1}A)*_0(\Delta*_{-1}G)$&$(H'*_{-1}\Gamma)*_0(\Delta*_{-1}G)$
  \end{tabular}
\end{center}

These cells have dimensions:
\begin{center}
  \begin{tabular}{l|llll}
    ${\rhc_{-1}},{\lhc_{-1}}$&$G$&$\alpha\from{}G\to{}G'$&$A\from\alpha\To\alpha'$&$\Gamma\from{}A\Tto{}A'$\\\hline
    $H$&0&1&2&3\\
    $\beta\from{}H\to{}H'$&1&1&2&3\\
    $B\from\beta\To\beta'$&2&2&3&3\\
    $\Delta\from{}B\Tto{}B'$&3&3&3&3
  \end{tabular}
\end{center}

\revise{For reference we spell out the data for
  $\beta\lhc_{-1}\alpha$: are given in \eqref{eq:betalhcalpha01},
  \eqref{eq:betalhcalpha2} and \eqref{eq:betalhcalphacoc}.}
\begin{equation}
  \label{eq:betalhcalpha01}
  \begin{xy}
    \xyboxmatrix{
      HGx\ar[r]^{H\alpha_x}\ar[d]_{HGf}&HG'x\ar[r]^{\beta_{G'x}}\ar[d]|{HG'f}\ar@2[dl]**{}?/-1ex/;?/1ex/^{H\alpha_f}&H'G'x\ar[d]^{H'G'f}\ar@2[dl]**{}?/-1ex/;?/1ex/^{\beta_{G'f}}\\
      HGy\ar[r]_{H\alpha_y}&HG'y\ar[r]_{\beta_{G'y}}&H'G'y }
  \end{xy}\,,
\end{equation}

\begin{sidewaysfigure}
  \begin{center}Details of $(\beta\lhc_{-1}\alpha)_\phi$\end{center}
  \begin{equation}
    \label{eq:betalhcalpha2}
    \begin{xy}
      \xyboxmatrix"A"@+.5cm{
        HGx\ar[r]^{H\alpha_x}\ar[d]|{HGf}="x"\ar@/_3pc/[d]|{HGf'}="y"\ar@2"x";"y"**{}?/-1ex/;?/1ex/^{HG\phi}&H'Gx\ar[r]^{\beta_{G'x}}\ar[d]|{H'Gf}\ar@2[dl]**{}?/-1ex/;?/1ex/^{H(\alpha_f)}&H'G'x\ar[d]^{H'G'f}\ar@2[dl]**{}?/-1ex/;?/1ex/^{\beta_{G'f}}\\
        HGy\ar[r]_{H\alpha_y}&H'Gy\ar[r]_{\beta_{G'y}}&H'G'y}
      \POS;(80,0) \xyboxmatrix"B"@+.5cm{
        HGx\ar[r]^{H\alpha_x}\ar[d]_{HGf'}&H'Gx\ar[r]^{\beta_{G'x}}\ar@/_1.5pc/[d]|{H'Gf'}="x"\ar@/^1.5pc/[d]|{H'Gf}="y"\ar@2"y";"x"**{}?/-1ex/;?/1ex/^{H'G\phi}\ar@2[dl]**{}?/1ex/;?/3ex/_{H(\alpha_{f'})}&H'G'x\ar[d]^{H'G'f}\ar@2[dl]**{}?/-3ex/;?/-1ex/^{\beta_{G'f}}\\
        HGy\ar[r]_{H\alpha_y}&H'Gy\ar[r]_{\beta_{G'y}}&H'G'y}
      \POS;(160,0) \xyboxmatrix"C"@+.5cm{
        HGx\ar[r]^{H\alpha_x}\ar[d]_{HGf'}&H'Gx\ar[r]^{\beta_{G'x}}\ar[d]|{H'Gf'}\ar@2[dl]**{}?/-1ex/;?/1ex/^{H\alpha_{f'}}&H'G'x\ar[d]|{H'G'f'}="x"\ar@/^3pc/[d]|{H'G'f}="y"\ar@2"y";"x"**{}?/-1ex/;?/1ex/^{H'G'\phi}\ar@2[dl]**{}?/-1ex/;?/1ex/^{\beta_{G'f'}}\\
        HGy\ar[r]_{H\alpha_y}&H'Gy\ar[r]_{\beta_{G'y}}&H'G'y}
      \ar@3"A";"B"**{}?/-1ex/;?/1ex/|{}="here"\ar@{.}"here";p+(0,30)*!C\labelbox{(\beta_{G'y}\#_0\underline{H\alpha_\phi})\\\#_1(\beta_{G'f}\#_0H\alpha_x)}
      \ar@3"B";"C"**{}?/-1ex/;?/1ex/|{}="here"\ar@{.}"here";p+(0,30)*!C\labelbox{(\beta_{G'y}\#_0H\alpha_{f'})\\\#_1(\underline{\beta_{G\phi}}\#_0H\alpha_x)}
    \end{xy}
  \end{equation}
\end{sidewaysfigure}

\begin{sidewaysfigure}
  \begin{center}Details of $(\beta\lhc_{-1}\alpha)^2_{f',f}$\end{center}
  \begin{equation}
    \label{eq:betalhcalphacoc}
    \begin{xy}
      \xyboxmatrix"A"{
        HGx\ar[r]^{H\alpha_x}\ar[d]_{HGf}&HG'x\ar[r]^{\beta_{G'x}}\ar[d]|{HG'f}\ar@2[dl]**{}?/-1ex/;?/1ex/_{H\alpha_f}|{}="x"&H'G'x\ar[d]^{H'G'f}\ar@2[dl]**{}?/-1ex/;?/1ex/_{\beta_{G'f}}\\
        HGy\ar[r]^{H\alpha_y}\ar[d]_{HGf'}&HG'y\ar[r]^{\beta_{G'x}}\ar[d]|{HG'f'}\ar@2[dl]**{}?/-1ex/;?/1ex/_{H\alpha_{f'}}&H'G'y\ar[d]^{H'G'f'}\ar@2[dl]**{}?/-1ex/;?/1ex/^{\beta_{G'f'}}|{}="y"\\
        HGz\ar[r]_{H\alpha_z}&HG'z\ar[r]_{\beta_{G'z}}&H'G'z
        \POS\tria "x";"y"
      }
      \POS;(70,0)
      \xyboxmatrix"B"{
        HGx\ar[r]^{H\alpha_x}\ar[d]_{HGf}&HG'x\ar[r]^{\beta_{G'x}}\ar[d]|{HG'f}\ar@2[dl]**{}?/-1ex/;?/1ex/_{H\alpha_f}|{}="x"&H'G'x\ar[d]^{H'G'f}\ar@2[dl]**{}?/-1ex/;?/1ex/_{\beta_{G'f}}\\
        HGy\ar[r]^{H\alpha_y}\ar[d]_{HGf'}&HG'y\ar[r]^{\beta_{G'x}}\ar[d]|{HG'f'}\ar@2[dl]**{}?/-1ex/;?/1ex/_{H\alpha_{f'}}&H'G'y\ar[d]^{H'G'f'}\ar@2[dl]**{}?/-1ex/;?/1ex/^{\beta_{G'f'}}|{}="y"\\
        HGz\ar[r]_{H\alpha_z}&HG'z\ar[r]_{\beta_{G'z}}&H'G'z
        \POS\tria "y";"x"
      }
      \POS+(70,0)
      \xyboxmatrix"C"{
        HGx\ar[r]^{H\alpha_x}\ar[dd]|{HG(f'\#_0f)}&HG'x\ar[r]^{\beta_{G'x}}\ar[dd]|{HG'(f'\#_0f)}\ar@2[ddl]**{}?/-3ex/;?/-1ex/_{H\alpha_{f'\#_0f}}|{}="x"&H'G'x\ar[dd]|{H'G'(f'\#_0f)}\ar@2[ddl]**{}?/-3ex/;?/-1ex/_{\beta_{G'(f'\#_0f)}}\\
        {}&{}&{}\\
        HGz\ar[r]_{H\alpha_z}&HG'z\ar[r]_{\beta_{G'z}}&H'G'z
      }
      \ar@3"A";"B"**{}?/-1ex/;?/1ex/|{}="here"\ar@{.}"here";p+(0,30)*!C\labelbox{(\beta_{G'z}\#_0H\alpha_{f'}\#_0HGf)\\\#_1(\underline{\beta_{G'f'}\ten{}H\alpha_f})\\\#_1(H'G'f'\#_0\beta_{G'f}\#_0H\alpha_x)}
      \ar@3"B";"C"**{}?/-1ex/;?/1ex/|{}="here"\ar@{.}"here";p+(0,30)*!C\labelbox{(\beta_{G'z}\#_0\underline{H\alpha^2_{f',f}})\\\#_1(\underline{\beta^2_{G'f',G'f}}\#_0H\alpha_x)}
    \end{xy}
  \end{equation}
\end{sidewaysfigure}

\revise{We also spell out the data for
  $\beta\rhc_{-1}\alpha$: in dimensions 0 and 1 in
  \eqref{eq:betarhcalpha01}, in dimension 2 in
  \eqref{eq:betarhcalpha2}, and for composable 1-cells in \eqref{eq:betarhcalphacoc}.
}
\begin{equation}
  \label{eq:betarhcalpha01}    
  \begin{xy}
    \xyboxmatrix{
      HGx\ar[r]^{\beta_{Gx}}\ar[d]_{HGf}&H'Gx\ar[r]^{H'\alpha_x}\ar[d]|{H'Gf}\ar@2[dl]**{}?/-1ex/;?/1ex/^{\beta_{Gf}}&H'G'x\ar[d]^{H'G'f}\ar@2[dl]**{}?/-1ex/;?/1ex/^{H'\alpha_f}\\
      HGy\ar[r]_{\beta_{Gy}}&H'Gy\ar[r]_{H'\alpha_y}&H'G'y }\,,
  \end{xy}
\end{equation}
\begin{sidewaysfigure}
  \begin{center}Details of $(\beta\rhc_{-1}\alpha)_{\phi}$\end{center}
  \begin{equation}
    \label{eq:betarhcalpha2}
    \begin{xy}
      \xyboxmatrix"A"@+.5cm{
        HGx\ar[r]^{\beta_{Gx}}\ar[d]|{HGf}="x"\ar@/_3pc/[d]|{HGf'}="y"\ar@2"x";"y"**{}?/-1ex/;?/1ex/^{HG\phi}&H'Gx\ar[r]^{H'\alpha_x}\ar[d]|{H'Gf}\ar@2[dl]**{}?/-1ex/;?/1ex/^{\beta_{Gf}}&H'G'x\ar[d]^{H'G'f}\ar@2[dl]**{}?/-1ex/;?/1ex/^{H'(\alpha_f)}\\
        HGy\ar[r]_{\beta_{Gy}}&H'Gy\ar[r]_{H'\alpha_y}&H'G'y}
      \POS;(80,0) \xyboxmatrix"B"@+.5cm{
        HGx\ar[r]^{\beta_{Gx}}\ar[d]_{HGf'}&H'Gx\ar[r]^{H'\alpha_x}\ar@/_1.5pc/[d]|{H'Gf'}="x"\ar@/^1.5pc/[d]|{H'Gf}="y"\ar@2"y";"x"**{}?/-1ex/;?/1ex/^{H'G\phi}\ar@2[dl]**{}?/1ex/;?/3ex/_{\beta_{Gf'}}&H'G'x\ar[d]^{H'G'f}\ar@2[dl]**{}?/-3ex/;?/-1ex/^{H'(\alpha_f)}\\
        HGy\ar[r]_{\beta_{Gy}}&H'Gy\ar[r]_{H'\alpha_y}&H'G'y}
      \POS;(160,0) \xyboxmatrix"C"@+.5cm{
        HGx\ar[r]^{\beta_{Gx}}\ar[d]_{HGf'}&H'Gx\ar[r]^{H'\alpha_x}\ar[d]|{H'Gf'}\ar@2[dl]**{}?/-1ex/;?/1ex/^{\beta_{Gf'}}&H'G'x\ar[d]|{H'G'f'}="x"\ar@/^3pc/[d]|{H'G'f}="y"\ar@2"y";"x"**{}?/-1ex/;?/1ex/^{H'G'\phi}\ar@2[dl]**{}?/-1ex/;?/1ex/^{H'(\alpha_f)}\\
        HGy\ar[r]_{\beta_{Gy}}&H'Gy\ar[r]_{H'\alpha_y}&H'G'y}
      \ar@3"A";"B"**{}?/-1ex/;?/1ex/|{}="here"\ar@{.}"here";p+(0,30)*!C\labelbox{(H'\alpha_y\#_0\underline{\beta_{G\phi}})\\\#_1(H'\alpha_f\#_0\beta_{Gx})}
      \ar@3"B";"C"**{}?/-1ex/;?/1ex/|{}="here"\ar@{.}"here";p+(0,30)*!C\labelbox{(H'\alpha_y\#_0\beta_{\alpha_{f'}})\\\#_1(\underline{H'\alpha_\phi}\#_0\beta_{Gf})}
    \end{xy}
  \end{equation}
\end{sidewaysfigure}
\begin{sidewaysfigure}
  \begin{center}Details of $(\beta\rhc_{-1}\alpha)^2_{f',f}$\end{center}
  \begin{equation}
    \label{eq:betarhcalphacoc}
    \begin{xy}
      \xyboxmatrix"A"{
        HGx\ar[r]^{\beta_{Gx}}\ar[d]_{HGf}&H'Gx\ar[r]^{H'\alpha_x}\ar[d]|{H'Gf}\ar@2[dl]**{}?/-1ex/;?/1ex/_{\beta_{Gf}}|{}="x"&H'G'x\ar[d]^{H'G'f}\ar@2[dl]**{}?/-1ex/;?/1ex/_{H'_{\alpha_f}}\\
        HGy\ar[r]^{\beta_{Gy}}\ar[d]_{HGf'}&H'Gy\ar[r]^{H'\alpha_y}\ar[d]|{H'Gf'}\ar@2[dl]**{}?/-1ex/;?/1ex/_{\beta_{Gf'}}&H'G'y\ar[d]^{H'G'f'}\ar@2[dl]**{}?/-1ex/;?/1ex/^{H'_{\alpha_{f'}}}|{}="y"\\
        HGz\ar[r]_{\beta_{Gz}}&H'Gz\ar[r]_{H'\alpha_z}&H'G'z
        \POS\tria "x";"y"
      }
      \POS;(70,0)
      \xyboxmatrix"B"{
        HGx\ar[r]^{\beta_{Gx}}\ar[d]_{HGf}&H'Gx\ar[r]^{H'\alpha_x}\ar[d]|{H'Gf}\ar@2[dl]**{}?/-1ex/;?/1ex/_{\beta_{Gf}}|{}="x"&H'G'x\ar[d]^{H'G'f}\ar@2[dl]**{}?/-1ex/;?/1ex/_{H'_{\alpha_f}}\\
        HGy\ar[r]^{\beta_{Gy}}\ar[d]_{HGf'}&H'Gy\ar[r]^{H'\alpha_y}\ar[d]|{H'Gf'}\ar@2[dl]**{}?/-1ex/;?/1ex/_{\beta_{Gf'}}&H'G'y\ar[d]^{H'G'f'}\ar@2[dl]**{}?/-1ex/;?/1ex/^{H'_{\alpha_{f'}}}|{}="y"\\
        HGz\ar[r]_{\beta_{Gz}}&H'Gz\ar[r]_{H'\alpha_z}&H'G'z
        \POS\tria "y";"x"
      }
      \POS+(70,0)
      \xyboxmatrix"C"{
        HGx\ar[r]^{\beta_{Gx}}\ar[dd]|{HG(f'\#_0f)}&H'Gx\ar[r]^{H'\alpha_x}\ar[dd]|{H'G(f'\#_0f)}\ar@2[ddl]**{}?/-3ex/;?/-1ex/_{\beta_{G(f'\#_f)}}|{}="x"&H'G'x\ar[dd]|{H'G'(f'\#_0f)}\ar@2[ddl]**{}?/-3ex/;?/-1ex/_{H'_{\alpha_{f'\#_0f}}}\\
        {}&{}&{}\\
        HGz\ar[r]_{\beta_{Gz}}&H'Gz\ar[r]_{H'\alpha_z}&H'G'z
      }
      \ar@3"A";"B"**{}?/-1ex/;?/1ex/|{}="here"\ar@{.}"here";p+(0,30)*!C\labelbox{(H'\alpha_z\#_0\beta_{\alpha_f}\#_0HGf)\\\#_1(\underline{H'\alpha_{f'}\ten\beta_{\alpha_f}})\\\#_1(H'G'f'\#_0H'\alpha_f\#_0\beta_{Gx})}
      \ar@3"B";"C"**{}?/-1ex/;?/1ex/|{}="here"\ar@{.}"here";p+(0,30)*!C\labelbox{(H'\alpha_z\#_0\underline{\beta^2_{Gf',Gf}})\\\#_1(H'\alpha^2_{f',f}\#_0\beta_{Gx})}
    \end{xy}
  \end{equation}
\end{sidewaysfigure}

\begin{thm}
  \label{thm:hcompthm}
  The compositions defined above result in the following types of cells:
  \begin{center}
    \begin{tabular}{ll|l}
      dimension& type&\\\hline
      0& $\Gray$-functor& $H*_{-1}G$\\
      1& pseudo-transformation& $H*_{-1}\alpha, \beta*_{-1}G$\\
      2& pseudo-modification& $H*_{-1}A, \beta*_{-1}\alpha, B*_{-1}G$\\
      3& perturbation& $H*_{-1}\Gamma, \beta*_{-1}A, B*_{-1}\alpha, \Delta*_{-1}G$
    \end{tabular}
  \end{center}
\end{thm}
\begin{prf}
  In dimension 0 $HG$ is trivially a $\Gray$-functor. Dimension 1 is
  equally trivial. In dimension 2 the only non-trivial case is
  $\beta*_{-1}\alpha$ which is proved in lemma
  \ref{lem:hcompbetaalpha}. In dimension 3 $H*_{-1}\Gamma,
  \Delta*_{-1}G$ are obviously perturbations. That $\beta*_{-1}A,
  B*_{-1}\alpha$ are such is shown in lemmata \ref{lem:hcompbetaA} and
  \ref{lem:hcompBalpha}. \qed
\end{prf}

\begin{lem}
  \label{lem:hcompbetaalpha}
  $\beta*_{-1}\alpha$ as defined above is a pseudo-modification
  \begin{multline}
    \beta*_{-1}\alpha\from\beta\lhc_{-1}\alpha\To\beta\rhc_{-1}\alpha\\=(\beta*_{-1}G')*_0(H*_{-1}\alpha)\To(H'*_{-1}\alpha)*_0(\beta*_{-1}G)\,.
  \end{multline}
\end{lem}
\begin{prf}
  We need to verify that $\beta*_{-1}\alpha$ satisfies the conditions
  \ref{def:psmod}. The conditions of definition
  \eqref{eq:modifex0cellunit} is obvious from
  \eqref{eq:betaalphaf}. 

  Next, we check that our definition satisfies
  \eqref{eq:modifexcoccomp}, that is, we need to verify that
  \eqref{eq:hcompbetaalphacocomp} commutes, which we do in
  \eqref{eq:hcompbetaalphacocompdetail}.
  \begin{sidewaysfigure}
    \begin{equation}
      \label{eq:hcompbetaalphacocomp}
      \begin{xy}
        0;(70,0):(0,-1.5)::
        0
        \xyboxmatrix"A"@ru@+.3cm{
          {}&{}&HGx\ar[r]^{\beta_{Gx}}\ar[d]|{H\alpha_x}\ar[dl]|{HGf}&H'Gx\ar[d]^{H'\alpha_x}\ar@2[dl]**{}?/-1ex/;?/1ex/^{\beta_{\alpha_x}}\\
          {}&HGy\ar[d]|{H\alpha_y}\ar[dl]|{HGf'}&HG'x\ar[r]|{\beta_{G'x}}&H'G'x\ar[dl]|{H'G'f}
          \ar@<+1ex>@2 {[dll];[ll]**{}?(.5)}**{}?(.7)/-1ex/;?(.7)/+1ex/_{(\beta\lhc_{-1}\alpha)_f}\\
          HGz\ar[d]_{H\alpha_z}&HG'y\ar[r]|{\beta_{G'y}}&H'G'y\ar[dl]|{H'G'f'}\ar@<+1ex>@2 {[dll];[ll]**{}?(.5)}**{}?(.7)/-1ex/;?(.7)/+1ex/_{(\beta\lhc_{-1}\alpha)_{f'}}&{}\\
          {}HG'z\ar[r]_{\beta_{G'z}}&H'G'z&{}&{}
        }
        \POS (1,0)
        \xyboxmatrix"B"@ru@+.3cm{
          {}&{}&HGx\ar[r]^{\beta_{Gx}}\ar[dl]|{HGf}&H'Gx\ar[d]^{H'\alpha_x}\ar@<+1ex>@2 {[];[d]**{}?(.5)};[dll]**{}?(.3)/-1ex/;?(.3)/+1ex/_{(\beta\rhc_{-1}\alpha)_f}\\
          {}&HGy\ar[d]|{H\alpha_y}\ar[r]|{\beta_{G'y}}\ar[dl]|{HGf'}&HG'y\ar[d]|{H'\alpha_y}\ar@2[dl]**{}?/-1ex/;?/1ex/^{\beta_{\alpha_y}}&H'G'x\ar[dl]|{H'G'f}\\
          HGz\ar[d]_{H\alpha_z}&HG'y\ar[r]|{\beta_{G'y}}&H'G'y\ar[dl]|{H'G'f'}\ar@<+1ex>@2 {[dll];[ll]**{}?(.5)}**{}?(.7)/-1ex/;?(.7)/+1ex/_{(\beta\lhc_{-1}\alpha)_{f'}}&{}\\
          {}HG'z\ar[r]_{\beta_{G'z}}&H'G'z&{}&{}
        }
        \POS (2,0)
        \xyboxmatrix"C"@ru@+.3cm{
          {}&{}&HGx\ar[r]^{\beta_{Gx}}\ar[dl]|{HGf}&H'Gx\ar[d]^{H'\alpha_x}
          \ar@<+1ex>@2 {[];[d]**{}?(.5)};[dll]**{}?(.3)/-1ex/;?(.3)/+1ex/_{(\beta\rhc_{-1}\alpha)_f}\\
          {}&HGy\ar[r]|{\beta_{Gy}}\ar[dl]|{HGf'}&H'Gy\ar[d]|{H'_{\alpha_y}}\ar@<+1ex>@2 {[];[d]**{}?(.5)};[dll]**{}?(.3)/-1ex/;?(.3)/+1ex/_{(\beta\rhc_{-1}\alpha)_{f'}}&H'G'x\ar[dl]|{H'G'f}\\
          HGz\ar[d]_{H\alpha_z}\ar[r]|{\beta_{Gz}}&H'Gz\ar[d]_{H'\alpha_z}\ar@2[dl]**{}?/-1ex/;?/1ex/^{\beta_{\alpha_z}}&H'G'y\ar[dl]|{H'G'f'}&{}\\
          {}HG'z\ar[r]_{\beta_{G'z}}&H'G'z&{}&{}
        }
        \POS (0,1)
        \xyboxmatrix"D"@ru@+.3cm{
          {}&HGx\ar[r]^{\beta_{Gx}}\ar[d]|{H\alpha_x}\ar[dl]_{HG(f'\#_0f)}&H'Gx\ar[d]^{H'\alpha_x}\ar@2[dl]**{}?/-1ex/;?/1ex/^{\beta_{\alpha_x}}\\
          HGz\ar[d]|{H\alpha_z}&HG'x\ar[r]|{\beta_{G'x}}&H'G'x\ar[dl]^{H'G'(f'\#_0f)}
          \ar@<+1ex>@2 {[dll];[ll]**{}?(.5)}**{}?(.7)/-1ex/;?(.7)/+1ex/|(-1){(\beta\lhc_{-1}\alpha)_{f'\#_0f}}\\
          HG'z\ar[r]|{\beta_{G'z}}&H'G'z&{}  
        }
        \POS (2,1)
        \xyboxmatrix"E"@ru@+.3cm{
          {}&HGx\ar[r]^{\beta_{Gx}}\ar[dl]_{HG(f'\#_0f)}&H'Gx\ar[d]^{H'\alpha_x}
          \ar@2 {[];[d]**{}?(.5)};[dll]**{}?(.4)/-1ex/;?(.4)/1ex/|(-1){(\beta\rhc_{-1}\alpha)_{f'\#_0f}}\\
          HGz\ar[d]_{H\alpha_z}\ar[r]|{\beta_{Gz}}&H'Gz\ar[d]|{H'\alpha_z}
          \ar@2[dl]**{}?/-1ex/;?/1ex/^{\beta_{\alpha_z}}&H'G'x\ar[dl]^{H'G'(f'\#_0f)}\\
          HG'z\ar[r]_{\beta_{G'z}}&H'G'z&{}
        }
        \ar@3"A";"B"**{}?/-1ex/;?/1ex/^{((\beta\lhc_{-1}\alpha)_f\#_0HGf)\\\#_1(H'G'f'\#_0\underline{(\beta*_{-1}\alpha)_f})}
        \ar@3"B";"C"**{}?/-1ex/;?/1ex/^{(\underline{(\beta*_{-1}\alpha)_{f'}}\#_0HGf)\\\#_1(H'G'f'\#_0(\beta\rhc_{-1}\alpha))}
        \ar@3"A";"D"**{}?<;?>**{}?/-1ex/;?/1ex/_{\underline{(\beta\lhc_{-1}\alpha)^2_{f',f}}\\\#_1(H'G'(f'\#_0f)\#_0\beta_{\alpha_x})}
        \ar@3"D";"E"**{}?/-1ex/;?/1ex/_{\underline{(\beta*_{-1}\alpha)_{f'\#_0f}}}
        \ar@3"C";"E"**{}?<;?>**{}?/-1ex/;?/1ex/^{(\beta_{\alpha_z}\#_0HG(f'\#_0f))\\\#_1\underline{(\beta\rhc_{-1}\alpha)^2_{f',f}}}
      \end{xy}
    \end{equation}
  \end{sidewaysfigure}

\begin{sidewaysfigure}
    \begin{equation}
      \label{eq:hcompbetaalphacocompdetail}      
      \def\internalbase{\POS 0;<1cm,.5cm>:<1cm,-.5cm>::0}
      \def\subdiagram#1#2{\POS*!C++\xybox{#2}="#1"}
      \begin{xy}
        0;<3.5cm,0cm>:(0,-1.5)::0
        \subdiagram{A}{
          \internalbase,
          (0,2);(0,3)**\dir{-},(1,1)**\dir{-},
          (0,3);(1,3)**\dir{-},
          (1,1);(1,2)**\dir{-},(2,0)**\dir{-},
          (1,2);(2,2)**\dir{-},
          (1,3);(2,2)**\dir{-},
          (2,0);(3,0)**\dir{-},(2,1)**\dir{-},
          (2,1);(3,1)**\dir{-},
          (2,2);(3,1)**\dir{-},
          (3,0);(3,1)**\dir{-},
          \ar@{}|{\beta_{\alpha_x}}(3,0);(2,1)
          \ar@{}|{(\beta\lhc_{-1}\alpha)_f}(2,1);(1,2)
          \ar@{}|{(\beta\lhc_{-1}\alpha)_{f'}}(1,2);(0,3)
        }
        +(1,0)\subdiagram{B}{
          \internalbase,
          (0,2);(0,3)**\dir{-},(1,1)**\dir{-},
          (0,3);(1,3)**\dir{-},
          (1,1);(1,2)**\dir{-},(2,0)**\dir{-},
          (1,2);(2,2)**\dir{-},
          (1,2);(2,1)**\dir{-},
          (1,3);(2,2)**\dir{-},
          (2,0);(3,0)**\dir{-},(2,1)**\dir{-},         
          (2,2);(3,1)**\dir{-},
          (3,0);(3,1)**\dir{-},
          \ar@{}|{\beta_{G'f\#_0\alpha_x}}(2,1);(3,1)
          \ar@{}|{H\alpha_f}(1,1);(2,1)
          \ar@{}|{(\beta\lhc_{-1}\alpha)_{f'}}(1,2);(0,3)
        }
        +(1,0)\subdiagram{C}{
          \internalbase,
          (0,2);(0,3)**\dir{-},(1,1)**\dir{-},
          (0,3);(1,3)**\dir{-},
          (1,1);(1,2)**\dir{-},(2,0)**\dir{-},
          (1,2);(2,2)**\dir{-},
          (1,3);(2,2)**\dir{-},
          (2,0);(3,0)**\dir{-},
          (2,1);(2,2)**\dir{-},
          (2,1);(3,0)**\dir{-},
          (2,2);(3,1)**\dir{-},
          (3,0);(3,1)**\dir{-},
          \ar@{}|{H'\alpha_f}(2,1);(3,1)
          \ar@{}|{\beta_{\alpha_y\#_0Gf}}(1,1);(2,1)
          \ar@{}|{(\beta\lhc_{-1}\alpha)_{f'}}(1,2);(0,3)
        }
        +(1,0)\subdiagram{D}{
          \internalbase,
          (0,2);(0,3)**\dir{-},(1,1)**\dir{-},
          (0,3);(1,3)**\dir{-},
          (1,1);(1,2)**\dir{-},(2,0)**\dir{-},
          (1,1);(2,1)**\dir{-},
          (1,2);(2,2)**\dir{-},
          (1,3);(2,2)**\dir{-},
          (2,0);(3,0)**\dir{-},
          (2,1);(2,2)**\dir{-},
          (2,2);(3,1)**\dir{-},
          (3,0);(3,1)**\dir{-},
          \ar@{}|{(\beta\rhc_{-1}\alpha)_f}(2,1);(3,0)
          \ar@{}|{\beta_{\alpha_y}}(1,2);(2,1)
          \ar@{}|{(\beta\lhc_{-1}\alpha)_{f'}}(1,2);(0,3)
        }
        +(1,0)\subdiagram{E}{
          \internalbase,
          (0,2);(0,3)**\dir{-},(1,1)**\dir{-},
          (0,3);(1,2)**\dir{-},(1,3)**\dir{-},
          (1,1);(1,2)**\dir{-},(2,0)**\dir{-},(2,1)**\dir{-},
          (1,3);(2,2)**\dir{-},
          (2,0);(3,0)**\dir{-},
          (2,1);(2,2)**\dir{-},
          (2,2);(3,1)**\dir{-},
          (3,0);(3,1)**\dir{-},
          \ar@{}|{(\beta\rhc_{-1}\alpha)_f}(2,1);(3,0)
          \ar@{}|{\beta_{G'f\#_0\alpha_y}}(1,2);(2,2)
          \ar@{}|{H\alpha_{f'}}(0,2);(1,2)
        }
        +(1,0)\subdiagram{F}{
          \internalbase,
          (0,2);(0,3)**\dir{-},(1,1)**\dir{-},
          (0,3);(1,3)**\dir{-},
          (1,1);(2,0)**\dir{-},(2,1)**\dir{-},
          (1,2);(2,1)**\dir{-},(1,3)**\dir{-},
          (1,3);(2,2)**\dir{-},
          (2,0);(3,0)**\dir{-},
          (2,1);(2,2)**\dir{-},
          (2,2);(3,1)**\dir{-},
          (3,0);(3,1)**\dir{-},
          \ar@{}|{(\beta\rhc_{-1}\alpha)_f}(2,1);(3,0)
          \ar@{}|{H'\alpha_{f'}}(1,2);(2,2)
          \ar@{}|{\beta_{\alpha_z\#_0Gf'}}(0,2);(1,2)
        }
        +(1,0)\subdiagram{G}{
          \internalbase,
          (0,2);(0,3)**\dir{-},(1,1)**\dir{-},(1,2)**\dir{-},
          (0,3);(1,3)**\dir{-},
          (1,1);(2,0)**\dir{-},(2,1)**\dir{-},
          (1,2);(1,3)**\dir{-},
          (1,3);(2,2)**\dir{-},
          (2,0);(3,0)**\dir{-},
          (2,1);(2,2)**\dir{-},
          (2,2);(3,1)**\dir{-},
          (3,0);(3,1)**\dir{-},
          \ar@{}|{(\beta\rhc_{-1}\alpha)_f}(2,1);(3,0)
          \ar@{}|{(\beta\rhc_{-1}\alpha)_{f'}}(1,2);(2,1)
          \ar@{}|{\beta_{\alpha_z}}(0,3);(1,2)
        }
        \POS 0+(0,1)
        \subdiagram{A1}{
          \internalbase,
          (0,2);(0,3)**\dir{-},(1,1)**\dir{-},
          (0,3);(1,3)**\dir{-},(1,2)**\dir{-},
          (1,1);(1,2)**\dir{-},(2,0)**\dir{-},
          (1,2);(2,2)**\dir{-},(2,1)**\dir{-},
          (1,3);(2,2)**\dir{-},
          (2,0);(3,0)**\dir{-},(2,1)**\dir{-},
          (2,1);(3,1)**\dir{-},
          (2,2);(3,1)**\dir{-},
          (3,0);(3,1)**\dir{-},
          \ar@{}|{\beta_{\alpha_x}}(2,1);(3,0)
          \ar@{}(1,1);(2,1)|{H\alpha_f}="x"
          \ar@{}|{\beta_{G'f}}(1,2);(3,1)
          \ar@{}|{H\alpha_{f'}}(0,2);(1,2)
          \ar@{}(0,3);(2,2)|{\beta_{G'f'}}="y"
          \tria"y";"x"
        }
        +(1,0)\subdiagram{B1}{
          \internalbase,
          (0,2);(0,3)**\dir{-},(1,1)**\dir{-},
          (0,3);(1,3)**\dir{-},(1,2)**\dir{-},
          (1,1);(1,2)**\dir{-},(2,0)**\dir{-},
          (1,2);(2,2)**\dir{-},(2,1)**\dir{-},
          (1,3);(2,2)**\dir{-},
          (2,0);(3,0)**\dir{-},(2,1)**\dir{-},
          (2,2);(3,1)**\dir{-},
          (3,0);(3,1)**\dir{-},
          \ar@{}|{\beta_{G'f\#_0\alpha_x}}(2,1);(3,1)
          \ar@{}(1,1);(2,1)|{H\alpha_f}="x"
          \ar@{}|{H\alpha_{f'}}(0,2);(1,2)
          \ar@{}(0,3);(2,2)|{\beta_{G'f'}}="y"
          \tria"y";"x"
        }
        +(1,0)\subdiagram{C1}{
          \internalbase,
          (0,2);(0,3)**\dir{-},(1,1)**\dir{-},
          (0,3);(1,3)**\dir{-},(1,2)**\dir{-},
          (1,1);(1,2)**\dir{-},(2,0)**\dir{-},
          (1,2);(2,1)**\dir{-},
          (1,3);(2,2)**\dir{-},
          (2,0);(3,0)**\dir{-},(2,1)**\dir{-},
          (2,2);(3,1)**\dir{-},
          (3,0);(3,1)**\dir{-},
          \ar@{}|{\beta_{G'(f\#_0f)\#_0\alpha_x}}(2,1);(3,1)
          \ar@{}(1,1);(2,1)|{H\alpha_f}
          \ar@{}|{H\alpha_{f'}}(0,2);(1,2)
        }
        +(1,0)\subdiagram{D1}{
          \internalbase,
          (0,2);(0,3)**\dir{-},(1,1)**\dir{-},
          (0,3);(1,3)**\dir{-},(1,2)**\dir{-},
          (1,1);(1,2)**\dir{-},(2,0)**\dir{-},
          (1,3);(2,2)**\dir{-},
          (2,0);(3,0)**\dir{-},
          (2,1);(3,0)**\dir{-},(2,2)**\dir{-},
          (2,2);(3,1)**\dir{-},
          (3,0);(3,1)**\dir{-},
          \ar@{}|{\beta_{G'f'\#_0\alpha_y\#_0Gf}}(1,2);(2,1)
          \ar@{}|{H'\alpha_f}(2,1);(3,1)
          \ar@{}|{H\alpha_{f'}}(0,2);(1,2)
        }
        +(1,0)\subdiagram{E1}{
          \internalbase,
          (0,2);(0,3)**\dir{-},(1,1)**\dir{-},
          (0,3);(1,3)**\dir{-},
          (1,1);(2,0)**\dir{-},
          (1,2);(1,3)**\dir{-},(2,1)**\dir{-},
          (1,3);(2,2)**\dir{-},
          (2,0);(3,0)**\dir{-},
          (2,1);(2,2)**\dir{-},(3,0)**\dir{-},
          (2,2);(3,1)**\dir{-},
          (3,0);(3,1)**\dir{-},
          \ar@{}|{\beta_{\alpha_z\#_0G'(f'\#_0f)}}(0,2);(1,2)
          \ar@{}|{H'\alpha_f}(2,1);(3,1)
          \ar@{}|{H'\alpha_{f'}}(1,2);(2,2)
        }
        +(1,0)\subdiagram{F1}{
          \internalbase,
          (0,2);(0,3)**\dir{-},(1,1)**\dir{-},
          (0,3);(1,3)**\dir{-},
          (1,1);(2,0)**\dir{-},(2,1)**\dir{-},
          (1,2);(1,3)**\dir{-},(2,1)**\dir{-},
          (1,3);(2,2)**\dir{-},
          (2,0);(3,0)**\dir{-},
          (2,1);(2,2)**\dir{-},(3,0)**\dir{-},
          (2,2);(3,1)**\dir{-},
          (3,0);(3,1)**\dir{-},
          \ar@{}|{\beta_{\alpha_z\#_0G'f'}}(0,2);(1,2)
          \ar@{}(1,1);(3,0)|{\beta_{Gf}}="x"
          \ar@{}|{H'\alpha_f}(2,1);(3,1)
          \ar@{}(1,2);(2,2)|{H'\alpha_{f'}}="y"
          \tria"y";"x"
        }
        +(1,0)\subdiagram{G1}{
          \internalbase,
          (0,2);(0,3)**\dir{-},(1,1)**\dir{-},(1,2)**\dir{-},
          (0,3);(1,3)**\dir{-},
          (1,1);(2,0)**\dir{-},(2,1)**\dir{-},
          (1,2);(1,3)**\dir{-},(2,1)**\dir{-},
          (1,3);(2,2)**\dir{-},
          (2,0);(3,0)**\dir{-},
          (2,1);(2,2)**\dir{-},(3,0)**\dir{-},
          (2,2);(3,1)**\dir{-},
          (3,0);(3,1)**\dir{-},
          \ar@{}|{\beta_{Gf'}}(0,2);(2,1)
          \ar@{}|{\beta_{\alpha_z}}(0,3);(1,2)
          \ar@{}(1,1);(3,0)|{\beta_{Gf}}="x"
          \ar@{}|{H'\alpha_f}(2,1);(3,1)
          \ar@{}(1,2);(2,2)|{H'\alpha_{f'}}="y"
          \tria"y";"x"
        }
        \POS 0+(0,2)
        \subdiagram{A2}{
          \internalbase,
          (0,2);(0,3)**\dir{-},(1,1)**\dir{-},
          (0,3);(1,3)**\dir{-},
          (1,1);(2,0)**\dir{-},
          (1,3);(2,2)**\dir{-},
          (2,0);(3,0)**\dir{-},(2,1)**\dir{-},
          (2,1);(3,1)**\dir{-},
          (2,2);(3,1)**\dir{-},
          (3,0);(3,1)**\dir{-},
          \ar@{}(2,1);(3,0)|{\beta_{\alpha_x}}
          \ar@{}(0,3);(2,1)|{(\beta\lhc_{-1}\alpha)_{f'\#_0f}}
        }
        +(1,0)\subdiagram{B2}{
        }
        +(1,0)\subdiagram{C2}{
          \internalbase,
          (0,2);(0,3)**\dir{-},(1,1)**\dir{-},
          (0,3);(1,3)**\dir{-},(2,1)**\dir{-},
          (1,1);(2,0)**\dir{-},
          (1,3);(2,2)**\dir{-},
          (2,0);(3,0)**\dir{-},(2,1)**\dir{-},
          (2,2);(3,1)**\dir{-},
          (3,0);(3,1)**\dir{-},
          \ar@{}(0,2);(2,1)|{H\alpha_{f'\#_0f}}
          \ar@{}(2,1);(3,1)|{\beta_{G'(f'\#_0f)\#_0\alpha_x}}
        }
        +(1,0)\subdiagram{D2}{
        }
        +(1,0)\subdiagram{E2}{
          \internalbase,
          (0,2);(0,3)**\dir{-},(1,1)**\dir{-},
          (0,3);(1,3)**\dir{-},
          (1,1);(2,0)**\dir{-},
          (1,2);(1,3)**\dir{-},(3,0)**\dir{-},
          (1,3);(2,2)**\dir{-},
          (2,0);(3,0)**\dir{-},
          (2,2);(3,1)**\dir{-},
          (3,0);(3,1)**\dir{-},
          \ar@{}(0,2);(1,2)|{\beta_{\alpha_z\#_0G(f'\#_0f)}}
          \ar@{}(1,2);(3,1)|{H'\alpha_{f'\#_0f}}
        }
        +(1,0)\subdiagram{F2}{
        }
        +(1,0)\subdiagram{G2}{
          \internalbase,
          (0,2);(0,3)**\dir{-},(1,1)**\dir{-},(1,2)**\dir{-},
          (0,3);(1,3)**\dir{-},
          (1,1);(2,0)**\dir{-},
          (1,2);(1,3)**\dir{-},
          (1,3);(2,2)**\dir{-},
          (2,0);(3,0)**\dir{-},
          (2,2);(3,1)**\dir{-},
          (3,0);(3,1)**\dir{-},
          \ar@{}(1,2);(3,0)|{(\beta\rhc_{-1}\alpha)_{f'\#_0f}}
          \ar@{}(0,3);(1,2)|{\beta_{\alpha_z}}
        }
        \ttar"A";"B"
        \ttar"A";"A1"
        \ttar"B";"C"\ttar"B";"B1"
        \ttar"C";"D"\ttar"C";"D1"\POS{"B";"C"**{}?};"C1"**{}?*\txt{\eqref{eq:pstransf12whiskright}}
        \ttar"D";"E"\POS"D";"D1"**{}?*\txt{\eqref{eq:pstransf2cocy}}
        \ttar"E";"F"\POS{"E";"F"**{}?};"E1"**{}?*\txt{\eqref{eq:pstransf12whiskleft}}
        \ttar"F";"G"\ttar"F";"F1"
        \ttar"G";"G1"
        \ttar"A1";"B1"\ttar"A1";"A2"
        \ttar"B1";"C1"
        \ttar"C1";"D1"\ttar"C1";"C2"\ar@3@/_7pc/"C1";"E1"\ar@{}@/_7pc/"C1";"E1"|{}="X"\POS"D1";"X"**{}?<>(.5)*\txt{\eqref{eq:pstransf2cellcomp}}\POS"D1"
        \ttar"D1";"E1"\ttar"D1";"E"
        \ttar"E1";"F1"\ttar"E1";"E2"
        \ttar"F1";"G1"
        \ttar"G1";"G2"
        \ttar"A2";"C2"
        \ttar"C2";"E2"\POS{"C2";"E2"?};"X"**{}?*\txt{\eqref{eq:pstransf3cell}}
        \ttar"E2";"G2"
      \end{xy}
    \end{equation}
    \begin{center}
      Verification of \eqref{eq:hcompbetaalphacocomp}. Unlabeled
      subdiagrams commute by naturality. All 3-cells shown are
      canonical and invertible.
    \end{center}
  \end{sidewaysfigure}

  Finally, we check that \eqref{eq:modifex2cell} holds for
  $\beta*_{-1}\alpha$, i.e.\ that \eqref{eq:hcompbetaalpha2cellcomp}
  holds, this is carried out in \eqref{eq:hcompbetaalpha2cellcompdetail}.
  \begin{sidewaysfigure}
    \begin{equation}
      \label{eq:hcompbetaalpha2cellcomp}
      \def\subdiagram#1#2{\xyboxmatrix"#1"@ru@+.3cm{#2}\POS}
      \begin{xy}
        0;(90,0):p
        \subdiagram{A}{
          {}&HGx\ar[r]^{\beta_{Gx}}\ar[d]|{H\alpha_x}\ar[dl]|{HGf}="x"\ar@/_3pc/[dl]|{HGf'}="y"\tarb{HG\phi}"x";"y"&H'Gx\ar[d]^{H'\alpha_x}\ar@2[dl]**{}?/-1ex/;?/1ex/^{\beta_{\alpha_x}}\\
          HGy\ar[d]_{H\alpha_y}&HG'x\ar[r]|{\beta_{G'x}}&H'G'x\ar[dl]|{H'G'f}\ar@2 {[dll];[ll]**{}?}**{}?(.6)/-1ex/;?(.6)/+1ex/_(1){(\beta\lhc_{-1}\alpha)_f}\\
          HG'\ar[r]_{\beta_{G'y}}&H'G'y&{}
        }+(1,0)
        \subdiagram{B}{
          {}&HGx\ar[r]^{\beta_{Gx}}\ar[dl]|{HGf}="x"\ar@/_3pc/[dl]|{HGf'}="y"\tarb{HG\phi}"x";"y"\POS="z"&H'Gx\ar[d]^{H'\alpha_x}\ar@2 {[];[d]**{}?};[dll]**{}?(.3)/-1ex/;?(.3)/1ex/^(1){(\beta\rhc_{-1}\alpha)_f}\\
          HGy\ar[d]_{H\alpha_y}\ar[r]|{\beta_{Gy}}&H'Gy\ar@2[dl]**{}?/-1ex/;?/+1ex/_{\beta_{\alpha_y}}="w"\ar[d]|{H'\alpha_x}&H'G'x\ar[dl]|{H'G'f}\\
          HG'\ar[r]_{\beta_{G'y}}&H'G'y&{}
          \tria"w";"z"
        }+(1,0)
        \subdiagram{C}{
          {}&HGx\ar[r]^{\beta_{Gx}}\ar[dl]|{HGf}="x"\ar@/_3pc/[dl]|{HGf'}="y"\tarb{HG\phi}"x";"y"\POS="z"&H'Gx\ar[d]^{H'\alpha_x}\ar@2 {[];[d]**{}?};[dll]**{}?(.3)/-1ex/;?(.3)/1ex/^(1){(\beta\rhc_{-1}\alpha)_f}\\
          HGy\ar[d]_{H\alpha_y}\ar[r]|{\beta_{Gy}}&H'Gy\ar@2[dl]**{}?/-1ex/;?/+1ex/_{\beta_{\alpha_y}}="w"\ar[d]|{H'\alpha_x}&H'G'x\ar[dl]|{H'G'f}\\
          HG'\ar[r]_{\beta_{G'y}}&H'G'y&{}
          \tria"z";"w"
        }
        \POS (0,-1)
        \subdiagram{A1}{
          {}&HGx\ar[r]^{\beta_{Gx}}\ar[d]|{H\alpha_x}\ar[dl]|{HGf'}&H'Gx\ar[d]^{H'\alpha_x}\ar@2[dl]**{}?/-1ex/;?/1ex/^{\beta_{\alpha_x}}="z"\\
          HGy\ar[d]_{H\alpha_y}&HG'x\ar[r]|{\beta_{G'x}}&H'G'x\ar[dl]|{H'G'f'}="x"\ar@/^3pc/[dl]|{H'G'f}="y"\tarb{H'G'\phi}"y";"x"\POS="w"\ar@2 {[dll];[ll]**{}?}**{}?(.6)/-1ex/;?(.6)/+1ex/_(1){(\beta\lhc_{-1}\alpha)_{f'}}\\
          HG'\ar[r]_{\beta_{G'y}}&H'G'y&{}
          \tria"z";"w"
        }+(1,0)
        \subdiagram{B1}{
          {}&HGx\ar[r]^{\beta_{Gx}}\ar[d]|{H\alpha_x}\ar[dl]|{HGf'}&H'Gx\ar[d]^{H'\alpha_x}\ar@2[dl]**{}?/-1ex/;?/1ex/^{\beta_{\alpha_x}}="z"\\
          HGy\ar[d]_{H\alpha_y}&HG'x\ar[r]|{\beta_{G'x}}&H'G'x\ar[dl]|{H'G'f'}="x"\ar@/^3pc/[dl]|{H'G'f}="y"\tarb{H'G'\phi}"y";"x"\POS="w"\ar@2 {[dll];[ll]**{}?}**{}?(.6)/-1ex/;?(.6)/+1ex/_(1){(\beta\lhc_{-1}\alpha)_{f'}}\\
          HG'\ar[r]_{\beta_{G'y}}&H'G'y&{}
          \tria"w";"z"
        }+(1,0)
        \subdiagram{C1}{
          {}&HGx\ar[r]^{\beta_{Gx}}\ar[dl]|{HGf'}&H'Gx\ar[d]^{H'\alpha_x}\ar@2 {[];[d]**{}?};[dll]**{}?(.3)/-1ex/;?(.3)/1ex/^(1){(\beta\rhc_{-1}\alpha)_{f'}}\\
          HGy\ar[d]_{H\alpha_y}\ar[r]|{\beta_{Gy}}&H'Gy\ar@2[dl]**{}?/-1ex/;?/+1ex/_{\beta_{\alpha_y}}\ar[d]|{H'\alpha_x}&H'G'x\ar[dl]|{H'G'f'}="x"\ar@/^3pc/[dl]|{H'G'f}="y"\tarb{H'G'\phi}"y";"x"\\
          HG'\ar[r]_{\beta_{G'y}}&H'G'y&{}
        }
        \ttarb{\underline{(\beta\lhc_{-1}\alpha)_\phi}\\\#_1(H'G'f\#_0\beta_{\alpha_x})}"A";"A1"
        \ttara{(\beta_{G'y}\#_0H(\alpha_y\#_0G\phi))\\\#_1\underline{(\beta*_{-1}\alpha)_f}}"A";"B"
        \ttarb{(\beta\lhc_{-1}\alpha)_{f'}\\\#_1\underline{(H'G'\phi\ten\beta_{\alpha_x})}}"A1";"B1"
        \ttara{\underline{(\beta_{\alpha_y}\ten{}HG\phi)^{-1}}\\\#_1(\beta\rhc_{-1}\alpha)_f}"B";"C"
        \ttarb{(H'(G'\phi\#_0\alpha_x)\#_0\beta_{Gx})\\\#_1\underline{(\beta*_{-1}\alpha)_{f'}}}"B1";"C1"
        \ttara{(\beta_{\alpha_y}\#_0HGf')\\\#_1\underline{(\beta\rhc_{-1}\alpha)_f}}"C";"C1"
      \end{xy}
    \end{equation}
  \end{sidewaysfigure}
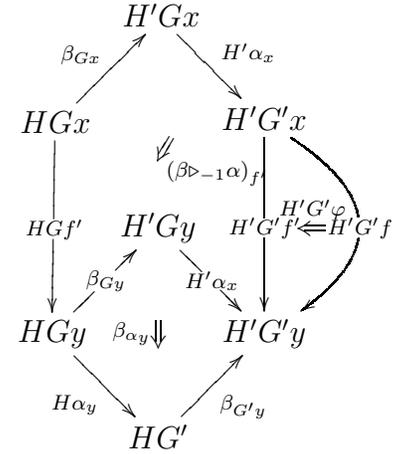
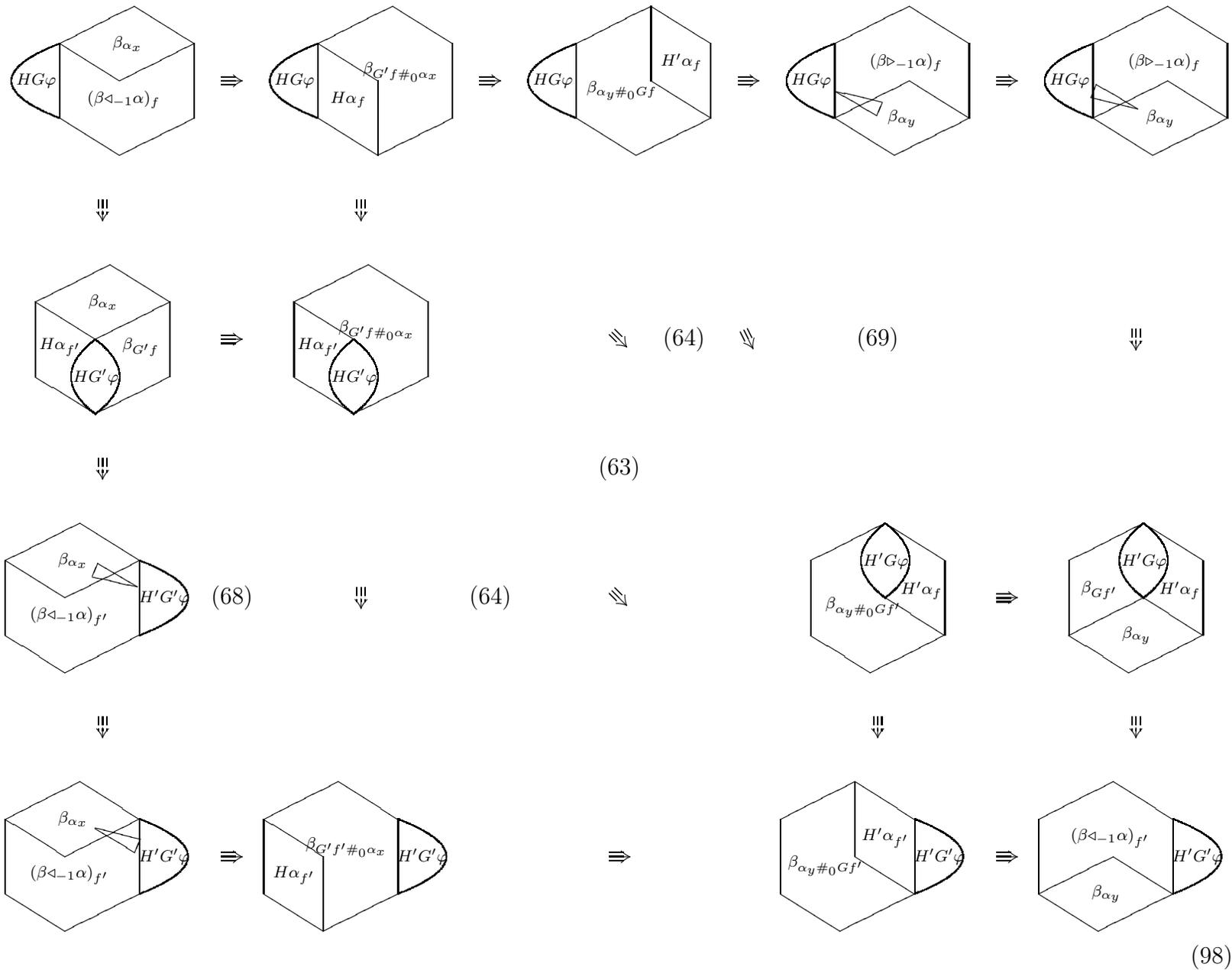
  \begin{sidewaysfigure}
    \begin{equation}
      \label{eq:hcompbetaalpha2cellcompdetail}
      \def\internalbase{\POS 0;<.13cm,0cm>:0,}
      \def\subdiagram#1#2{\POS*!C++\xybox{#2}="#1",}
      \tbdef{(10,5)}{(8,-5)}{(0,-10)}
      \begin{xy}
        0;<4.5cm,0cm>:(0,-1)::0
        \subdiagram{A0}{
          \internalbase
          \ar@{-}{\t(0,0,0)*{}="x"};{\t(1,0,0)*{}="y"}
          \ar@{-}"y";{\t(1,1,0)}*{}="w"
          \ar@{-}"w";{\t(0,1,0)}*{}="z"
          \ar@{-}"x";"z"
          \ar@{}{\t(0,0,1)}*{}="x'";{\t(1,0,1)}*{}="y'"
          \ar@{}"y'";{\t(1,1,1)}*{}="w'"
          \ar@{-}"w'";{\t(0,1,1)}*{}="z'"
          \ar@{-}"x'";"z'"
          \ar@{-}"x";"x'"_{}="X"
          \ar@{-}@/_2pc/"x";"x'"^{}="Y"
          \ar@{}"y";"y'"
          \ar@{}"z";"z'"
          \ar@{-}"w";"w'"
          \ar@{}"X";"Y"|{HG\phi}
          \ar@{}"y";"z"|{\beta_{\alpha_x}}
          \ar@{}"z";{"x'";"w'"**{}?}|{(\beta\lhc_{-1}\alpha)_f}
        }
        +(1,0)        
        \subdiagram{B0}{
          \internalbase
          \ar@{-}{\t(0,0,0)*{}="x"};{\t(1,0,0)*{}="y"}
          \ar@{-}"y";{\t(1,1,0)}*{}="w"
          \ar@{}"w";{\t(0,1,0)}*{}="z"
          \ar@{-}"x";"z"
          \ar@{}{\t(0,0,1)}*{}="x'";{\t(1,0,1)}*{}="y'"
          \ar@{}"y'";{\t(1,1,1)}*{}="w'"
          \ar@{-}"w'";{\t(0,1,1)}*{}="z'"
          \ar@{-}"x'";"z'"
          \ar@{-}"x";"x'"_{}="X"
          \ar@{-}@/_2pc/"x";"x'"^{}="Y"
          \ar@{}"y";"y'"
          \ar@{-}"z";"z'"
          \ar@{-}"w";"w'"
          \ar@{}"X";"Y"|{HG\phi}
          \ar@{}"z";"x'"|{H\alpha_f}
          \ar@{}"z";{"y";"w'"**{}?}|{\beta_{G'f\#_0\alpha_x}}
        }
        +(1,0)        
        \subdiagram{C0}{
          \internalbase
          \ar@{-}{\t(0,0,0)*{}="x"};{\t(1,0,0)*{}="y"}
          \ar@{-}"y";{\t(1,1,0)}*{}="w"
          \ar@{}"w";{\t(0,1,0)}*{}="z"
          \ar@{}"x";"z"
          \ar@{}{\t(0,0,1)}*{}="x'";{\t(1,0,1)}*{}="y'"
          \ar@{-}"y'";{\t(1,1,1)}*{}="w'"
          \ar@{-}"w'";{\t(0,1,1)}*{}="z'"
          \ar@{-}"x'";"z'"
          \ar@{-}"x";"x'"_{}="X"
          \ar@{-}@/_2pc/"x";"x'"^{}="Y"
          \ar@{-}"y";"y'"
          \ar@{}"z";"z'"
          \ar@{-}"w";"w'"
          \ar@{}"X";"Y"|{HG\phi}
          \ar@{}"w";"y'"|{H'\alpha_f}
          \ar@{}"z";{"x";"z'"**{}?}|{\beta_{\alpha_y\#_0Gf}}}
        +(1,0)        
        \subdiagram{D0}{
          \internalbase
          \ar@{-}{\t(0,0,0)*{}="x"};{\t(1,0,0)*{}="y"}
          \ar@{-}"y";{\t(1,1,0)}*{}="w"
          \ar@{}"w";{\t(0,1,0)}*{}="z"
          \ar@{}"x";"z"
          \ar@{-}{\t(0,0,1)}*{}="x'";{\t(1,0,1)}*{}="y'"
          \ar@{-}"y'";{\t(1,1,1)}*{}="w'"
          \ar@{-}"w'";{\t(0,1,1)}*{}="z'"
          \ar@{-}"x'";"z'"
          \ar@{-}"x";"x'"_{}="X"
          \ar@{-}@/_2pc/"x";"x'"^{}="Y"
          \ar@{}"y";"y'"
          \ar@{}"z";"z'"
          \ar@{-}"w";"w'"
          \ar@{}"X";"Y"|{HG\phi}="Z"
          \ar@{}"x'";"w'"|{\beta_{\alpha_y}}="W"
          \ar@{}"y";{"x'";"w'"**{}?}|{(\beta\rhc_{-1}\alpha)_f}
          \tria"W";"Z"
        }      
        +(1,0)
        \subdiagram{E0}{
          \internalbase
          \ar@{-}{\t(0,0,0)*{}="x"};{\t(1,0,0)*{}="y"}
          \ar@{-}"y";{\t(1,1,0)}*{}="w"
          \ar@{}"w";{\t(0,1,0)}*{}="z"
          \ar@{}"x";"z"
          \ar@{-}{\t(0,0,1)}*{}="x'";{\t(1,0,1)}*{}="y'"
          \ar@{-}"y'";{\t(1,1,1)}*{}="w'"
          \ar@{-}"w'";{\t(0,1,1)}*{}="z'"
          \ar@{-}"x'";"z'"
          \ar@{-}"x";"x'"_{}="X"
          \ar@{-}@/_2pc/"x";"x'"^{}="Y"
          \ar@{}"y";"y'"
          \ar@{}"z";"z'"
          \ar@{-}"w";"w'"
          \ar@{}"X";"Y"|{HG\phi}="Z"
          \ar@{}"x'";"w'"|{\beta_{\alpha_y}}="W"
          \ar@{}"y";{"x'";"w'"**{}?}|{(\beta\rhc_{-1}\alpha)_f}
          \tria"Z";"W"
        }
        ,(0,1)
        \subdiagram{A1}{
          \internalbase
          \ar@{-}{\t(0,0,0)*{}="x"};{\t(1,0,0)*{}="y"}
          \ar@{-}"y";{\t(1,1,0)}*{}="w"
          \ar@{-}"w";{\t(0,1,0)}*{}="z"
          \ar@{-}"x";"z"
          \ar@{}{\t(0,0,1)}*{}="x'";{\t(1,0,1)}*{}="y'"
          \ar@{}"y'";{\t(1,1,1)}*{}="w'"
          \ar@{-}"w'";{\t(0,1,1)}*{}="z'"
          \ar@{-}"x'";"z'"
          \ar@{-}"x";"x'"
          \ar@{}"y";"y'"
          \ar@{-}@/^1pc/"z";"z'"_{}="X"
          \ar@{-}@/_1pc/"z";"z'"^{}="Y"
          \ar@{-}"w";"w'"
          \ar@{}"X";"Y"|{HG'\phi}
          \ar@{}"x";"z'"|(.4){H\alpha_{f'}}
          \ar@{}"y";"z"|{\beta_{\alpha_x}}
          \ar@{}"w";"z'"|(.4){\beta_{G'f}}
        }
        +(1,0)        
        \subdiagram{B1}{
          \internalbase
          \ar@{-}{\t(0,0,0)*{}="x"};{\t(1,0,0)*{}="y"}
          \ar@{-}"y";{\t(1,1,0)}*{}="w"
          \ar@{}"w";{\t(0,1,0)}*{}="z"
          \ar@{-}"x";"z"
          \ar@{}{\t(0,0,1)}*{}="x'";{\t(1,0,1)}*{}="y'"
          \ar@{}"y'";{\t(1,1,1)}*{}="w'"
          \ar@{-}"w'";{\t(0,1,1)}*{}="z'"
          \ar@{-}"x'";"z'"
          \ar@{-}"x";"x'"
          \ar@{}"y";"y'"
          \ar@{-}@/^1pc/"z";"z'"_{}="X"
          \ar@{-}@/_1pc/"z";"z'"^{}="Y"
          \ar@{-}"w";"w'"
          \ar@{}"X";"Y"|{HG'\phi}
          \ar@{}"x";"z'"|(.4){H\alpha_{f'}}
          \ar@{}"z";{"y";"w'"**{}?}|{\beta_{G'f\#_0\alpha_x}}
        }
        +(1,0)        
        +(1,0)        
        +(1,0)
        ,(0,2)
        \subdiagram{A2}{
          \internalbase
          \ar@{-}{\t(0,0,0)*{}="x"};{\t(1,0,0)*{}="y"}
          \ar@{-}"y";{\t(1,1,0)}*{}="w"
          \ar@{-}"w";{\t(0,1,0)}*{}="z"
          \ar@{-}"x";"z"
          \ar@{}{\t(0,0,1)}*{}="x'";{\t(1,0,1)}*{}="y'"
          \ar@{}"y'";{\t(1,1,1)}*{}="w'"
          \ar@{-}"w'";{\t(0,1,1)}*{}="z'"
          \ar@{-}"x'";"z'"
          \ar@{-}"x";"x'"
          \ar@{}"y";"y'"
          \ar@{}"z";"z'"
          \ar@{-}"w";"w'"_{}="X"
          \ar@{-}@/^2pc/"w";"w'"^{}="Y"
          \ar@{}"X";"Y"|{H'G'\phi}="W"
          \ar@{}"y";"z"|{\beta_{\alpha_x}}="Z"
          \ar@{}"z";{"x'";"w'"**{}?}|{(\beta\lhc_{-1}\alpha)_{f'}}
          \tria"Z";"W"
        }
        +(1,0)        
        +(1,0)        
        +(1,0)
        \subdiagram{D2}{
          \internalbase
          \ar@{-}{\t(0,0,0)*{}="x"};{\t(1,0,0)*{}="y"}
          \ar@{-}"y";{\t(1,1,0)}*{}="w"
          \ar@{}"w";{\t(0,1,0)}*{}="z"
          \ar@{}"x";"z"
          \ar@{}{\t(0,0,1)}*{}="x'";{\t(1,0,1)}*{}="y'"
          \ar@{-}"y'";{\t(1,1,1)}*{}="w'"
          \ar@{-}"w'";{\t(0,1,1)}*{}="z'"
          \ar@{-}"x'";"z'"
          \ar@{-}"x";"x'"
          \ar@{}"y";"y'"
          \ar@{-}@/^1pc/"y";"y'"_{}="X"
          \ar@{-}@/_1pc/"y";"y'"^{}="Y"
          \ar@{-}"w";"w'"
          \ar@{}"X";"Y"|{H'G\phi}
          \ar@{}"y";"w'"|(.6){H'\alpha_{f}}
          \ar@{}"y'";{"x";"z'"**{}?}|{\beta_{\alpha_y\#_0Gf'}}
        }
        +(1,0)
        \subdiagram{E2}{
          \internalbase
          \ar@{-}{\t(0,0,0)*{}="x"};{\t(1,0,0)*{}="y"}
          \ar@{-}"y";{\t(1,1,0)}*{}="w"
          \ar@{}"w";{\t(0,1,0)}*{}="z"
          \ar@{}"x";"z"
          \ar@{-}{\t(0,0,1)}*{}="x'";{\t(1,0,1)}*{}="y'"
          \ar@{-}"y'";{\t(1,1,1)}*{}="w'"
          \ar@{-}"w'";{\t(0,1,1)}*{}="z'"
          \ar@{-}"x'";"z'"
          \ar@{-}"x";"x'"
          \ar@{}"y";"y'"
          \ar@{-}@/^1pc/"y";"y'"_{}="X"
          \ar@{-}@/_1pc/"y";"y'"^{}="Y"
          \ar@{-}"w";"w'"
          \ar@{}"X";"Y"|{H'G\phi}
          \ar@{}"y";"w'"|(.6){H'\alpha_{f}}
          \ar@{}"x'";"y"|(.4){\beta_{Gf'}}
          \ar@{}"y'";"z'"|{\beta_{\alpha_y}}
        }
        ,(0,3)
        \subdiagram{A3}{
          \internalbase
          \ar@{-}{\t(0,0,0)*{}="x"};{\t(1,0,0)*{}="y"}
          \ar@{-}"y";{\t(1,1,0)}*{}="w"
          \ar@{-}"w";{\t(0,1,0)}*{}="z"
          \ar@{-}"x";"z"
          \ar@{}{\t(0,0,1)}*{}="x'";{\t(1,0,1)}*{}="y'"
          \ar@{}"y'";{\t(1,1,1)}*{}="w'"
          \ar@{-}"w'";{\t(0,1,1)}*{}="z'"
          \ar@{-}"x'";"z'"
          \ar@{-}"x";"x'"
          \ar@{}"y";"y'"
          \ar@{}"z";"z'"
          \ar@{-}"w";"w'"_{}="X"
          \ar@{-}@/^2pc/"w";"w'"^{}="Y"
          \ar@{}"X";"Y"|{H'G'\phi}="W"
          \ar@{}"y";"z"|{\beta_{\alpha_x}}="Z"
          \ar@{}"z";{"x'";"w'"**{}?}|{(\beta\lhc_{-1}\alpha)_{f'}}
          \tria"W";"Z"
        }
        +(1,0)        
        \subdiagram{B3}{
          \internalbase
          \ar@{-}{\t(0,0,0)*{}="x"};{\t(1,0,0)*{}="y"}
          \ar@{-}"y";{\t(1,1,0)}*{}="w"
          \ar@{}"w";{\t(0,1,0)}*{}="z"
          \ar@{-}"x";"z"
          \ar@{}{\t(0,0,1)}*{}="x'";{\t(1,0,1)}*{}="y'"
          \ar@{}"y'";{\t(1,1,1)}*{}="w'"
          \ar@{-}"w'";{\t(0,1,1)}*{}="z'"
          \ar@{-}"x'";"z'"
          \ar@{-}"x";"x'"
          \ar@{}"y";"y'"
          \ar@{-}"z";"z'"
          \ar@{-}"w";"w'"_{}="X"
          \ar@{-}@/^2pc/"w";"w'"^{}="Y"
          \ar@{}"X";"Y"|{H'G'\phi}
          \ar@{}"z";{"y";"w'"**{}?}|{\beta_{G'f'\#_0\alpha_x}}
          \ar@{}"x'";"z"|{H\alpha_{f'}}
        }
        +(1,0)        
        +(1,0)        
        \subdiagram{D3}{
          \internalbase
          \ar@{-}{\t(0,0,0)*{}="x"};{\t(1,0,0)*{}="y"}
          \ar@{-}"y";{\t(1,1,0)}*{}="w"
          \ar@{}"w";{\t(0,1,0)}*{}="z"
          \ar@{}"x";"z"
          \ar@{}{\t(0,0,1)}*{}="x'";{\t(1,0,1)}*{}="y'"
          \ar@{-}"y'";{\t(1,1,1)}*{}="w'"
          \ar@{-}"w'";{\t(0,1,1)}*{}="z'"
          \ar@{-}"x'";"z'"
          \ar@{-}"x";"x'"
          \ar@{-}"y";"y'"
          \ar@{}"z";"z'"
          \ar@{-}"w";"w'"_{}="X"
          \ar@{-}@/^2pc/"w";"w'"^{}="Y"
          \ar@{}"X";"Y"|{H'G'\phi}
          \ar@{}"w";"y'"|{H'\alpha_{f'}}
          \ar@{}"z";{"x";"z'"**{}?}|{\beta_{\alpha_y\#_0Gf'}}
        }
        +(1,0)
        \subdiagram{E3}{
          \internalbase
          \ar@{-}{\t(0,0,0)*{}="x"};{\t(1,0,0)*{}="y"}
          \ar@{-}"y";{\t(1,1,0)}*{}="w"
          \ar@{}"w";{\t(0,1,0)}*{}="z"
          \ar@{}"x";"z"
          \ar@{-}{\t(0,0,1)}*{}="x'";{\t(1,0,1)}*{}="y'"
          \ar@{-}"y'";{\t(1,1,1)}*{}="w'"
          \ar@{-}"w'";{\t(0,1,1)}*{}="z'"
          \ar@{-}"x'";"z'"
          \ar@{-}"x";"x'"
          \ar@{}"y";"y'"
          \ar@{}"z";"z'"
          \ar@{-}"w";"w'"_{}="X"
          \ar@{-}@/^2pc/"w";"w'"^{}="Y"
          \ar@{}"X";"Y"|{H'G'\phi}
          \ar@{}"x'";"w'"|{\beta_{\alpha_y}}
          \ar@{}"y";{"x'";"w'"**{}?}|{(\beta\lhc_{-1}\alpha)_{f'}}
        }
        \ttar"A0";"B0"
        \ttar"A1";"A2"
        \ttar"A1";"B1"
        \ttar"A2";"A3"
        \ttar"A3";"B3"
        \ttar"B0";"B1"
        \ttar"B0";"C0"
        \ttar"B0";"D2"
        \ttar"B1";"B3"
        \ttar"B1";"D3"        
        \ttar"B3";"D3"
        \ttar"C0";"D0"
        \ttar"C0";"D2"
        \ttar"D0";"E0"
        \ttar"D2";"D3"
        \ttar"D2";"E2"
        \ttar"D3";"E3"
        \ttar"E0";"E2"
        \ttar"E2";"E3"
        \ttar"A0";"A1"
        \POS"A1";"B3"**{}?*\txt{\eqref{eq:pstransf12whiskleft}}
        \POS{"B0";"C0"**{}?};"D2"**{}?*\txt{\eqref{eq:pstransf2cellcomp}}
        \POS"C0";"E2"**{}?*\txt{\eqref{eq:pstransf12whiskright}}
        \POS"B0";"D3"**{}?*\txt{\eqref{eq:pstransf3cell}}
        \POS"B1";{"D3";"B3"**{}?}**{}?*\txt{\eqref{eq:pstransf2cellcomp}}
      \end{xy}
    \end{equation}
    \begin{center}
      Verfification of \eqref{eq:hcompbetaalpha2cellcomp}. Unlabeled
      subdiagrams commute by naturality. All 3-cells shown are
      canonical and invertible.
    \end{center}
  \end{sidewaysfigure}
  \qed
\end{prf}

\begin{lem}
  \label{lem:hcompbetaA}
  Our definition of $\_*_{-1}\_$ makes $\beta*_{-1}A$, a perturbation.
\end{lem}
\begin{prf}
  \revise{complete the proof}
\end{prf}

\begin{lem}
  \label{lem:hcompBalpha}
  Our definition of $\_*_{-1}\_$ makes $B*_{-1}\alpha$, a perturbation.
\end{lem}
\begin{prf}
  \revise{complete the proof}
\end{prf}

\revise{the derived horizontal compositions $\rhc$, $\lhc$…} 

\begin{thm}
  \label{thm:pasteunit}
  The $*_{-1}$-composition preserves units, i.e.\ $\beta*_{-1}\id_G=\id_{\beta*_{-1}G}$.
\end{thm}
\begin{prf}
  We see that
  $(\beta*_{-1}\id_G)_x=\beta_{\id_Gx}=(\id_{\beta*_{-1}G})_x$ by
  \eqref{eq:pstransf0id} and
  $(\beta*_{-1}\id_G)_f=\beta_{\id_{Gf}}=(\id_{\beta*_{-1}G})_f$ using
  \eqref{eq:betaalphaf} and \eqref{eq:pstransf2cocnorm}. \qed
\end{prf}

\begin{thm}
  \label{thm:betabetaalphapaste}
  For $\Gray$-functors and pseudo-transformations 
  \begin{equation*}
    \begin{xy}
      \xyboxmatrix{
        \G\ar@/^1.5pc/[r]^{G}="x"\ar@/_1.5pc/[r]_{G'}="y"\tara{\alpha}"x";"y"&\H\ar@/^2pc/[r]^{H}="x"\ar[r]|{H'}="y"\ar@/_2pc/[r]_{H''}="z"\tara{\beta}"x";"y"\tara{\beta'}"y";"z"&\K
      }
    \end{xy}
  \end{equation*}
we have the following compatibility of $*_{-1}$ and $*_0$:
  \begin{equation}
    \label{eq:betabetaalphapaste}
    \begin{xy}
      \xyboxmatrix@+.5cm@ur{
        HG\ar[r]^{\beta*_{-1}G}\ar[d]_{H*_{-1}\alpha}&H'G\ar[r]^{\beta'*_{-1}G}\ar[d]|{H'*_{-1}\alpha}\ar@<-1ex>@2[dl]**{}?/-1ex/;?/1ex/^{\beta*_{-1}\alpha}&H''G\ar[d]^{H''*_{-1}\alpha}\ar@<-1ex>@2[dl]**{}?/-1ex/;?/1ex/^{\beta'*_{-1}\alpha}\\
        HG'\ar[r]_{\beta*_{-1}G'}&H'G'\ar[r]_{\beta'*_{-1}G'}&H''G'
      }
    \end{xy}=
    \begin{xy}
      \xyboxmatrix@+.5cm@ur{
        HG\ar[r]^{\beta*_{-1}G}\ar[d]_{H*_{-1}\alpha}&H'G\ar[r]^{\beta'*_{-1}G}&H''G\ar[d]^{H''*_{-1}\alpha}\ar@<-1ex>@2[dll]**{}?/-1ex/;?/1ex/^(-1){(\beta'*_0\beta)*_{-1}\alpha}\\
        HG'\ar[r]_{\beta*_{-1}G'}&H'G'\ar[r]_{\beta'*_{-1}G'}&H''G'
      }
    \end{xy}\,.
  \end{equation}
\end{thm}
\begin{prf} 
  We need to show that the two pseudo-modifications defined
  in \eqref{eq:betabetaalphapaste} are equal.

  The left-hand and the right-hand side of
  \eqref{eq:betabetaalphapaste}can be evaluated according the
  definitions of $*_0$ and $*_{-1}$ given in this section and in
  \ref{sec:comp-pseu-transf}. On 0-cells we get
  \begin{multline*}
    (((\beta'*_{-1}G')*_0(\beta*_{-1}\alpha))*_1((\beta'*_{-1}\alpha)*_0(\beta*_{-1}G)))_x\\
    =(\beta'_{G'x}\#_0\beta_{\alpha_x})\#_1(\beta'_{\alpha_x}\#_0\beta_{Gx})
    =(\beta'*_0\beta)_{\alpha_x}=((\beta'*_0\beta)*_{-1}\alpha)_x.
  \end{multline*}

  On 1-cells the left hand side is given by
  \eqref{eq:betaalphabetapastef}; the right-hand side of
  \eqref{eq:betabetaalphapaste} is given by
  \eqref{eq:betabetaalphapastef}. We verify their equality in
    \eqref{eq:betabetaalphapastefdetail}.
  \begin{sidewaysfigure}
    \begin{equation}
      \label{eq:betaalphabetapastef}
      \def\outerbase{%
        \POS 0;<5cm,0cm>:(0,-1)::0%
      }
      \def\internalbase{\POS 0;<.13cm,0cm>:0,}
      \def\subdiagram#1#2{%
        \POS*!C++\xybox{%
          \internalbase%
          \vertices%
          \outerframe%
          #2%
        }="#1",%
      }
      \tbdef{(10,5)}{(8,-5)}{(0,-10)}
      \def\drawvertex#1{%
        \POS*{}="#1"%
      }
      \def\drawvertexat#1#2{%
        \t(#2)%
        \drawvertex{#1}%
      }
      \def\vertices{%
        \drawvertexat{a}{0,0,0}\drawvertexat{b}{1,0,0}\drawvertexat{c}{2,0,0}
        \drawvertexat{d}{0,1,0}\drawvertexat{e}{1,1,0}\drawvertexat{f}{2,1,0}
        \drawvertexat{a'}{0,0,1}\drawvertexat{b'}{1,0,1}\drawvertexat{c'}{2,0,1}
        \drawvertexat{d'}{0,1,1}\drawvertexat{e'}{1,1,1}\drawvertexat{f'}{2,1,1}
      }
      \def\outerframe{%
        \ar@{-}"a";"a'"
        \ar@{-}"a";"b"
        \ar@{-}"a'";"d'"
        \ar@{-}"b";"c"
        \ar@{-}"c";"f"
        \ar@{-}"d'";"e'"
        \ar@{-}"e'";"f'"
        \ar@{-}"f";"f'"
      }
      \begin{xy}
        \outerbase
        \subdiagram{A0}{
          \ar@{-}"a";"d"
          \ar@{-}"d";"e"
          \ar@{-}"e";"f"
          \pla{(\beta'*_0\beta)_{\alpha_x}}"c";"d"
          \pla{((\beta'*_0\beta)\rhc_{-1}\alpha)_f}"e";"d'"
        }+(2,0)
        \subdiagram{C0}{
          \ar@{-}"a'";"b'"
          \ar@{-}"b";"b'"
          \ar@{-}"b";"e"
          \ar@{-}"b'";"e'"
          \ar@{-}"e";"e'"
          \ar@{-}"e";"f"
          \pla{\beta_{Gf}}"b";"a'"
          \pla{\beta_{\alpha_y}}"b'";"d'"
          \pla{H'\alpha_f}"e";"b'"
          \pla{\beta'_{\alpha_x}}"c";"e"
          \pla{\beta'_{G'f}}"f";"e'"
        }+(2,0)
        \subdiagram{E0}{
          \ar@{-}"a'";"b'"
          \ar@{-}"b'";"c'"
          \ar@{-}"c'";"f'"
          \pla{((\beta'*_0\beta)\lhc_{-1}\alpha)_f}"c";"b'"
          \pla{(\beta'*_0\beta)_{\alpha_y}}"c'";"d'"
        },(0,1)
        \subdiagram{A1}{
          \ar@{-}"a";"d"
          \ar@{-}"b";"e"
          \ar@{-}"d";"d'"
          \ar@{-}"d";"e"
          \ar@{-}"e";"e'"
          \ar@{-}"e";"f"        
          \pla{\beta'_{\alpha_x}}"c";"e"
          \pla{\beta_{G'f}}"e";"d'"
          \pla{H\alpha_f}"d";"a'"
          \tria{\pla{\beta_{\alpha_x}}"b";"d"};{\pla{\beta'_{G'f}}"f";"e'"}
        }+(1,0)
        \subdiagram{B1}{
          \ar@{-}"a";"d"
          \ar@{-}"b";"e"
          \ar@{-}"d";"d'"
          \ar@{-}"d";"e"
          \ar@{-}"e";"e'"
          \ar@{-}"e";"f"        
          \pla{\beta'_{\alpha_x}}"c";"e"
          \pla{\beta_{G'f}}"e";"d'"
          \pla{H\alpha_f}"d";"a'"
          \tria{\pla{\beta'_{G'f}}"f";"e'"};{\pla{\beta_{\alpha_x}}"b";"d"}
        }+(1,0)
        \subdiagram{C1}{
          \ar@{-}"a'";"b'"
          \ar@{-}"b";"b'"
          \ar@{-}"b";"e"
          \ar@{-}"b'";"e'"
          \ar@{-}"e";"e'"
          \ar@{-}"e";"f"
          \pla{\beta_{Gf}}"b";"a'"
          \pla{\beta_{\alpha_y}}"b'";"d'"
          \pla{H'\alpha_f}"e";"b'"
          \pla{\beta'_{\alpha_x}}"c";"e"
          \pla{\beta'_{G'f}}"f";"e'"
        }+(1,0)
        \subdiagram{D1}{
          \ar@{-}"a'";"b'"
          \ar@{-}"b";"b'"
          \ar@{-}"b'";"c'"
          \ar@{-}"b'";"e'"
          \ar@{-}"c";"c'"
          \ar@{-}"c'";"f'"
          \pla{\beta'_{Gf}}"c";"b'"
          \pla{\beta_{\alpha_y}}"b'";"d'"
          \pla{H''\alpha_f}"c";"f'"
          \tria{\pla{\beta'_{\alpha_y}}"c'";"e'"};{\pla{\beta_{Gf}}"b";"a'"}
        }+(1,0)
        \subdiagram{E1}{          
          \ar@{-}"a'";"b'"   
          \ar@{-}"b";"b'"
          \ar@{-}"b'";"c'"
          \ar@{-}"b'";"e'"
          \ar@{-}"c";"c'"
          \ar@{-}"c'";"f'"
          \pla{\beta'_{Gf}}"c";"b'"
          \pla{\beta_{\alpha_y}}"b'";"d'"
          \pla{H''\alpha_f}"c";"f'"
          \tria{\pla{\beta_{Gf}}"b";"a'"};{\pla{\beta'_{\alpha_y}}"c'";"e'"}
        }+(1,0)
        \areq"A0";"A1"
        \ttara{
          \underline{((\beta'*_{-1}G)*_0(\beta*_{-1}\alpha))_f}\\
          \#_1(H''G'f\#_0\beta'_{\alpha_x}\#_0\beta_{Gx})
        }"A0";"C0"
        \ttar"A1";"B1"\offlabel{
          (\beta'_{G'y}\#_0\beta_{G'y}\#_0H\alpha_f)\\
          \#_1(\beta'_{G'y}\#_0\beta_{G'f}\#_0H\alpha_x)\\
          \#_1(\underline{\beta'_{G'f}\ten\beta_{\alpha_x}})\\
          \#_1(H''G'f\#_0\beta'_{\alpha_x}\#_0\beta_{Gx})
        };p+<0cm,-3cm>
        \ttar"B1";"C1"\offlabel{
          (\beta'_{G'y}\#_0\underline{(\beta*_{-1}\alpha)_f})\\
          \#_1(\beta'_{G'f}\#_0H'\alpha_x\#_0\beta_{Gx})\\
          \#_1(H''G'f\#_0\beta'_{\alpha_x}\#_0\beta_{Gx})
        };p+<0cm,-3cm>
        \areq"C0";"C1"
        \ttara{
          (\beta'_{G'y}\#_0\beta_{\alpha_y}\#_0HG_f)\\
          \#_1\underline{((\beta'*_{-1}G)*_0(\beta*_{-1}\alpha))_f}
        }"C0";"E0"
        \ttar"C1";"D1"\offlabel{%
          (\beta'_{G'y}\#_0\beta_{\alpha_y}\#_0HGf)\\
          \#_1(\beta'_{G'y}\#_0H'\alpha_y\#_0\beta_{Gf})\\
          \#_1(\underline{(\beta'*_{-1}\alpha)_f}\#_0\beta_{Gx})
        };p+<0cm,-3cm>
        \ttar"D1";"E1"\offlabel{
          (\beta'_{G'y}\#_0\beta_{\alpha_y}\#_0HGf)\\
          \#_1(\underline{\beta'_{\alpha_y}\ten\beta_{Gf}})^{-1}\\
          \#_1(H''\alpha_y\#_0\beta'_{Gf}\#_0\beta_{Gx})\\
          \#_1(H''\alpha_f\#_0\beta'_{Gx}\#_0\beta_{Gx})
        };p+<0cm,-3cm>
        \areq"E0";"E1"
        \pla{\text{\eqref{eq:psmodwhiskrightf}}}"C0";"A1"
        \pla{\text{\eqref{eq:psmodwhiskleftf}}}"E0";"C1"
      \end{xy}
    \end{equation}
    \begin{center}
      The left-hand side of \eqref{eq:betabetaalphapaste}, i.e.\
      $(((\beta'*_{-1}G')*_0(\beta*_{-1}\alpha))*_{-1}((\beta'*_{-1}\alpha)*_0(\beta*_{-1}G)))_f$,
      according to \eqref{eq:psmodcomp1}, \eqref{eq:psmodwhiskrightf}
      and \eqref{eq:psmodwhiskrightf}.
    \end{center}
  \end{sidewaysfigure}

  \begin{sidewaysfigure}
    \begin{center}
      $((\beta'*_0\beta)*_{-1}\alpha)_f$. See
      \eqref{eq:betabetatwoHfalphax}, \eqref{eq:betabetatwoalphayHf},
      and \eqref{eq:betabetaalphaf} for a breakdown of the 3-cells.
    \end{center}
    \begin{equation}
      \label{eq:betabetaalphapastef}
      \begin{xy}
        \POS;(0,70)
        \xyboxmatrix"A"@ru@+.3cm{
          {}&HGx\ar[r]^{(\beta'*_0\beta)_{Gx}}\ar[d]|{H\alpha_x}\ar[dl]|{HGf}&H''Gx\ar[d]^{H''\alpha_x}\ar@2[dl]**{}?/-1ex/;?/1ex/^{(\beta'*_0\beta)_{\alpha_x}}\\
          HGy\ar[d]_{H\alpha_y}&HG'x\ar[r]|{(\beta'*_0\beta)_{G'x}}\ar[dl]|{HG'f}\ar@2[l]**{}?/-1ex/;?/+1ex/^{H\alpha_f}&H''G'x\ar[dl]|{H''G'f}\ar@2[dll]**{}?/-1ex/;?/+1ex/^{(\beta'*_0\beta)_{G'f}}\\
          HG'\ar[r]_{(\beta'*_0\beta)_{G'y}}&H''G'y&{}
        };(55,70)
        \xyboxmatrix"B"@ru@+.3cm{
          {}&HGx\ar[r]^{(\beta'*_0\beta)_{Gx}}\ar[d]|{H\alpha_x}\ar[dl]|{HGf}&H''Gx\ar[d]^{H''\alpha_x}\\
          HGy\ar[d]_{H\alpha_y}&HG'x\ar[dl]|{HG'f}\ar@2[l]**{}?/-1ex/;?/+1ex/^{H\alpha_f}&H''G'x\ar[dl]|{H''G'f}\ar@2[l]**{}?/-1ex/;?/+1ex/_{(\beta'*_0\beta)_{G'f\#_0\alpha_x}}\\
          HG'\ar[r]_{(\beta'*_0\beta)_{G'y}}&H''G'y&{}
        };(110,70)
        \xyboxmatrix"C"@ru@+.3cm{
          {}&HGx\ar[r]^{(\beta'*_0\beta)_{Gx}}\ar[dl]|{HGf}&H''Gx\ar[d]^{H''\alpha_x}\ar[dl]|{H''Gf}\\
          HGy\ar[d]_{H\alpha_y}&H''Gx\ar[d]|{H''\alpha_y}\ar@2[l]**{}?/-1ex/;?/1ex/^{(\beta'*_0\beta)_{\alpha_y\#_0Gf}}&H''G'x\ar[dl]|{H''G'f}\ar@2[l]**{}?/-1ex/;?/+1ex/^{H''\alpha_f}\\
          HG'\ar[r]_{(\beta'*_0\beta)_{G'y}}&H''G'y&{}
        };(165,70) 
        \xyboxmatrix"D"@ru@+.3cm{
          {}&HGx\ar[r]^{(\beta'*_0\beta)_{Gx}}\ar[dl]|{HGf}&H''Gx\ar[d]^{H''\alpha_x}\ar[dl]|{H''Gf}\ar@2[dll]**{}?/-1ex/;?/1ex/^{(\beta'*_0\beta)_{Gf}}\\
          HGy\ar[d]_{H\alpha_y}\ar[r]|{(\beta'*_0\beta)_{Gy}}&H''Gx\ar[d]|{H''\alpha_y}\ar@2[dl]**{}?/-1ex/;?/1ex/^{(\beta'*_0\beta)_{\alpha_y}}&H''G'x\ar[dl]|{H''G'f}\ar@2[l]**{}?/-1ex/;?/+1ex/^{H''\alpha_f}\\
          HG'\ar[r]_{(\beta'*_0\beta)_{G'y}}&H''G'y&{} }
        \ar@3"A";"B"**{}?/-1ex/;?/1ex/|{}="here"\ar@{.}"here";p+(0,35)*!C\labelbox{((\beta'*_0\beta)_{G'_y}\#_0H\alpha_f)\\\#_1\underline{(\beta'*_0\beta)^2_{G'f,\alpha_x}}}
        \ar@3"B";"C"**{}?/-1ex/;?/1ex/|{}="here"\ar@{.}"here";p+(0,35)*!C\labelbox{\underline{(\beta'*_0\beta)_{\alpha_f}}}
        \ar@3"C";"D"**{}?/-1ex/;?/1ex/|{}="here"\ar@{.}"here";p+(0,35)*!C\labelbox{\underline{((\beta'*_0\beta)^2_{\alpha_y\#_0Gf})^{-1}}\\\#_1(H\alpha'_f,(\beta'*_0\beta)_{Gx})}        
      \end{xy}
    \end{equation}
  \end{sidewaysfigure}
  \begin{sidewaysfigure}
    \begin{center}Details of $(\beta'*_0\beta)^2_{G'f,\alpha_x}$\end{center}  
    \begin{equation}
      \label{eq:betabetatwoHfalphax}
      \begin{xy}
        \xyboxmatrix"A"@+.3cm@ur{
          {}&HGx\ar[r]^{\beta_{Gx}}\ar[d]_{H\alpha_x}&H'Gx\ar[r]^{\beta'_{Gx}}\ar[d]|{H'\alpha_x}\ar@2[dl]**{}?/-1ex/;?/1ex/_/+1.5ex/{\beta_{\alpha_x}}|{}="x"&H''Gx\ar[d]^{H''\alpha_x}\ar@2[dl]**{}?/-1ex/;?/1ex/^{\beta'_{\alpha_x}}\\
          {}&HG'x\ar[r]|{\beta_{G'x}}\ar[dl]_{HG'f}&H'G'x\ar[r]|{\beta'_{G'x}}\ar[dl]|{H'G'f}\ar@2@<-1ex>[dll]**{}?/-1ex/;?/1ex/^{\beta_{G'f}}&H''G'x\ar[dl]^{H''G'f}\ar@2@<-1ex>[dll]**{}?/-1ex/;?/1ex/^{\beta'_{G'f}}|{}="y"\\
          HG'y\ar[r]_{\beta_{G'y}}&H'G'y\ar[r]_{\beta'_{G'y}}&H''G'y&{}
          \tria"x";"y"
        };(75,0)
        \xyboxmatrix"B"@+.3cm@ur{
          {}&HGx\ar[r]^{\beta_{Gx}}\ar[d]_{H\alpha_x}&H'Gx\ar[r]^{\beta'_{Gx}}\ar[d]|{H'\alpha_x}\ar@2[dl]**{}?/-1ex/;?/1ex/_/+1.5ex/{\beta_{\alpha_x}}|{}="x"&H''Gx\ar[d]^{H''\alpha_x}\ar@2[dl]**{}?/-1ex/;?/1ex/^{\beta'_{\alpha_x}}\\
          {}&HG'x\ar[r]|{\beta_{G'x}}\ar[dl]_{HG'f}&H'G'x\ar[r]|{\beta'_{G'x}}\ar[dl]|{H'G'f}\ar@2@<-1ex>[dll]**{}?/-1ex/;?/1ex/^{\beta_{G'f}}&H''G'x\ar[dl]^{H''G'f}\ar@2@<-1ex>[dll]**{}?/-1ex/;?/1ex/^{\beta'_{G'f}}|{}="y"\\
          HG'y\ar[r]_{\beta_{G'y}}&H'G'y\ar[r]_{\beta'_{G'y}}&H''G'y&{}
          \tria"y";"x"
        };(150,0) 
        \xyboxmatrix"C"@+.3cm@ur{
          {}&HGx\ar[r]^{\beta_{Gx}}\ar[d]_{H\alpha_x}&H'Gx\ar[r]^{\beta'_{Gx}}\ar[d]|{H'\alpha_x}&H''Gx\ar[d]^{H''\alpha_x}\\
          {}&HG'x\ar[dl]_{HG'f}&H'G'x\ar[dl]|{H'G'f}\ar@2@<-1ex>[l]**{}?/-1ex/;?/1ex/_{\beta_{G'f\#_0\alpha_x}}&H''G'x\ar[dl]^{H''G'f}\ar@2@<-1ex>[l]**{}?/-1ex/;?/1ex/_{\beta'_{G'f\#_0\alpha_x}}\\
          HG'y\ar[r]_{\beta_{G'y}}&H'G'y\ar[r]_{\beta'_{G'y}}&H''G'y&{}
        }
        \ar@3"A";"B"**{}?/-1ex/;?/1ex/|{}="here"\ar@{.}"here";p+(0,40)*!C\labelbox{(\beta'_{G'y}\#_0\beta_{G'f}\#_0H\alpha_x)\\\#_1(\underline{\beta'_{G'f}\ten\beta_{\alpha_x}})\\\#_1(H''G'f\#_0\beta'_{\alpha_x}\#_0\beta_{Gx})}
        \ar@3"B";"C"**{}?/-1ex/;?/1ex/|{}="here"\ar@{.}"here";p+(0,40)*!C\labelbox{(\beta'_{G'y}\#_0\underline{\beta^2_{G'f,\alpha_x}})\\\#_1(\underline{\beta'^2_{G'f,\alpha_x}}\#_0\beta_{Gx})}
      \end{xy}
    \end{equation}
  \end{sidewaysfigure}
  \begin{sidewaysfigure}
    \begin{center}Details of $((\beta'*_0\beta)^2_{\alpha_y\#_0Gf})^{-1}$\end{center}  
    \begin{equation}
      \label{eq:betabetatwoalphayHf}
      \begin{xy}
        \xyboxmatrix"A"@+.1cm@ur{
          {}&HGx\ar[r]^{\beta_{Gx}}\ar[dl]_{HGf}&H'Gx\ar[r]^{\beta'_{Gx}}\ar[dl]|{H'Gf}&H''Gx\ar[dl]^{H''Gf}\\
          HGy\ar[d]_{H\alpha_y}&H'Gy\ar[d]|{H'\alpha_y}\ar@2[l]**{}?/-1ex/;?/1ex/^{\beta_{\alpha_y\#_0Gf}}&H''Gy\ar[d]^{H''\alpha_y}\ar@2[l]**{}?/-1ex/;?/1ex/^{\beta'_{\alpha_y\#_0Gf}}&{}\\
          HG'y\ar[r]_{\beta_{G'y}}&H'G'y\ar[r]_{\beta'_{G'y}}&H''G'y&{}
        };(75,0)
        \xyboxmatrix"B"@+.1cm@ur{
          {}&HGx\ar[r]^{\beta_{Gx}}\ar[dl]_{HGf}&H'Gx\ar[r]^{\beta'_{Gx}}\ar[dl]|{H'Gf}\ar@2@<1ex>[dll]**{}?/-1ex/;?/1ex/_{\beta_{Gf}}|{}="x"&H''Gx\ar[dl]^{H''Gf}\ar@2@<-1ex>[dll]**{}?/-1ex/;?/1ex/^{\beta'_{Gf}}\\
          HGy\ar[d]_{H\alpha_y}\ar[r]|{\beta_{Gy}}&H'Gy\ar[r]|{\beta'_{Gy}}\ar[d]|{H'\alpha_y}\ar@2[dl]**{}?/-1ex/;?/1ex/^{\beta_{\alpha_x}}&H''Gy\ar[d]^{H''\alpha_y}\ar@2[dl]**{}?/-1ex/;?/1ex/^{\beta'_{\alpha_y}}|{}="y"&{}\\
          HG'y\ar[r]_{\beta_{G'y}}&H'G'y\ar[r]_{\beta'_{G'y}}&H''G'y&{}
          \tria"y";"x"
        };(150,0)
        \xyboxmatrix"C"@+.1cm@ur{
          {}&HGx\ar[r]^{\beta_{Gx}}\ar[dl]_{HGf}&H'Gx\ar[r]^{\beta'_{Gx}}\ar[dl]|{H'Gf}\ar@2@<1ex>[dll]**{}?/-1ex/;?/1ex/_{\beta_{Gf}}|{}="x"&H''Gx\ar[dl]^{H''Gf}\ar@2@<-1ex>[dll]**{}?/-1ex/;?/1ex/^{\beta'_{Gf}}\\
          HGy\ar[d]_{H\alpha_y}\ar[r]|{\beta_{Gy}}&H'Gy\ar[r]|{\beta'_{Gy}}\ar[d]|{H'\alpha_y}\ar@2[dl]**{}?/-1ex/;?/1ex/^{\beta_{\alpha_x}}&H''Gy\ar[d]^{H''\alpha_y}\ar@2[dl]**{}?/-1ex/;?/1ex/^{\beta'_{\alpha_y}}|{}="y"&{}\\
          HG'y\ar[r]_{\beta_{G'y}}&H'G'y\ar[r]_{\beta'_{G'y}}&H''G'y&{}
          \tria"x";"y"
        }
        \ar@3"A";"B"**{}?/-1ex/;?/1ex/|{}="here"\ar@{.}"here";p+(0,40)*!C\labelbox{(\beta'_{G'f}\#_0(\underline{\beta^2_{\alpha_y\#_0Gf}})^{-1})\\\#_1((\underline{\beta'^2_{\alpha_y\#_0Gf}})^{-1}\#_0\beta_{Gx})} 
        \ar@3"B";"C"**{}?/-1ex/;?/1ex/|{}="here"\ar@{.}"here";p+(0,40)*!C\labelbox{(\beta'_{G'y}\#_0\beta_{\alpha_x}\#_0HGf)\\\#_1(\underline{\beta'_{\alpha_y}\ten\beta_{\beta_{Gf}}})^{-1}\\\#_1(H''\alpha_y\#_0\beta'_{Gf}\#_0\beta_{Gx})}
      \end{xy}
    \end{equation}
  \end{sidewaysfigure}
  \begin{sidewaysfigure}
    \begin{center}Details of $(\beta'*_0\beta)_{\alpha_f}$\end{center}  
    \begin{equation}
      \label{eq:betabetaalphaf}
      \begin{xy}
        \xyboxmatrix"A"@+.1cm@ur{
          {}&HGx\ar[r]^{\beta_{Gx}}\ar[d]|{H\alpha_x}\ar[dl]_{HGf}&H'Gx\ar[r]^{\beta'_{Gx}}\ar[d]|{H'\alpha_x}&H''Gx\ar[d]^{H''\alpha_x}\\
          HGy\ar[d]_{H\alpha_y}&HG'x\ar[dl]|{HG'f}\ar@2@<-1ex>[l]**{}?/-1ex/;?/1ex/^{H\alpha_f}&H'G'x\ar[dl]|{H'G'f}\ar@2@<-1ex>[l]**{}?/-1ex/;?/1ex/^{\beta_{G'f\#_0\alpha_x}}&H''G'x\ar[dl]^{H''G'f}\ar@2@<-1ex>[l]**{}?/-1ex/;?/1ex/^{\beta'_{G'f\#_0\alpha_x}}\\
          HG'y\ar[r]_{\beta_{G'y}}&H'G'y\ar[r]_{\beta'_{G'y}}&H''G'y&{}
        };(85,0)
        \xyboxmatrix"B"@+.1cm@ur{
          {}&HGx\ar[r]^{\beta_{Gx}}\ar[dl]_{HGf}&H'Gx\ar[r]^{\beta'_{Gx}}\ar[d]|{H'\alpha_x}\ar[dl]|{H'Gf}&H''Gx\ar[d]^{H''\alpha_x}\\
          HGy\ar[d]_{H\alpha_y}&H'Gy\ar[d]|{H'\alpha_y}\ar@2[l]**{}?/-1ex/;?/1ex/^{\beta_{\alpha_y\#_0Gf}}&H'G'x\ar[dl]|{H'G'f}\ar@2@<-1ex>[l]**{}?/-1ex/;?/1ex/^{H'\alpha_f}&H''G'x\ar[dl]^{H''G'f}\ar@2@<-1ex>[l]**{}?/-1ex/;?/1ex/^{\beta'_{G'f\#_0\alpha_x}}\\
          HG'y\ar[r]_{\beta_{G'y}}&H'G'y\ar[r]_{\beta'_{G'y}}&H''G'y&{}
        };(170,0) 
        \xyboxmatrix"C"@+.1cm@ur{
          {}&HGx\ar[r]^{\beta_{Gx}}\ar[dl]_{HGf}&H'Gx\ar[r]^{\beta'_{Gx}}\ar[dl]|{H'Gf}&H''Gx\ar[d]^{H''\alpha_x}\ar[dl]|{H''Gf}\\
          HGy\ar[d]_{H\alpha_y}&H'Gy\ar[d]|{H'\alpha_y}\ar@2[l]**{}?/-1ex/;?/1ex/^{\beta_{\alpha_y\#_0Gf}}&H''Gy\ar[d]|{H''\alpha_y}\ar@2[l]**{}?/-1ex/;?/1ex/^{\beta'_{\alpha_y\#_0Gf}}|{}="y"&H''G'x\ar[dl]^{H''G'f}\ar@2@<-1ex>[l]**{}?/-1ex/;?/1ex/^{H''\alpha_f}\\
          HG'y\ar[r]_{\beta_{G'y}}&H'G'y\ar[r]_{\beta'_{G'y}}&H''G'y&{}
        }
        \ar@3"A";"B"**{}?/-1ex/;?/1ex/^{(\beta'_{G'f}\#_0\underline{\beta_{\alpha_f}})\\\#_1(\beta'_{G'f\#_0\alpha_x}\#_0\beta_{Gx})} 
        \ar@3"B";"C"**{}?/-1ex/;?/1ex/^{(\beta'_{G'y}\#_0\beta_{\alpha_y\#_0Gf})\\\#_0(\underline{\beta'_{\alpha_f}}\#_0\beta_{Gx})}
      \end{xy}
    \end{equation}
  \end{sidewaysfigure}

  \begin{sidewaysfigure}
    \begin{equation}
      \label{eq:betabetaalphapastefdetail}
      \def\outerbase{%
        \POS 0;<5cm,0cm>:(0,-1)::0%
      }
      \def\internalbase{\POS 0;<.13cm,0cm>:0,}
      \def\subdiagram#1#2{%
        \POS*!C++\xybox{%
          \internalbase%
          \vertices%
          \outerframe%
          #2%
        }="#1",%
      }
      \tbdef{(10,5)}{(8,-5)}{(0,-10)}
      \def\drawvertex#1{%
        \POS*{}="#1"%
      }
      \def\drawvertexat#1#2{%
        \t(#2)%
        \drawvertex{#1}%
      }
      \def\vertices{%
        \drawvertexat{a}{0,0,0}\drawvertexat{b}{1,0,0}\drawvertexat{c}{2,0,0}
        \drawvertexat{d}{0,1,0}\drawvertexat{e}{1,1,0}\drawvertexat{f}{2,1,0}
        \drawvertexat{a'}{0,0,1}\drawvertexat{b'}{1,0,1}\drawvertexat{c'}{2,0,1}
        \drawvertexat{d'}{0,1,1}\drawvertexat{e'}{1,1,1}\drawvertexat{f'}{2,1,1}
      }
      \def\outerframe{%
        \ar@{-}"a";"a'"
        \ar@{-}"a";"b"
        \ar@{-}"a'";"d'"
        \ar@{-}"b";"c"
        \ar@{-}"c";"f"
        \ar@{-}"d'";"e'"
        \ar@{-}"e'";"f'"
        \ar@{-}"f";"f'"
      }
      \begin{xy}
        \outerbase
        ,(-2,0)
        \subdiagram{A0}{
          \ar@{-}"a";"d"
          \ar@{-}"b";"e"
          \ar@{-}"d";"d'"
          \ar@{-}"d";"e"
          \ar@{-}"e";"e'"
          \ar@{-}"e";"f"        
          \pla{\beta'_{\alpha_x}}"c";"e"
          \pla{\beta_{G'f}}"e";"d'"
          \pla{H\alpha_f}"d";"a'"
          \tria{\pla{\beta_{\alpha_x}}"b";"d"};{\pla{\beta'_{G'f}}"f";"e'"}
        }+(1,0)
        \subdiagram{B0}{
          \ar@{-}"a";"d"
          \ar@{-}"b";"e"
          \ar@{-}"d";"d'"
          \ar@{-}"d";"e"
          \ar@{-}"e";"e'"
          \ar@{-}"e";"f"        
          \pla{\beta'_{\alpha_x}}"c";"e"
          \pla{\beta_{G'f}}"e";"d'"
          \pla{H\alpha_f}"d";"a'"
          \tria{\pla{\beta'_{G'f}}"f";"e'"};{\pla{\beta_{\alpha_x}}"b";"d"}
        }+(1,0)
        \subdiagram{C0}{
          \ar@{-}"a'";"b'"
          \ar@{-}"b";"b'"
          \ar@{-}"b";"e"
          \ar@{-}"b'";"e'"
          \ar@{-}"e";"e'"
          \ar@{-}"e";"f"
          \pla{\beta_{Gf}}"b";"a'"
          \pla{\beta_{\alpha_y}}"b'";"d'"
          \pla{H'\alpha_f}"e";"b'"
          \pla{\beta'_{\alpha_x}}"c";"e"
          \pla{\beta'_{G'f}}"f";"e'"
        }+(1,0)
        \subdiagram{D0}{
          \ar@{-}"a'";"b'"
          \ar@{-}"b";"b'"
          \ar@{-}"b'";"c'"
          \ar@{-}"b'";"e'"
          \ar@{-}"c";"c'"
          \ar@{-}"c'";"f'"
          \pla{\beta'_{Gf}}"c";"b'"
          \pla{\beta_{\alpha_y}}"b'";"d'"
          \pla{H''\alpha_f}"c";"f'"
          \tria{\pla{\beta'_{\alpha_y}}"c'";"e'"};{\pla{\beta_{Gf}}"b";"a'"}
        }+(1,0)
        \subdiagram{E0}{          
          \ar@{-}"a'";"b'"   
          \ar@{-}"b";"b'"
          \ar@{-}"b'";"c'"
          \ar@{-}"b'";"e'"
          \ar@{-}"c";"c'"
          \ar@{-}"c'";"f'"
          \pla{\beta'_{Gf}}"c";"b'"
          \pla{\beta_{\alpha_y}}"b'";"d'"
          \pla{H''\alpha_f}"c";"f'"
          \tria{\pla{\beta_{Gf}}"b";"a'"};{\pla{\beta'_{\alpha_y}}"c'";"e'"}
        },(-1.5,1)
        \subdiagram{A1}{
          \ar@{-}"a";"d"
          \ar@{-}"b";"e"
          \ar@{-}"d";"d'"
          \ar@{-}"e";"e'"
          \ar@{-}"e";"f"
          \pla{H\alpha_f}"d";"a'"
          \pla{\beta_{G'f\#_0\alpha_x}}"e";"d"
          \pla{\beta'_{\alpha_x}}"c";"e"
          \pla{\beta'_{G'f}}"f";"e'"
        }+(1,0)
        \subdiagram{B1}{
          \ar@{-}"b";"b'"
          \ar@{-}"b";"e"
          \ar@{-}"b'";"e'"
          \ar@{-}"e";"e'"
          \ar@{-}"e";"f"
          \pla{\beta_{\alpha_y\#_0Gf}}"b'";"a'"
          \pla{H'\alpha_f}"e";"b'"
          \pla{\beta'_{\alpha_x}}"c";"e"
          \pla{\beta'_{G'f}}"f";"e'"
        }+(1,0)
        \subdiagram{C1}{
          \ar@{-}"a'";"b'"
          \ar@{-}"b";"b'"
          \ar@{-}"b";"e"
          \ar@{-}"b'";"e'"
          \ar@{-}"e";"e'"
          \pla{\beta_{Gf}}"b";"a'"
          \pla{\beta_{\alpha_y}}"b'";"d'"
          \pla{H'\alpha_f}"e";"b'"
          \pla{\beta'_{G'f\#_0\alpha_x}}"f";"e"
        }+(1,0)
        \subdiagram{D1}{
          \ar@{-}"a'";"b'"
          \ar@{-}"b";"b'"
          \ar@{-}"b'";"e'"
          \ar@{-}"c";"c'"
          \ar@{-}"c'";"f'"
          \pla{\beta_{Gf}}"b";"a'"
          \pla{\beta_{\alpha_y}}"b'";"d'"
          \pla{H''\alpha_f}"f";"c'"
          \pla{\beta'_{\alpha_y\#_0Gf}}"c'";"b'"
        }
        ,(-2,2)
        \subdiagram{A2}{
          \ar@{-}"a";"d"
          \ar@{-}"b";"e"
          \ar@{-}"d";"d'"
          \ar@{-}"e";"e'"
          \pla{H\alpha_f}"d";"a'"
          \pla{\beta_{G'f\#_0\alpha_x}}"e";"d"
          \pla{\beta'_{G'f\#_0\alpha_x}}"f";"e"
        }+(2,0)
        \subdiagram{B2}{
          \ar@{-}"b";"b'"
          \ar@{-}"b";"e"
          \ar@{-}"b'";"e'"
          \ar@{-}"e";"e'"
          \pla{\beta'_{G'f\#_0\alpha_x}}"f";"e"
          \pla{\beta_{\alpha_y\#_0Gf}}"b'";"a'"
          \pla{H'\alpha_f}"e";"b'"
        }+(2,0)
        \subdiagram{C2}{
          \ar@{-}"b";"b'"
          \ar@{-}"b'";"e'"
          \ar@{-}"c";"c'"
          \ar@{-}"c'";"f'"
          \pla{\beta_{\alpha_y\#_0Gf}}"b'";"a'"
          \pla{H''\alpha_f}"f";"c'"
          \pla{\beta'_{\alpha_y\#_0Gf}}"c'";"b'"
        }
        \ttar"A0";"B0"
        \ttar"B0";"C0"
        \ttar"C0";"D0"
        \ttar"D0";"E0"
        \ttar"A1";"B1"
        \ttar"A2";"B2"
        \ttar"A2";"A1"        
        \ttar"B0";"A1"
        \ar@3@/_6pc/"B0";"A2"
        \ttar"B1";"C0"
        \ttar"B2";"C2"
        \ttar"B2";"B1"
        \ttar"B2";"C1"
        \ttar"C0";"C1"
        \ttar"C1";"D1"
        \ttar"C2";"D1"
        \ttar"D1";"D0"
        \ar@3@/_6pc/"C2";"D0"
        \pla{\text{\eqref{eq:betaalphaf}}}"A1";"C0"
        \pla{\text{\eqref{eq:betaalphaf}}}"C1";"D0"
      \end{xy}
    \end{equation}
    \begin{center}
      Verification of the equality of \eqref{eq:betaalphabetapastef}
      (top row) and \eqref{eq:betabetaalphapastef} (remaining
      outline). Unlabeled subdiagrams commute by naturality
    \end{center}
  \end{sidewaysfigure}
  \qed
\end{prf}

\clearpage
\label{sec:bibl} \bibliographystyle{plainnat} \bibliography{skewgray}

\end{document}